\documentclass[a4paper,11pt]{article}

\usepackage[utf8]{inputenc} 
\usepackage[T1]{fontenc}
\usepackage[english]{babel}
\usepackage{amsmath, amsfonts, amsthm,amssymb,here,dsfont, hyperref}
\usepackage{mathtools}
\usepackage{graphicx,color,subfigure,multirow, here}
\usepackage{natbib}
\usepackage{enumitem}
\usepackage{ulem}
\usepackage{comment}

\newtheorem{theo}{Theorem}[section]

\newtheorem{prop}[theo]{Proposition}

\newtheorem{lemma}[theo]{Lemma}
\newtheorem{remark}[theo]{Remark}

\newtheorem{assumption}[theo]{Assumption}

\newcommand{\w}{\widehat}

\newcommand{\one}{\mathds{1}}

\newcommand{\E}{\mathbb{E}}

\renewcommand{\P}{\mathbb{P}}

\newcommand{\R}{\mathbb{R}}

\title{Minimax rates of convergence for the nonparametric estimation of the diffusion coefficient from time-homogeneous SDE paths}

\author{Eddy-Michel Ella-Mintsa$^{(1,2)}$}

\begin{document}

\maketitle
\begin{center}
(1) LAMA, Université Gustave Eiffel\\ 
(2) MIA Paris-Saclay, AgroParisTech
\end{center}

\begin{abstract}
Consider a diffusion process $X = (X_t)_{t \in [0,1]}$, solution of a time-homogeneous stochastic differential equation observed at discrete times and high frequency. We assume that the drift and diffusion coefficients of the process $X$ are unknown. In this paper, we study the minimax convergence rates of projection estimators of the square of the diffusion coefficient. Two observation schemes are considered depending on the estimation interval. The square of the diffusion coefficient is estimated on the real line from repeated observations of the process $X$, where the number of diffusion paths tends to infinity. For the case of a compact estimation interval, we study projection estimators of the square of the diffusion coefficient constructed from a single diffusion path on one side and from repeated observations on the other side, where the number of trajectories tends to infinity. In each of these cases, we establish minimax convergence rates of the worst risk of estimation over a space of H\"older functions. 
\end{abstract}



\section{Introduction}

Consider the following stochastic differential equation.
\begin{equation}\label{eq:model}
    dX_t = b(X_t)dt + \sigma(X_t)dW_t, ~~ t \in [0,T], ~~ X_0 = x_0,
\end{equation}
where the drift coefficient $b$ and the diffusion coefficient $\sigma$ are unknown, $(W_t)_{t \geq 0}$ is the standard Brownian motion and $T>0$ is the time horizon. The diffusion process $X = (X_t)_{t \in [0,T]}$, solution of Equation~\eqref{eq:model}, belongs to a probability space $(\Omega, \mathcal{F}, \mathbb{P}_X)$, where $\mathbb{P}_X$ is the distribution of the process $X$. We denote by $\left(\mathcal{F}_t\right)_{t \in [0,T]}$ the natural filtration of the diffusion process $X$. 

In this paper, we suppose to have a sample $\left\{X_{k\Delta_n}^j, ~ j = 1, \ldots, N, ~ k = 0, \ldots, n\right\}$ composed of $Nn$ observations of the process $X$ collected at discrete times with time step $\Delta_n = T/n$. We assume that the observations are collected at high frequency, that is, the time step $\Delta_n$ tends to zero ($n \rightarrow \infty$). The objective of this paper is to study the minimax convergence rate of the risk of estimation of the square of the diffusion coefficient $\sigma^2$ considering two simulation schemes. We focus on the estimation of $\sigma^2$ both on any compact interval and on the real line $\mathbb{R}$. Estimation on the real line is only considered in the setting where $N,n \rightarrow \infty$. For the case of the estimation of $\sigma^2$ on a compact interval, two statistical settings are studied. First, we investigate the minimax rate of the estimator of $\sigma^2$ from a single diffusion path ($N = 1$ and $n \rightarrow \infty$). Second, we extend the study to the case of $N$ independent discrete observations of the diffusion process $X$, with $N \rightarrow \infty$.

\subsection{Motivation}

Diffusion processes are increasingly used for the modeling of random phenomena in areas ranging from biology (see, e.g., \cite{sbalzarini2006analysis}, \cite{bressloff2024cellular}), physics (see, e.g., \cite{romanczuk2012active}), to quantitative finance (see, e.g., \cite{el1997backward}). More specifically, diffusion paths are at the core of the modeling of random phenomena on topics such as population genetics (see, e.g., \cite{crow2017introduction}), populations dynamics (see, e.g., \cite{nagai1983asymptotic}), evolution over time of the price of a financial asset (see, e.g., \cite{lamberton2011introduction}) or the interest rate (see, e.g., \cite{cairns2004interest}, \cite{wu2019interest}). 

Short-time diffusion processes $X = (X_t)_{t \in [0,T]}$, (with $T < \infty$) are used to study random phenomena that occur in a short time. We can mention, for example, the study of the price of financial assets such as European-style options, whose price can be modeled by the diffusion process $(S_t)_{t \in [0,T]}$, solution of the Black \& Scholes model which reads as follows:
\begin{equation*}
    dS_t = \mu S_t dt + \sigma S_t dW_t, ~~ t \in [0,T], ~ S_0 > 0, 
\end{equation*}
where $S_0>0$ is the actual price of the asset, meaning the price at $t=0$, and the time horizon $T>0$ is the time until the option expires (see \textit{e.g.} \cite{lamberton2011introduction}). The above diffusion model can be considered to study the risky or risk-free status of the financial asset based on the evolution of its volatility characterized by the diffusion coefficient $x \mapsto \sigma x$ since the quadratic variation of the process $(S_t)_{t \in [0,T]}$ satisfies
$$\dfrac{d\left<S,S\right>_t}{dt} = \sigma^2 S_t^2,$$
or to discriminate between high-price and low-price options considering the range of their respective prices at time $T<\infty$, which depends on the drift and diffusion coefficients. Solving each of these problems can lead to the construction of a classification procedure based on the estimation of the drift and diffusion coefficients (see \cite{denis2020consistent}, \cite{gadat2020optimal}, \cite{denis2024nonparametric}). More generally, if we go back to Equation~\eqref{eq:model}, then the coefficients $b$ and $\sigma$ satisfy
\begin{equation*}
    b(x) := \underset{h \rightarrow 0}{\lim}{h^{-1}\mathbb{E}\left(X_{t+h} - x | X_t = x\right)}, ~~ x \in \mathbb{R}, ~ t \in [0,T),
\end{equation*}
and conditional on the event $\{X_t = x\}$, with $x \in \mathbb{R}$ and $t \in [0,T)$,
\begin{equation*}
    \sigma^2(x) := \underset{h \rightarrow 0}{\lim}{\dfrac{\left<X,X\right>_{t+h} - \left<X,X\right>_t}{h}}.
\end{equation*}
From the above definitions of the two coefficients, it is impossible to build a consistent estimator of the drift function $b$, interpreted as the instantaneous mean of the process $X$, from a single discrete observation of the short-time diffusion process $X = (X_t)_{t \in [0, T]}$ since $n\Delta_n = T < \infty$. This condition justifies the use of repeated observations, that is, $N$ independent copies $X^1, \ldots, X^N$ of the diffusion process $X$, with $N \rightarrow \infty$, for the construction of a consistent estimator of $b$ (see \textit{e.g.} \cite{comte2020nonparametric}, \cite{denis2020ridge}, \cite{denis2024nonparametric}). In contrast, the square $\sigma^2$ of the diffusion coefficient can be estimated from a single copy of the short-time diffusion process $X$ (see \textit{e.g.} \cite{florens1993estimating}, \cite{hoffmann1999lp}, \cite{ella2024nonparametric}). More importantly, faster convergence rates can be reached when considering $N$ independent observations of $X$ with $N \rightarrow \infty$ (see \textit{e.g.} \cite{ella2024nonparametric}).

In the sequel, and for the sake of simplicity, the time horizon $T<\infty$ is set to $T = 1$.

\subsection{Related works}

There exists an important literature on statistical inference for stochastic differential equations. This statistical inference can be parametric, where the coefficients of the diffusion process depend on an unknown parameter to be estimated (see, e.g., \cite{jacod1993random}, \cite{genon1993estimation}, \cite{clement1997}, \cite{gloter2000discrete}, \cite{sorensen2002estimation}). For the nonparametric setting, $\sigma^2$ is estimated on a compact interval (see, e.g., \cite{hoffmann1997minimax}, \cite{hoffmann1999lp}, \cite{comte2007penalized}), or on the real line $\mathbb{R}$ (see \cite{florens1993estimating}, \cite{denis2024nonparametric}, \cite{ella2024nonparametric}). Focusing on the nonparametric estimation of the square $\sigma^2$ of the diffusion coefficient on a compact interval, \cite{comte2007penalized} proposed non-adaptive and adaptive estimators of $\sigma^2$ on the compact interval $[0,1]$. Not far from our framework, \cite{hoffmann1999lp} proposed a projection estimator of $\sigma^2$ on the compact interval $[0,1]$ from one observation of a short-time diffusion process and established a minimax convergence rate of order $n^{-s/(2s+1)}$ over a Besov space of smoothness parameter $s>3/2$, taking advantage of the optimality of the wavelet basis. \cite{ella2024nonparametric} extended the study to other bases such as the spline basis or the Fourier basis, proposed a projection estimator of $\sigma^2$ on a compact interval, from a single diffution path ($N=1$ and $n \rightarrow \infty$) and from repeated observations ($N,n \rightarrow \infty$). The author established convergence rates of order $n^{-\beta/(2\beta+1)}$ (for $N=1$) and $(Nn)^{-\beta/(2\beta+1)}$ (for $N \rightarrow \infty$) over the H\"older space of smoothness parameter $\beta > 3/2$. The first nonparametric estimator of $\sigma^2$ on the real line and from a single diffusion path is the consistent kernel estimator proposed in \cite{florens1993estimating}. \cite{denis2024nonparametric} followed with the construction of a projection estimator $\sigma^2$ on the real line and from $N$ independent copies of the short-time diffusion process $X$. They established a convergence rate of order $N^{-1/5}$ over the space of Lipschitz functions. \cite{ella2024nonparametric} obtained, when $N$ and $n$ tend to infinity, a faster convergence rate of order $(Nn)^{-\beta/(4\beta+1)}$ over the H\"older space of smoothness parameter $\beta \geq 1$, where the risk of estimation is defined from an empirical pseudo-norm, and a rate of order $N^{-\beta/2(\beta+1)}$ when the risk of estimation is defined from an empirical norm.

\subsection{Main contributions}

To our knowledge, the only paper that studied the minimax convergence rate for estimating the square of a space-dependent diffusion coefficient $\sigma^2$ is \cite{hoffmann1999lp}. The author proposed, from a single diffusion path observed at discrete times with time step $\Delta_n = 1/n$, a projection estimator of $\sigma^2$, restricted to the compact interval $[0,1]$, on a finite-dimensional space spanned by the wevelet basis. The optimality of the wavelet basis has been consequential in establishing the minimax rate. In \cite{ella2024nonparametric}, the study of the nonparametric estimation of $\sigma^2$ on a compact interval and from a single diffusion path is extended to bases of Lipschitz functions such as the spline basis or the Fourier basis, and a convergence rate of order $n^{-\beta/(2\beta+1)}$ of the $L^2-$risk of estimation has been established over the H\"older class of smoothness parameter $\beta > 3/2$. The author also considered the case of repeated observations $(N,n \rightarrow \infty)$ and established a convergence rate of the $L^2-$risk of estimation of order $(Nn)^{-\beta/(2\beta+1)}$.

The main contributions of this paper are the following.
\begin{enumerate}
    \item In the context of nonparametric estimation of $\sigma^2$ on a compact interval, we establish the optimality of the convergence rates $n^{-\beta/(2\beta+1)}$ (for $N=1$ and $n \rightarrow \infty$) and $(Nn)^{-\beta/(2\beta+1)}$ (for $N,n \rightarrow \infty$) over the H\"older space of smoothness parameter $\beta > 3/2$ . More specifically, in addition to the upper bounds of the $L^2-$risks of estimation of $\sigma^2$ obtained in \cite{ella2024nonparametric}, we establish, for each simulation scheme, a lower bound of the same order using Theorem 2.5 in \cite{tsybakov2008introduction}, and, under the assumptions on the coefficients $b$ and $\sigma$ presented in the next section, the exact formula of the transition density of the diffusion process $X$, solution of Equation~\eqref{eq:model}, provided in \cite{dacunha1986estimation}. These lower bounds will require the diffusion coefficients $b$ and $\sigma$ to be smooth enough, with $\beta > 2$ when $\sigma^2$ is estimated from a single diffusion path, and $\beta \geq 4$ for the estimation of $\sigma^2$ from repeated observations of the diffusion process $X$.
    \item Focusing on the nonparametric estimation of $\sigma^2$ on a growing interval $[-A_N,A_N]$ depending on the number $N$ of diffusion paths and such that $A_N > 0$ and the compact interval $[-A_N, A_N]$ tends to the real line $\mathbb{R}$ as $N$ tends to infinity, we establish an upper bound of the estimation risk of the projection estimator of $\sigma^2$ of order $\log^{\beta}(N)(Nn)^{-2\beta/(2\beta+1)}$ over the H\"older space of smoothness parameter $\beta$. This rate is obtained under specific assumptions on the Gram matrix of the basis that generates the approximation space as described in Section~\ref{subsec:Nonparametric_estimation_sigma}. Moreover, the use of the explicit formula of the transition density of the unique strong solution $X$ of Equation~\eqref{eq:model}, under regularity assumptions on $b$ and $\sigma$, allows us to derive the desired rate of convergence with $A_N \propto \sqrt{\log(N)}$. This condition on $A_N$ allows to establish a convergence rate of order $N^{-3\beta/2(2\beta+1)}$ (up to an extra factor or order $\exp\left(c\sqrt{\log(N)}\right)$) for the estimation risk of the nonparametric estimator of $\sigma^2$ on $\mathbb{R}$. Note that the establishment of these two upper bounds imposes the smoothness parameter $\beta$ to be large enough, as described in Section~\ref{subsec:Nonparametric_estimation_sigma}.
    \item We establish a lower bound of order $(Nn)^{-2\beta/(2\beta+1)}$ for the estimation risk of the projection estimator of $\sigma^2$ on the growing interval $[-A_N, A_N]$ and on the real line $\mathbb{R}$.   
\end{enumerate}

\subsection{Outline of the paper}

In Section~\ref{sec:framework and assumptions}, we present the statistical setting for the construction of nonparametric estimators of $\sigma^2$. We study the risk bounds of a projection estimator of $\sigma^2$ on the real line in Section~\ref{sec:upper-bound}, and the lower bounds are investigated in Section~\ref{sec:lower-bound}, considering the estimation of $\sigma^2$ on a compact interval and on the real line. Section~\ref{sec:conclusion} and 
\ref{sec:proofs} are respectively devoted to the conclusion and the proofs of out main results.

\section{Framework and assumptions}
\label{sec:framework and assumptions}

Consider a diffusion process $X=(X_t)_{t\in[0,1]}$, solution of Equation \eqref{eq:model}, and its discrete-time observation $\bar{X}=(X_{k\Delta_n})_{0\leq k\leq n}$ where the time step $\Delta_n = 1/n$. The objective is to construct a nonparametric estimator of $\sigma^{2}_{|I} = \sigma^2 \mathds{1}_{I}$, where $I \subset \mathbb{R}$ is either a compact interval or the real line. For this purpose, we make the following assumptions.

\begin{assumption}\label{ass:LipFunctions}
    There exists a constant $L_0 > 0$ such that for all $x \in \mathbb{R}$,
    \begin{align*}
        \left|b(x) - b(y)\right| + \left|\sigma(x) - \sigma(y)\right| \leq L_0 |x-y|.
    \end{align*}
\end{assumption}

\begin{assumption}\label{ass:Ellipticity}
    There exist constants $\kappa_0, \kappa_1 > 0$ such that $\kappa_0 < \kappa_1$ and
    $$\forall ~ x \in \mathbb{R}, ~~ 0 < \kappa_0 \leq \sigma(x) \leq \kappa_1.$$
\end{assumption}

\begin{assumption}
    \label{ass:RegularityBis}
    $\sigma\in\mathcal{C}^{2}\left(\mathbb{R}\right)$, and there exist constants $\gamma, C > 0$ such that
    \begin{equation*}
        \forall ~ x \in \mathbb{R}, ~~ |\sigma^{\prime}(x)| + |\sigma^{\prime\prime}(x)| \leq C(1+|x|^{\gamma}).
    \end{equation*}
\end{assumption}

\begin{assumption}
    \label{ass:Restrict-Model} 
    $b \in \mathcal{C}^{1}\left(\mathbb{R}\right)$, and the functions $b$ and $b^{\prime}$ are bounded on $\mathbb{R}$. Moreover,
    $\sigma \in \mathcal{C}^{2}\left(\mathbb{R}\right)$ such that $\sigma^{\prime}$ and $\sigma^{\prime\prime}$ are bounded on $\mathbb{R}$, and
    $$\forall ~ A \in \mathbb{R}\setminus\{0\}, ~~ \exp\left[-\dfrac{1}{2}\left(\int_{0}^{A}{\dfrac{1}{\sigma(u)}du}\right)^2\right] \leq \exp\left(-\dfrac{A^2}{2}\right).$$
\end{assumption}
Under Assumption~\ref{ass:LipFunctions}, there exists a unique strong solution $X=(X_t)_{t\in[0,1]}$ of Equation \eqref{eq:model}. Under Assumptions~\ref{ass:LipFunctions} and \ref{ass:Ellipticity}, the unique strong solution $X$ admits a transition density $(s, t, x, y)\mapsto p_X(s, t, x, y)$ which can be approximated by Gaussian densities when the numerical constants $\kappa_0$ and $\kappa_1$ satisfy $\kappa_0 = 1/\kappa_1$ and $\kappa_1 > 1$ (see, for example, \cite{gobet2002lan}). Furthermore, we draw from Assumptions~\ref{ass:LipFunctions} and \ref{ass:Ellipticity} the following important results:
\begin{equation}\label{eq:Cons-Ass2.1}
    \forall q\geq 1, \ \ \mathbb{E}\left[\underset{t\in[0,1]}{\sup}{|X_t|^q}\right]<\infty, \ \ \ \forall s,t>0, \ \mathbb{E}\left|X_t-X_s\right|^{2q}\leq C|t-s|^{q},
\end{equation}
where $C>0$ is a constant. In the present paper, we cannot rely, as in \cite{ella2024nonparametric}, on the approximation of the transition density of the process $X$ provided in \cite{gobet2002lan} for the study of a minimax convergence rate of the risk of estimation of $\sigma^2$ on a compact interval and on the growing interval $[-A_N, A_N]$, or for the establishment of a better convergence rate on the real line compared to the rate obtained in \cite{ella2024nonparametric}. In fact, \cite{gobet2002lan} provides the following approximation of the transition density $p_X$:
\begin{equation}\label{eq:approx-gobet}
    \dfrac{1}{C\sqrt{t-s}}\exp\left(-c_{\sigma}\dfrac{(y-x)^2}{t-s}\right) \leq p_X(s,t,x,y) \leq \dfrac{C}{\sqrt{t-s}}\exp\left(-\dfrac{(y-x)^2}{c_{\sigma}(t-s)}\right),
\end{equation}
where $x,y \in \mathbb{R}, ~~ s,t \in [0,1]$ such that $s<t$, and $C,c_{\sigma}>1$ are unknown constants, with $c_{\sigma}>1$ depending on $\sigma$. The drawback with this approximation is that the Gaussian densities that approximate $p_X$ depend on the unknown parameter $c_{\sigma}>1$ that cannot be controlled. This condition makes it difficult to establish an optimal convergence rate when $\sigma^{2}$ is estimated on the real line, or similarly, on the growing interval $[-A_N, A_N]$ with $A_N \propto \sqrt{\log(N)}$, a condition on $A_N$ that is proven to be optimal (see \cite{denis2024nonparametric}). A further discussion on this topic is provided in Section~\ref{subsec:Nonparametric_estimation_sigma}. We rather use the exact formula of the transition density $p_X$ established in \cite{dacunha1986estimation} under Assumptions \ref{ass:LipFunctions} and \ref{ass:Ellipticity}, which does not require additional conditions on the numerical constants $\kappa_0$ and $\kappa_1$, and which depends only on parameters that can be controlled under Assumptions~\ref{ass:LipFunctions}, \ref{ass:Ellipticity}, and \ref{ass:Restrict-Model}. 

Assumption~\ref{ass:RegularityBis} is used for the construction, from Equation~\eqref{eq:model}, of a regression model for the nonparametric estimation of $\sigma^2$, and for the establishment of upper bounds of its estimation risks for each of the cases $N \rightarrow \infty$ and $N = 1$.

Finally, Assumption~\ref{ass:Restrict-Model} is only used for the study of a convergence rate for the estimation of $\sigma^2$ on the real line. In fact, we show that the risk bound of the projection estimator of $\sigma^2$ on the real line is dominated by
\begin{equation*}
    \log^{\beta}(N)(Nn)^{-\beta/(2\beta+1)} + \log(N)\underset{t \in [0,1]}{\sup}{\mathbb{P}\left(|X_t| > A_N\right)},
\end{equation*}
where the smoothness parameter $\beta \geq 1$ of the H\"older class $\Sigma(\beta, R)$ is assumed to be large enough, and $\underset{t \in [0,1]}{\sup}{\mathbb{P}\left(|X_t| > A_N\right)}$ is the exit probability of the process $X = (X_t)_{t \in [0,1]}$ from the growing interval $[-A_N,A_N]$. Then, from the exact formula of the transition density together with Assumption~\ref{ass:Restrict-Model}, we establish a rate of convergence of the risk of estimation, which we prove to be faster than the rate provided in \cite{ella2024nonparametric}. Moreover, Assumption~\ref{ass:Restrict-Model} implies Assumptions~\ref{ass:LipFunctions} and \ref{ass:RegularityBis}, and requires Assumption~\ref{ass:Ellipticity} with $\kappa_1 = 1$.

\subsection{Definitions and notations}
\label{subsec:definitions and notations}

Consider the sample paths $\left\{\bar{X}^{j}, \ j=1,\cdots,N\right\}$ composed of $N$ independent copies of the discrete observation $\bar{X}$ of the diffusion process $X$. For simplicity, we denote by $\mathbb{P}$ the joint distribution of $\left(\bar{X}^{1},\cdots,\bar{X}^{N},\bar{X}\right)$ and $\mathbb{E}$ its corresponding expectation. We also denote by $\mathbb{P}_X$ the marginal distribution of the diffusion process $X$ and $\mathbb{E}_X$ its related expectation. For all continuous functions $h$ such that for all $t\in[0,1], ~ \mathbb{E}_X\left[h^{2}(X_t)\right]<\infty$, we define the following empirical and pseudo norms:
\begin{align*}
    \|h\|^{2}_{n}:=\mathbb{E}_X\left[\frac{1}{n}\sum_{k=0}^{n-1}{h^{2}\left(X_{k\Delta_n}\right)}\right], ~~~ \|h\|^{2}_{n,N}:=\frac{1}{Nn}\sum_{j=1}^{N}{\sum_{k=0}^{n-1}{h^{2}(X^{j}_{k\Delta_n})}}, ~~~ \|h\|_X^2 := \int_{0}^{1}{h^{2}(X_t)dt}.
\end{align*}
Then, for any deterministic function $h$ such that $\underset{t\in[0,1]}{\sup}{\mathbb{E}_X\left[h^{2}(X_t)\right]}<\infty$, one has $\mathbb{E}\left[\|h\|^{2}_{n,N}\right] = \|h\|^{2}_{n}$ and,
\begin{align*}
    \|h\|^{2}_{n} = \dfrac{h^{2}(x_0)}{n} + \int_{\R}{h^{2}(x)\frac{1}{n}\sum_{k=1}^{n-1}{p_X(0, k\Delta_n,x_0,x)}dx} = \dfrac{h^{2}(x_0)}{n} + \int_{\R}{h^{2}(x)f_n(x)dx},
\end{align*}
where $f_n: x\mapsto\frac{1}{n}\sum_{k=1}^{n-1}{p_X(0, k\Delta_n,x_0,x)}$. 

For any square matrix $M \in \mathcal{M}_r(\mathbb{R})$ of size $r \geq 2$, we denote by $M^{\prime}$, the transpose matrix of $M$. Moreover, when $M$ is a symmetrical and positive definite matrix, that is, for all $v \in \mathbb{R}^r, ~ v^{\prime}Mv > 0$, the operator norm of $M$ is defined by
\begin{equation*}
    \left\|M^{-1}\right\|_{\mathrm{op}} := \dfrac{1}{\min\left\{\lambda_1, \ldots, \lambda_d\right\}},
\end{equation*}
where $\{\lambda_1, \ldots, \lambda_d\}$ is the spectrum of $M$, with $d \leq r$.

For any integers $p,q \in \mathbb{N}$ such that $q > p$, we denote by $[\![p, q]\!]$ the set $\{p, p+1, \ldots, q\}$. Finally, for simplicity and without loss of generality, we set $X_0 = x_0 = 0$. In fact, for any diffusion process $X$, solution of Equation~\eqref{eq:model} with $X_0 = x_0 \in \mathbb{R}\setminus\{0\}$, we consider the diffusion process $Y = (X_t-x_0)_{t \in [0,1]}$, solution of the stochastic differential equation
\begin{equation*}
    dY_t = \widetilde{b}(Y_t)dt + \widetilde{\sigma}(Y_t)dW_t, ~~ Y_0 = 0,
\end{equation*}
where $\widetilde{b}: x \mapsto b(x+x_0)$ and $\widetilde{\sigma}: x \mapsto \sigma(x+x_0)$.

\subsection{Approximation spaces}
\label{subsec:data dependent approximation spaces}

In this paper, we assume that $\sigma^{2}$ belongs to the H\"older space $\Sigma_I(\beta, R)$ of smoothness parameter $\beta \geq 1$, with $R > 0$, defined as follows:
\begin{equation*}
    \Sigma_I(\beta, R) := \left\{f \in \mathcal{C}^{d}(I, \mathbb{R}), ~ \left|f^{(d)}(x) - f^{(d)}(y)\right| \leq R|x - y|^{\beta - d}, ~~ x,y \in \mathbb{R}\right\},
\end{equation*}
where $d = \lfloor \beta \rfloor$ is the highest integer strictly smaller than $\beta$, $I \subset \mathbb{R}$, and for $I = \mathbb{R}$, we set $\Sigma_I(\beta, R) = \Sigma(\beta, R)$. We approximate the space $\Sigma_{I}(\beta, R)$ by a $m$-dimensional space spanned by a basis $\left(\phi_0, \ldots, \phi_{m-1}\right)$ of Lipschitz functions defined on $I$, with $m \geq 1$, and satisfying the following condition:
\begin{equation}\label{eq:condition-basis}
    \exists r \in \mathbb{N}, ~ \mathcal{L}(m) := \underset{x \in I}{\sup}{\sum_{\ell = 0}^{m-1}{\phi_{\ell}^2(x)}} = \mathrm{O}\left(m^r\right), ~~ \exists s \in \mathbb{N}, ~ \mathcal{R}(m) := \underset{x \in I}{\sup}{\sum_{\ell = 0}^{m-1}{\phi_{\ell}^{\prime 2}(x)}} = \mathrm{O}\left(m^s\right).
\end{equation}
\begin{remark}
    Note that in the context of the estimation of $\sigma_{|I}^2$ from a single diffusion path ($N=1$), Condition~\eqref{eq:condition-basis} is required to establish a convergence rate of order $n^{-\beta/(2\beta+1)}$ over the H\"older space $\Sigma_I(\beta, R)$ where $I$ is a compact interval (see \cite{ella2024nonparametric}, \textit{Theorem 3.2 and Corollary 3.3}). For the sake of simplicity, only bases that satisfy Condition~\eqref{eq:condition-basis} are considered in the paper.
\end{remark} 
In the sequel, we focus on the spline basis and the Fourier basis which both satisfy Condition~\eqref{eq:condition-basis}. Note that the results obtained with the Fourier basis can easily be extended to any other orthonormal basis satisfying Condition~\eqref{eq:condition-basis}.

\paragraph{Spline basis ([B]).}

Let $K \geq 1$ and $M \geq 1$ be two integers. We approximate the H\"older space $\Sigma_I(\beta, R)$ by a $K+M$-dimensional subspace $\mathcal{S}_{K}$ spanned by the \textbf{B}-spline basis, given as follows:
\begin{equation}
    \label{eq:approximation subspace}
    \mathcal{S}_{K} := \mathrm{Span}\left(B_{\ell}, \ \ \ell = -M,\cdots, K-1\right).
\end{equation}
Let $A, B > 0$ be two real numbers such that $A<B$. The \textbf{B}-spline basis $\left\{B_{-M}, \ldots, B_{K-1}\right\}$ is built on the compact interval $[A, B]$, with $\mathbf{u} = \left(u_{-M}, \ldots, u_{K+M}\right)$ the knots vector such that $u_{-M} = \cdots = u_{-1} = u_0 = A$, $u_{K+1} = \cdots = u_{K+M} = B$, and for all $i = 0, \cdots, K$,
\begin{align*}
    u_i = A + i\frac{B-A}{K}.
\end{align*}
The \textbf{B}-spline functions $B_{\ell}, ~ \ell = -M, \ldots, K-1$ are positive and compactly supported piecewise polynomials of degree $M$ that satisfy the following property:  
$$\forall ~ x \in [A, B], ~~ \sum_{k = -M}^{K-1}{B_{\ell}(x)} = 1.$$
Thus, each function $h = \sum_{\ell = -M}^{K-1}{a_{\ell}B_{\ell}} \in \mathcal{S}_{K}$ is $M-1$ times continuously differentiable. For more details on the \textbf{B}-spline basis, we refer the reader to \cite{deboor78}, \cite{gyorfi2006distribution}.

Furthermore, we need to control the $\ell^2-$norm of the coordinate vectors of elements of $\mathcal{S}_{K}$, which leads to the following constrained subspace,
\begin{equation*}
    \mathcal{S}_{K, A,B}:=\left\{h=\sum_{\ell=-M}^{K-1}{a_{\ell}B_{\ell}}, \ \ \sum_{\ell=-M}^{K-1}{a^{2}_{\ell}} \leq (K+M)(B-A)^2\log(Nn)\right\}.
\end{equation*}
The control of the coordinate vectors leads to the construction of ridge estimators of $\sigma^{2}$ (see \cite{massart2006risk}). It is also required to establish upper bounds of the risk of estimation of the projection estimator. 

Finally, from \cite{denis2020ridge}, \textit{Lemma D.2}, the approximation of the H\"older space $\Sigma_I(\beta, R)$ by the finite-dimensional space $\mathcal{S}_{K, A,B}$ induces the following bound for the bias term:
\begin{equation}\label{eq:Bias-Term}
    \underset{h \in \mathcal{S}_{K, A, B}}{\inf}{\|h - \sigma_{|I}^2\|_n^2} \leq C(B-A)^{2\beta}K^{-2\beta},
\end{equation}
where $C>0$ is a constant and $I = [A, B]$ is the estimation interval.

\paragraph{Fourier basis ([F]).}

Let $A, B, D \geq 1$ be three integers, and $I = [A,B]$. We denote by $\mathcal{S}_D$ the space of approximation of the H\"older space $\Sigma_{I}(\beta, R)$ of dimension $m = 2D + 1$ spanned by the Fourier basis $\left\{1, f_{\ell}, g_{\ell}, ~ \ell = 1, \ldots, D \right\}$ dilated in the compact interval $I = [A,B]$, where for all $\ell \in [\![1, D]\!]$, the functions $f_{\ell}$ and $g_{\ell}$ read as follows:
\begin{equation*}
    f_{\ell}(x) := \dfrac{\sqrt{2}}{B-A}\cos\left(2\pi\ell\dfrac{x-A}{B-A}\right), ~~ g_{\ell}(x) := \dfrac{\sqrt{2}}{B-A}\sin\left(2\pi\ell\dfrac{x-A}{B-A}\right), ~~ x \in I.
\end{equation*}
We then deduce the following constrained approximation subspace,
\begin{equation*}
    \mathcal{S}_{D,A,B} := \left\{h = \gamma_0 + \sum_{\ell = 1}^{D}{\left(\alpha_{\ell}f_{\ell} + \lambda_{\ell}g_{\ell}\right)}, ~~ \gamma_0^2 + \sum_{\ell = 1}^{D}{\alpha_{\ell}^2} + \sum_{\ell = 1}^{D}{\lambda_{\ell}^2} \leq (2D+1)(B-A)^2\log(Nn)\right\}.
\end{equation*}
From \cite{barron1999risk}, \textit{Lemma 12} adapted to the compact interval $[A, B]$, there exists a constant $C>0$ such that 
\begin{equation}\label{eq:Bias-Fourier}
    \underset{h \in \mathcal{S}_{D,A,B}}{\inf}{\left\|h-\sigma_{|I}^2\right\|_n^2} \leq C(B-A)^{-1}D^{-2\beta}.
\end{equation}
More generally, we consider the following subspaces. 
\begin{equation}\label{eq:constrainded approximation subspace}
    \mathcal{S}_{m,A,B} := \left\{h = \sum_{\ell = 0}^{m-1}{\mu_{\ell}\phi_{\ell}}, ~~ \sum_{\ell = 0}^{m-1}{a_{\ell}^2} \leq m(B-A)^2\log(Nn)\right\}, ~ m \geq 1,
\end{equation}
where for the spline basis $[\mathbf{B}]$, $m = K+M$ with 
$$\mu_{\ell} = a_{\ell - M} ~~ \mathrm{and} ~~ \phi_{\ell} = B_{\ell - M}, ~~ \ell \in [\![0,m-1]\!],$$ 
and for the Fourier basis $[\mathbf{F}]$, $m = 2D+1$ with, 
\begin{align*}
    &~ \mu_0 = \gamma_0 ~~ \mathrm{and} ~~ \phi_0 = 1,\\
    &~ \forall \ell \in [\![1, D]\!], ~~ \mu_{\ell} = \alpha_{\ell} ~~ \mathrm{and} ~~ \phi_{\ell} = f_{\ell},\\
    &~ \forall \ell \in [\![D+1, 2D]\!], ~~ \mu_{\ell} = \lambda_{\ell} ~~ \mathrm{and} ~~ \phi_{\ell} = g_{\ell}.
\end{align*}
In the sequel, since two simulation schemes are considered for the estimation of $\sigma^2$, the dimension $m$ depends on $n$ when $N = 1$ and $n \rightarrow \infty$. For repeated observations ($N,n \rightarrow \infty$), we assume that $N$ and $n$ are proportional, and then the dimension $m$ can only depend on $N$.

Projection estimators of $\sigma^2$ on the bases $[\mathbf{B}]$ and $[\mathbf{F}]$ are studied in \cite{ella2024nonparametric}. In this paper, and for the estimation of $\sigma^2$ on a compact interval, we extend the study carried out in \cite{ella2024nonparametric} establishing the optimality of the obtained rates for each of the two simulation schemes $N=1, ~ n \rightarrow \infty$ and $N,n \rightarrow \infty$. When it comes to the estimation of $\sigma^2$ on $\mathbb{R}$, the risk bounds were established in \cite{ella2024nonparametric} for the bases $[\mathbf{B}]$ and $[\mathbf{F}]$ and the Hermite basis. In this paper, we prove that the obtained rates are sub-optimal, being the result of a rough bounding of the estimation risk. However, the proof technique used in this article for the study of a minimax rate of convergence on the growing interval $[-A_N, A_N]$ and on the real line requires an approximation space spanned by a basis of positive functions and whose associated Gram matrix satisfies certain key assumptions that are presented in Section~\ref{subsec:Nonparametric_estimation_sigma}. Therefore, only the spline basis is considered, satisfying the required assumptions.

\subsection{Nonparametric estimation of the square of the diffusion coefficient}
\label{subsec:Nonparametric_estimation_sigma}

In this section, we propose a projection estimator of $\sigma^2$ over a finite-dimensional subspace $\mathcal{S}_{m, A, B}$ given by Equation~\eqref{eq:constrainded approximation subspace} in Section~\ref{subsec:data dependent approximation spaces}. For this purpose, we suppose to have at our disposal a sample 
\begin{equation}\label{eq:sample}
    \mathfrak{S}_{n,N} = \left\{X^{1}_{\Delta_n}, \ldots, X^{1}_{n\Delta_n}, \ldots, X^{N}_{\Delta_n}, \ldots, X^{N}_{n\Delta_n}\right\}
\end{equation} 
constituted of $N\times n$ observations coming from $N$ independent copies $\bar{X}^{1} = (X^{1}_{k\Delta_n})_{0\leq k\leq n}, \ldots, \bar{X}^{N} = (X^{N}_{k\Delta_n})_{0\leq k\leq n}$ of $\bar{X} = (X_{k\Delta_n})_{0\leq k\leq n}$, discrete observation of the diffusion process $X$, solution of Equation~\eqref{eq:model}, with time-step $\Delta_n = 1/n$.  

From Assumption~\ref{ass:LipFunctions} and Equation~\eqref{eq:model}, we obtain from the sample $\mathfrak{S}_{n,N}$ that the square $\sigma^2$ of the diffusion coefficient is a solution of the following regression model:
\begin{equation*}
  U^{j}_{k\Delta_n}=\sigma^{2}(X^{j}_{k\Delta_n})+\zeta^{j}_{k\Delta_n}+R^{j}_{k\Delta_n}, ~~ \forall (j,k)\in[\![1,N]\!]\times[\![0,n-1]\!],
\end{equation*}
where the response variable $U^{j}_{k\Delta_n}$ is given by
$$ U^{j}_{k\Delta_n} := \frac{\left(X^{j}_{(k+1)\Delta_n} - X^{j}_{k\Delta_n}\right)^2}{\Delta_n},$$
and the error terms are respectively given by $\zeta^{j}_{k\Delta_n}=\zeta^{j,1}_{k\Delta_n}+\zeta^{j,2}_{k\Delta_n}+\zeta^{j,3}_{k\Delta_n}$, with:
\begin{equation*}
    \zeta^{j,1}_{k\Delta_n}=\frac{1}{\Delta_n}\left[\left(\int_{k\Delta_n}^{(k+1)\Delta_n}{\sigma(X^{j}_{s})dW^{j}_{s}}\right)^2-\int_{k\Delta_n}^{(k+1)\Delta_n}{\sigma^{2}(X^{j}_{s})ds}\right],
\end{equation*}
\begin{equation*}
    \zeta^{j,2}_{k\Delta_n}=\frac{2}{\Delta_n}\int_{k\Delta_n}^{(k+1)\Delta_n}{((k+1)\Delta_n-s)\sigma^{\prime}(X^{j}_{s})\sigma^{2}(X^{j}_{s})dW^{j}_{s}},
\end{equation*}
\begin{equation*}
\zeta^{j,3}_{k\Delta_n}=2b(X^{j}_{k\Delta_n})\int_{k\Delta_n}^{(k+1)\Delta_n}{\sigma\left(X^{j}_{s}\right)dW^{j}_{s}},
\end{equation*}
and $R^{j}_{k\Delta_{n}}
=R^{j,1}_{k\Delta_{n}}+R^{j,2}_{k\Delta_{n}}+R^{j,3}_{k\Delta_{n}}$, with:
\begin{equation*}
    R^{j,1}_{k\Delta_n}=\frac{1}{\Delta_n}\left(\int_{k\Delta_n}^{(k+1)\Delta_n}{b(X^{j}_{s})ds}\right)^2, ~~~~ R^{j,3}_{k\Delta_n} = \frac{1}{\Delta_n}\int_{k\Delta_n}^{(k+1)\Delta_n}{((k+1)\Delta_n-s)\Phi(X^{j}_{s})ds},
\end{equation*}
\begin{equation*}
R^{j,2}_{k\Delta_n}=\frac{2}{\Delta_n}\left(\int_{k\Delta_n}^{(k+1)\Delta_n}{\left(b(X^{j}_{s})-b(X^{j}_{k\Delta_n})\right)ds}\right)\left(\int_{k\Delta_n}^{(k+1)\Delta_n}{\sigma(X^{j}_{s})dW^{j}_{s}}\right),
\end{equation*}
where 
\begin{equation*} \Phi:=2b\sigma^{\prime}\sigma+\left[\sigma^{\prime\prime}\sigma+\left(\sigma^{\prime}\right)^2\right]\sigma^{2}.
\end{equation*}
The proof of this result is provided in \cite{ella2024nonparametric}, \textit{Lemma 10}.

Let $\gamma_{n,N}: h \in \mathcal{S}_{m} \mapsto \gamma_{n,N}(h)$ be a contrast function defined from the sample $\mathfrak{S}_{n,N}$ by
\begin{equation*}
    \forall ~ h \in \mathcal{S}_{m}, ~~ \gamma_{n,N}(h) := \dfrac{1}{Nn}\sum_{j=1}^{N}{\sum_{k=0}^{n-1}{\left(U^{j}_{k\Delta_n} - h(X^{j}_{k\Delta_n})\right)^2}}.
\end{equation*}
Thus, the nonparametric estimators $\widehat{\sigma}^{2}_{m}$ of $\sigma^2$ are built by minimizing the contrast function $\gamma_{n,N}$ on the subspace $\mathcal{S}_{m, A, B}$, that is,
\begin{equation}\label{eq:ProjEstim}
    \forall ~ m \geq 1, ~~ \widehat{\sigma}^{2}_{m} = \underset{h \in \mathcal{S}_{m,A,B}}{\arg\min}{~\gamma_{n,N}(h)}.
\end{equation}
Consider the Gram matrix $\mathbf{\Psi}_{m}$ related to the basis $\left(\phi_0, \ldots, \phi_{m-1}\right)$ and given by
\begin{equation}\label{eq:GramMatrix}
    \mathbf{\Psi}_{m} := \left(\mathbb{E}_X\left[\dfrac{1}{n}\sum_{k=0}^{n-1}{\phi_{\ell}(X_{k\Delta_n})\phi_{\ell^{\prime}}(X_{k\Delta_n})}\right]\right)_{0 \leq \ell, \ell^{\prime} \leq m-1}.
\end{equation}
The Gram matrix $\Psi_{m}$, which is symmetric, is proven to be positive definite for each of the two bases $[\mathbf{B}]$ and $[\mathbf{F}]$ (see the proof of Lemma~\ref{lm:OperatorNorm-Psi}), which implies that $\Psi_m$ is invertible. We obtain the following result.

\begin{lemma}\label{lm:OperatorNorm-Psi}
    Under Assumptions~\ref{ass:LipFunctions}, \ref{ass:Ellipticity} and \ref{ass:Restrict-Model} and for the basis $[\mathbf{B}]$ and $[\mathbf{F}]$, there exist constants $C,c>0$ such that
    \begin{equation*}
        \mathcal{L}(m)\left\|\mathbf{\Psi}^{-1}_{m}\right\|_{\mathrm{op}} \leq C\dfrac{m\log(N)}{A_{N}}\exp\left(\dfrac{1}{2(1 - \log^{-1}(N))}\left(\int_{0}^{A_{N}}{\dfrac{1}{\sigma(u)}du}\right)^2  + cA_{N}\right),
    \end{equation*}
    where for all $m \geq 1$,
    \begin{equation*}
        \mathcal{L}(m) := \underset{x \in \mathbb{R}}{\sup}{\sum_{\ell=0}^{m-1}{\phi_{\ell}^2(x)}}.
    \end{equation*}
\end{lemma}
The establishment of convergence rates of projection estimators of $\sigma^2$ requires the study of an equivalence relation between the empirical norm $\|.\|_n$ and the pseudo-norm $\|.\|_{n,N}$ defined in Section~\ref{subsec:definitions and notations}. In the context of estimation of $\sigma^2$ on the real line $\mathbb{R}$, Lemma~\ref{lm:OperatorNorm-Psi} provides a key result that is crucial to prove that the probability of the complementary of the following random event 
\begin{equation*}
    \Omega_{n,N, m} := \underset{h \in \mathcal{S}_{m} \setminus \{0\}}{\bigcap}{\left\{\left| \dfrac{\|h\|_{n,N}^2}{\|h\|_n^2} - 1\right| \leq \dfrac{1}{2}\right\}}
\end{equation*}
is upper bounded by a negligible term in the risk bound of the estimation risk of $\widehat{\sigma}_{m}^2$. Consider, for instance, a diffusion model whose coefficients are given by
\begin{equation*}
    b: x \mapsto \dfrac{1}{2 + \cos(x)}, ~~ \sigma: x \mapsto \sqrt{\dfrac{4}{5}} + \dfrac{1}{4\pi + x^2}.
\end{equation*}
The functions $b$ and $\sigma$ satisfy Assumptions~\ref{ass:LipFunctions}, \ref{ass:Ellipticity} and \ref{ass:RegularityBis}. Moreover, if $\|\sigma\|_{\infty} \leq 1$, then Assumption~\ref{ass:Restrict-Model} is also satisfied. Setting $n \propto N$, $A_N = \sqrt{\frac{3\beta}{2\beta+1}\log(N)}$, $m = N^{2/(2\beta+1)}\log^{-5/2}(N)$ with $\beta > 8$, and $N \rightarrow \infty$, there exist constants $C_1>0$ and $C_2>0$ such that
\begin{multline*}
   \mathcal{L}(m)\left\|\mathbf{\Psi}^{-1}_{m}\right\|_{\mathrm{op}} \leq C\dfrac{m\log(N)}{A_{N}}\exp\left(\dfrac{5}{8}\dfrac{3\beta}{2\beta-1}\dfrac{2\beta-1}{2\beta+1}\log(N) + cA_N\right)\\
    \leq C_1\dfrac{N}{\log^{2}(N)}\exp\left(\left(\dfrac{15\beta}{16\beta-8} - 1\right)\log(N) + cA_N\right) \leq C_2\dfrac{N}{\log^2(N)}.
\end{multline*}
If we set $A_N = \sqrt{\frac{1}{2\beta+1}\log(N)}$ and $m = N^{2/(2\beta+1)}\log^{-5/2}(N)$, then the same result is obtained with any smoothness parameter $\beta \geq 1$.

Note that the novelty in this paper is the use of the exact formula of the transition density of the unique strong solution $X$ of Equation~\eqref{eq:model} provided in \cite{dacunha1986estimation}, which contains parameters that can be controlled. In contrast, in \cite{ella2024nonparametric}, the upper bounding of $\left\|\Psi_{m}^{-1}\right\|_{\mathrm{op}}$ is based on the approximation of the transition density $(s,t,x,y) \mapsto p_X(s,t,x,y)$ provided in \cite{gobet2002lan} (see Equation~\eqref{eq:approx-gobet}). This approximation leads to the following result:
\begin{equation*}
    \mathcal{L}(m)\left\|\mathbf{\Psi}^{-1}_{m}\right\|_{\mathrm{op}} \leq C\dfrac{m\log(N)}{A_{N}}\exp\left(c_{\sigma}A_N^2\right).
\end{equation*}
Since the unknown constant $c_{\sigma} > 1$ cannot be controlled, the condition 
\begin{equation}\label{eq:Gram-Condition}
    \mathcal{L}(m)\left\|\mathbf{\Psi}^{-1}_{m}\right\|_{\mathrm{op}}  = \mathrm{O}\left(\dfrac{N}{\log^2(N)}\right)
\end{equation}
holds only when $A_N/\sqrt{\log(N)}$ tends to zero as $N$ tends to infinity, implying that $A_N = o\left(\sqrt{\log(N)}\right)$. This condition on $A_N$ cannot allow us to establish, under the assumptions of this paper, a better convergence rate for the estimation of $\sigma^2$ on the real line compared to the rate obtained in \cite{ella2024nonparametric}.

In the sequel, we study the minimax convergence rates of the projection estimators $\widehat{\sigma}^{2}_{m}, ~ m \geq 1$ both on any compact interval and on the real line $\mathbb{R}$ under the assumption on the Gram matrix $\Psi_m, ~ m\geq 1$ given by Equation~\eqref{eq:Gram-Condition}.

\section{Upper bound of the risk of estimation of the the diffusion coefficient on the real line}
\label{sec:upper-bound}

This section is devoted to the study of the upper bound of the risk of the nonparametric estimator $\widehat{\sigma}_{m}^2$ of $\sigma^{2}$ when $N,n \rightarrow \infty$. Two different cases are considered. First, we establish a convergence rate of the risk 
\begin{align*}
     \mathbb{E}\left[\left\|\widehat{\sigma}_{m}^2 - \sigma_{A_N}^2\right\|_n^2\right],
\end{align*} 
where $\sigma_{A_N}^2 = \sigma^2\mathds{1}_{[-A_N, A_N]}$ is the restriction of $\sigma^2$ on the growing interval $[-A_N, A_N]$. We then extend the study to the establishment of a rate of convergence of the risk
\begin{align*}
     \mathbb{E}\left[\left\|\widehat{\sigma}_{m}^2 - \sigma^2\right\|_n\right].
\end{align*}
As we stated in Section~\ref{subsec:data dependent approximation spaces}, the proof technique requires a basis of positive functions whose related Gram matrix $\Psi_m$ is invertible and satisfies Lemma~\ref{lm:OperatorNorm-Psi}. Thus, in this section, only the spline basis is considered. We establish below the first risk bound for the projection estimator of $\sigma^{2}$ on the growing interval $[-A_N,A_N]$.

\begin{theo}\label{thm:upper-bound}
    Suppose that $A_N \propto \sqrt{\log(N)}$, $n \propto N, ~ K_N \propto N^{2/(2\beta+1)}\log^{-5/2}(N), ~ N \rightarrow \infty$, and the Gram matrix $\mathbf{\Psi}_{K_N}$ given in Equation~\eqref{eq:GramMatrix} satisfies
    \begin{equation*}
        \mathcal{L}(K_N+M)\left\|\mathbf{\Psi}^{-1}_{K_N}\right\|_{\mathrm{op}} = \mathrm{O}\left(\dfrac{N}{\log^{2}(N)}\right).
    \end{equation*}
    Under Assumptions~\ref{ass:LipFunctions}, \ref{ass:Ellipticity}, and  \ref{ass:Restrict-Model}, there exists a constant $C>0$ such that
    \begin{align*}
      \underset{\sigma^2 \in \Sigma(\beta, R)}{\mathrm{sup}}{~\mathbb{E}\left[\left\|\widehat{\sigma}^{2}_{K_N}-\sigma^{2}_{A_N}\right\|_{n}^2\right]} \leq C\log^{\beta}(N)(Nn)^{-2\beta/(2\beta+1)},
    \end{align*}
    where $I = [-A_N, A_N]$ and $\beta \geq 1$.
\end{theo}

The above result is established on the growing interval $[-A_N, A_N]$, where $A_N$ is proportional to $\sqrt{\log(N)}$, contrary to the result established in \cite{ella2024nonparametric} assuming that
\begin{equation*}
    \dfrac{A_N}{\sqrt{\log(N)}} \longrightarrow 0 ~~ \mathrm{as} ~~ N \rightarrow \infty.
\end{equation*}
Note that when the dimension is set to $m = N^{2/(2\beta+1)}\log^{-5/2}(N)$,  the choice of $A_N$ depends on the smoothness of the square $\sigma^2$ of the diffusion coefficient. For example, choosing $A_N = \sqrt{\log(N)/(2\beta+1)}$ is suitable for any smoothness parameter $\beta \geq 1$, while choosing $A_N = \sqrt{3\beta\log(N)/(2\beta+1)}$ is appropriate when the smoothness parameter $\beta$ is large enough. However, the assumption $A_N = \sqrt{3\beta\log(N)/(2\beta+1)}$ is the one that allows us to derive a convergence rate, where $\sigma^2$ is estimated on the real line. For this purpose, consider, for $K_N \geq 1$, the truncated estimator $\widetilde{\sigma}_{K_N}^2$ defined from $\widehat{\sigma}_{K_N}^2$ and for $x \in \mathbb{R}$ by
\begin{equation*}
    \widetilde{\sigma}_{K_N}^2(x) = \widehat{\sigma}_{K_N}^2(x)\mathds{1}_{\widehat{\sigma}_{K_N}^2(x) \leq \log(N)} + \log(N)\mathds{1}_{\widehat{\sigma}_{K_N}^2(x) \geq \log(N)}. 
\end{equation*}
From the definition of the truncated estimator $\widetilde{\sigma}_{K_N}^2$, we have
\begin{equation}\label{eq:RelateTrunc}
    \mathbb{E}\left[\left\|\widetilde{\sigma}^{2}_{K_N}-\sigma^{2}\right\|_{n}\right] \leq \mathbb{E}\left[\left\|\widehat{\sigma}^{2}_{K_N}-\sigma^{2}\right\|_{n}\right].  
\end{equation}
We derive below a convergence rate of the worst risk of $\widetilde{\sigma}_{K_N}^2$ over the H\"older space of smoothness parameter $\beta$.

\begin{prop}\label{prop:Rate-realline}
    Suppose that there exists $\beta_0 > 1$ such that for $\beta \geq \beta_0$, $A_N = \sqrt{\dfrac{3\beta}{2\beta+1}\log(N)}$, $n \propto N, ~ K_N \propto \dfrac{N^{2/(2\beta+1)}}{\log^{5/2}(N)}$, $N \rightarrow \infty$, and the Gram matrix $\mathbf{\Psi}_{K_N}$ given by Equation~\eqref{eq:GramMatrix} satisfies
    \begin{equation*}
        \mathcal{L}(K_N+M)\left\|\mathbf{\Psi}^{-1}_{K_N}\right\|_{\mathrm{op}} = \mathrm{O}\left(\dfrac{N}{\log^{2}(N)}\right).
    \end{equation*}
    Under Assumptions~\ref{ass:LipFunctions}, \ref{ass:Ellipticity}, and  \ref{ass:Restrict-Model}, there exist constants $C, c>0$ such that
    \begin{align*}
      \underset{\sigma^2 \in \Sigma(\beta, R)}{\sup}{~\mathbb{E}\left[\left\|\widetilde{\sigma}^{2}_{K_N}-\sigma^{2}\right\|_{n}\right]} \leq C\exp\left(c\sqrt{\log(N)}\right)N^{-3\beta/2(2\beta+1)}.
    \end{align*}
\end{prop}

The above result provides a convergence rate of the risk of estimation of $\sigma^2$ on the real line. This rate is utterly faster than the one established in \cite{ella2024nonparametric}, which is of order $N^{-\beta/2(\beta+1)}$ up to a logarithmic factor. However, the obtained rate requires a sufficiently large smoothness parameter $\beta$ for the H\"older space $\Sigma(\beta, R)$ that contains $\sigma^2$. The constraint on $\beta$ is removed for $K_N = N^{1/(2\beta+1)}\log^{-5/2}(N)$ and $A_N = \sqrt{\frac{3\beta}{2\beta+1}\log(N)}$, but these assumptions lead to a slower rate of order $N^{-\beta/(2\beta+1)}$. However, this rate remains slower than that provided by Theorem~\ref{thm:upper-bound}. We show in the proof of Proposition~\ref{prop:Rate-realline} that the upper bound of the risk of estimation is dominated by the upper bound of 
\begin{equation}\label{eq:ExitProb}
    \underset{t \in [0,1]}{\sup}{~\mathbb{P}(|X_t| > A_N)},
\end{equation}
which is the probability of exit of the process $X$ from the growing interval $[-A_N, A_N]$. A convergence rate of the same order as that provided by Theorem~\ref{thm:upper-bound} cannot be established under Assumption~\ref{ass:Restrict-Model} as it would require $A_N \geq \sqrt{\frac{4\beta}{2\beta+1}\log(N)}$, which is not possible, as Equation~\eqref{eq:Gram-Condition} would no longer be satisfied from the result of Lemma~\ref{lm:OperatorNorm-Psi}. 

Note that establishing a convergence rate of order $\exp\left(c\sqrt{\log(N)}\right)N^{-3\beta/2(2\beta+1)}$ was possible under the assumptions of this article considering the risk of estimation $\mathbb{E}\left[\left\|\widetilde{\sigma}_{K_N}^2 - \sigma^2\right\|_n\right]$. If, on the other hand, we consider the risk of estimation $\mathbb{E}\left[\left\|\widetilde{\sigma}_{K_N}^2 - \sigma^2\right\|_n^2\right]$, then we obtain
\begin{align*}
    \underset{\sigma^2 \in \Sigma(\beta, R)}{\sup}{~\mathbb{E}\left[\left\|\widetilde{\sigma}^{2}_{K_N}-\sigma^{2}\right\|_{n}^2\right]} \leq C\exp\left(c\sqrt{\log(N)}\right)N^{-3\beta/2(2\beta+1)}.
\end{align*}
Thus, we derive a slower convergence rate of order $\exp\left(c\sqrt{\log(N)}\right)N^{-3\beta/4(2\beta+1)}$ which is no longer faster than the rate of order $N^{-\beta/2(\beta+1)}$ (up to a logarithmic factor) established in \cite{ella2024nonparametric}.

The next section focuses on the study of the lower bound of the risk of estimation of $\sigma^2$.

\section{Lower bound of the risk of estimation of the diffusion coefficient}
\label{sec:lower-bound}

In this section, we study the lower bound of the worst risk of estimation of the square of the diffusion coefficient from the sample $\mathfrak{S}_{n,N}$. In Section~\ref{subsec:lower-bound-compact}, we establish the lower bound of the risk of estimation of $\sigma^2$ on a compact interval $[A, B]$, where $A$ and $B$ are two real numbers independent of $N$ and $n$ and such that $A<B$. Section~\ref{subsec:lower-bound-Npaths-R} extends the result of Section~\ref{subsec:lower-bound-compact} to the estimation of $\sigma^{2}$ on the real line $\mathbb{R}$ as $N$ tends to infinity. 

\subsection{Lower bound of the risk of estimation of the diffusion coefficient on a compact interval.}
\label{subsec:lower-bound-compact}

We study the optimality of the convergence rates of the nonparametric estimators of $\sigma^2$ on a compact interval established in \cite{ella2024nonparametric} in the respective settings where $N,n \rightarrow \infty$ and where $N = 1$ and $n \rightarrow \infty$. In fact, in \cite{ella2024nonparametric}, the author built projection estimators of $\sigma^2$ on a fixed and compact interval of the form $I = [A,B]$ where $A,B \in \mathbb{R}$ such that $B > A$. We mean by fixed interval, the simple fact that $A$ and $B$ do not depend on the number $N$ of diffusion paths or the time step $\Delta_n = 1/n$. The author then established in \textit{Theorem 3.2 with Corollary 3.3, and Theorem 4.1 with Corollary 4.2}, under Assumptions~\ref{ass:LipFunctions}, \ref{ass:RegularityBis}, and Assumption~\ref{ass:Ellipticity} with $\kappa_1 > 1$ and $\kappa_0 = 1/\kappa_1$, the following results:
\begin{align*}
    \mathbb{E}\left[\left\|\widehat{\sigma}_{m}^2 - \sigma^{2}_{|I}\right\|_n^2\right] = &~  \mathrm{O}\left((Nn)^{-2\beta/(2\beta+1)}\right) ~~ \mathrm{when} ~~ N,n \rightarrow \infty ~~ \mathrm{and} ~~ N \propto n,\\
    \mathbb{E}\left[\left\|\widehat{\sigma}_{m}^2 - \sigma^{2}_{|I}\right\|_n^2\right] = &~  \mathrm{O}\left(n^{-2\beta/(2\beta+1)}\right) ~~ \mathrm{when} ~~ N = 1 ~~ \mathrm{and} ~~ n \rightarrow \infty,
\end{align*}
where $\beta > 3/2$, $m \propto (Nn)^{1/(2\beta+1)}$ for $N \rightarrow \infty$, and $m \propto n^{1/(2\beta+1)}$ for $N = 1$. Since the empirical norm $\|.\|_n$ and the $L^2-$ norm are equivalent in the space $L^2([A, B])$ of square integrable functions on the compact interval $I = [A, B]$, from the above results we deduce that
\begin{align*}
    \mathbb{E}\left[\left\|\widehat{\sigma}_{m}^2 - \sigma^{2}_{|I}\right\|^2\right] = &~  \mathrm{O}\left((Nn)^{-2\beta/(2\beta+1)}\right) ~~ \mathrm{when} ~~ N,n \rightarrow \infty ~~ \mathrm{and} ~~ N \propto n,\\
    \mathbb{E}\left[\left\|\widehat{\sigma}_{m}^2 - \sigma^{2}_{|I}\right\|^2\right] = &~  \mathrm{O}\left(n^{-2\beta/(2\beta+1)}\right) ~~ \mathrm{when} ~~ N = 1 ~~ \mathrm{and} ~~ n \rightarrow \infty.
\end{align*}
We derive below the lower bound of the worst risk of estimators of $\sigma^{2}$, constructed from the sample $\mathfrak{G}_{n,N}$, and defined with the $L^2$ norm $\|.\|$ when $N,n \rightarrow \infty$ and when $N = 1$ and $n \rightarrow \infty$.

\begin{theo}\label{thm:lower-bound-compact}
    Consider the estimation interval $I = [A, B]$, where $A, B \in \mathbb{R}$ such that $B > A$. Under Assumptions~\ref{ass:LipFunctions}, \ref{ass:Ellipticity} and \ref{ass:RegularityBis} , we obtain the following results.
    \begin{enumerate}
        \item For $N = 1, ~ n \rightarrow \infty$ and for $\beta > 2$, there exists a constant $c>0$ depending on $\beta$ such that
        \begin{align*}
            \underset{\w{\sigma}^{2}}{\inf}{\underset{\sigma^2 \in \Sigma_{I}(\beta,R)}{\sup}{~\E\left[\|\w{\sigma}^{2}-\sigma^{2}_{|I}\|^2\right]}} \geq cn^{-2\beta/(2\beta+1)}.
        \end{align*}
       \item For $N,n \rightarrow \infty$ with $n \propto N^2$, and for $\beta \geq 4$, there exists a constant $c>0$ depending on $\beta$ such that
       \begin{align*}
           \underset{\w{\sigma}^{2}}{\inf}{\underset{\sigma^2 \in \Sigma_{I}(\beta,R)}{\sup}{~\E\left[\|\w{\sigma}^{2}-\sigma^{2}_{|I}\|^2\right]}} \geq c(Nn)^{-2\beta/(2\beta+1)}.
       \end{align*}
    \end{enumerate}
\end{theo}
The result of Theorem~\ref{thm:lower-bound-compact} shows that, since the estimation interval $I = [A,B]$ does not depend on the sample size, the lower bound of the risk of estimation is of the same order as the upper bound established in \cite{ella2024nonparametric}. More precisely, considering the estimators $\widehat{\sigma}^2$ of $\sigma^2 \in \Sigma_I(\beta,R)$ defined by Equation~\eqref{eq:ProjEstim}, where $\sigma$ satisfies Assumptions~\ref{ass:LipFunctions}, \ref{ass:RegularityBis}, and Assumption \ref{ass:Ellipticity} with $\kappa_1 > 1$ and $\kappa_0 = 1/\kappa_1$, we deduce that for $\beta > 2, ~ N = 1$ and $n \rightarrow \infty$, there exist constants $C, c > 0$ such that
\begin{align*}
    cn^{-2\beta/(2\beta+1)} \leq \underset{\w{\sigma}^{2}}{\inf}{\underset{\sigma^2 \in \Sigma_I(\beta,R)}{\sup}{~\E\left[\|\w{\sigma}^{2}-\sigma^{2}_{|I}\|^2\right]}} \leq Cn^{-2\beta/(2\beta+1)}.
\end{align*}
Note that for the cases $N = 1$ and $n \rightarrow \infty$, an optimal convergence rate is established in \cite{hoffmann1999lp} where the square of the diffusion coefficient is built by projection on a finite-dimensional space spanned by wavelet basis. The optimal rate of order $n^{-s/(2s+1)}$ is obtained over a Besov space of smoothness parameter $s > 2$ taking advantage of the optimality of wavelet functions. \cite{ella2024nonparametric} extended the study to any basis of Lipschitz functions on a compact interval and established a convergence rate of the same order over the space of H\"older functions. In this paper, we have established the optimality of the rate derived in \cite{ella2024nonparametric}, \textit{Corollary 3.3}, in the context of the estimation of $\sigma^2$ on a compact interval and from a single diffusion path.

For the estimation of $\sigma^2$ on a compact interval and from repeated observations of the diffusion process $X$, different assumptions on the time step $\Delta_n$ are required for the establishment of the lower bound in this paper and the establishment of the upper bound in \cite{ella2024nonparametric}. In fact, in \cite{ella2024nonparametric} (\textit{Lemma 8.8, Theorem 4.1 and Corollary 4.2}), the study of an upper bound of order $(Nn)^{-\beta/(2\beta+1)}$ required the dimension $m$ of the approximation space to belong to the finite set $\left\{1, \ldots, \left\lfloor N^{1/2}\log^{-1}(Nn)\right\rfloor\right\}$, which implied that $m = (Nn)^{1/(2\beta+1)}$ with $n \propto N$ and $\beta > 3/2$. However, from the result of Theorem~\ref{thm:lower-bound-compact}, the establishment of a lower bound of order $(Nn)^{-\beta/(2\beta+1)}$ imposes that $n \propto N^2$. This condition is the result of the complexity of the diffusion model characterized by the explicit formula of the transition density $(s,t,x,y) \mapsto p_X(s,t,x,y)$ of the diffusion process $X$ provided in \cite{dacunha1986estimation}. We then have the following results under Assumptions~\ref{ass:LipFunctions}, \ref{ass:Ellipticity} and \ref{ass:RegularityBis}:
\begin{enumerate}
    \item[(i)] For $n \propto N$ and $\beta > 3/2$, there exists a constant $C>0$ such that (see \cite{ella2024nonparametric}, \textit{Corollary 4.2}):
    \begin{align*}
        \underset{\sigma^2 \in \Sigma_{I}(\beta,R)}{\sup}{~\E\left[\|\w{\sigma}^{2}-\sigma^{2}_{|I}\|^2\right]} \leq C(Nn)^{-2\beta/(2\beta+1)}.
    \end{align*}
    \item[(ii)] For $n \propto N^2$ and $\beta \geq 4$, there exists a constant $c>0$ such that:
    \begin{align*}
        \underset{\w{\sigma}^{2}}{\inf}{\underset{\sigma^2 \in \Sigma_{I}(\beta,R)}{\sup}{~\E\left[\|\w{\sigma}^{2}-\sigma^{2}_{|I}\|^2\right]}} \geq c(Nn)^{-2\beta/(2\beta+1)}.
    \end{align*}
\end{enumerate}

The hypotheses used for the study of the lower bound of the estimation risk are built from a $\beta-$H\"older function proposed in \cite{tsybakov2008introduction}. The strong assumption on the time step $\Delta_n = 1/n$ comes from the participation of the first and second derivatives of the hypotheses by the transition density of the diffusion process $X$ on one side, and the It\^o formula on the other side. 

\subsection{Lower bound of the risk of estimation of the diffusion coefficient on the real line from i.i.d. diffusion paths}
\label{subsec:lower-bound-Npaths-R}

This section is devoted to the establishment of a lower bound of the risk of estimation of $\sigma^2$ on the real line $\mathbb{R}$ and on the growing interval $[-A_N, A_N]$, where $N \rightarrow +\infty$. We obtain the following result.   

\begin{theo}\label{thm:lower-bound-R-Nn}
Suppose that $n \propto N^2$, $N \rightarrow \infty$ and grant the assumptions~\ref{ass:LipFunctions} and \ref{ass:Ellipticity}. The following holds,
\begin{equation*}
   \underset{\w{\sigma}^{2}}{\inf}{\underset{\sigma^2 \in \Sigma(\beta,R)}{\sup}{~\E\left[\left\|\w{\sigma}^2-\sigma^{2}\right\|_n^{2}\right]}} \geq c(Nn)^{-2\beta/(2\beta+1)},
\end{equation*}
where the constant $c>0$ depends on the smoothness parameter $\beta \geq 4$ of the H\"older space $\Sigma(\beta,R)$. Moreover, suppose that $A_N \propto \sqrt{\log(N)}$. There exists a constant $c>0$ depending on $\beta$ such that
\begin{equation*}
   \underset{\w{\sigma}^{2}}{\inf}{\underset{\sigma^2 \in \Sigma(\beta,R)}{\sup}{~\E\left[\left\|\w{\sigma}^2-\sigma^{2}_{A_N}\right\|_n^{2}\right]}} \geq c(Nn)^{-2\beta/(2\beta+1)},
\end{equation*}
where $\sigma_{A_N}^2 = \sigma^2\mathds{1}_{[-A_N, A_N]}$.
\end{theo}
We obtain a lower bound of the same order as that provided by Theorem~\ref{thm:lower-bound-compact}, where $\sigma^2$ is estimated on a compact interval. Then, combining Theorem~\ref{thm:upper-bound} and Theorem~\ref{thm:lower-bound-R-Nn}, we obtain the following results.
\begin{enumerate}
    \item[(i)] Under Assumptions~\ref{ass:LipFunctions}, \ref{ass:Ellipticity} and \ref{ass:Restrict-Model}, and for $n \propto N, ~ N \rightarrow \infty$ and $\beta$ large enough, there exists a constant $C>0$ such that
    \begin{equation*}
        \underset{\sigma^2 \in \Sigma(\beta, R)}{\mathrm{sup}}{~\mathbb{E}\left[\left\|\widehat{\sigma}^{2}_{K_N}-\sigma^{2}_{A_N}\right\|_{n}^2\right]} \leq C\log^{\beta}(N)(Nn)^{-2\beta/(2\beta+1)},
    \end{equation*}
    where $K_N = N^{2/(2\beta+1)}\log^{-5/2}(N)$ and $A_N \propto \sqrt{\log(N)}$.
    \item[(ii)] Under Assumptions~\ref{ass:LipFunctions} and \ref{ass:Ellipticity} with $n \propto N^2$ and $\beta \geq 4$, there exists a constant $c>0$ such that 
    \begin{equation*}
        \underset{\w{\sigma}^{2}}{\inf}{\underset{\sigma^2 \in \Sigma(\beta,R)}{\sup}{~\E\left[\left\|\w{\sigma}^2-\sigma^{2}_{A_N}\right\|_n^{2}\right]}} \geq c(Nn)^{-2\beta/(2\beta+1)},
    \end{equation*}
    where $A_N \propto \sqrt{\log(N)}$.
\end{enumerate} 
As we can see, the upper bound and the lower bound of the worst estimation risk differ from a logarithmic factor. The logarithmic factor on the upper bound of the worst estimation risk comes from the growing interval $[-A_N, A_N]$ on which $\sigma^2$ is estimated. For the study of the lower bound of the worst estimation risk, the hypotheses do not depend on the estimation interval $[-A_N, A_N]$, but are defined on the real line and are constant outside of the compact interval $[-1,1]$. 

For the case of estimation of $\sigma^2$ on the real line, the obtained lower bound of order $(Nn)^{-\beta/(2\beta+1)}$, with $n \propto N^2$ and $\beta \geq 4$, is faster than the rate provided by Proposition~\ref{prop:Rate-realline}. The sub-optimality of the rate provided by Proposition~\ref{prop:Rate-realline} is a direct consequence of the choice of an approximation space spanned by a basis whose functions are not defined on the real line. Then, we have to deal with the probability of exit of the process $X$ from the growing interval $[-A_N, A_N]$, given by Equation~\eqref{eq:ExitProb}, and which cannot reach a rate of order $(Nn)^{-\beta/(2\beta+1)}$ for $n \propto N$ under the assumptions of the present paper. The Hermite basis could have been suitable for the study of an optimal upper bound of the risk of estimation. Unfortunately, under the assumptions made in this paper, we cannot prove that the Gram matrix $\Psi_m$ related to the Hermite basis satisfies Equation~\eqref{eq:Gram-Condition}. In fact, it is proven in \cite{comte2019regression}, \textit{Proposition 2.5}, that if the density function $f$ of the law that generates the observations satisfies $f(x) \geq c/(1+x^2)^k$ for $x \in \mathbb{R}$ and $k \geq 1$, then the Gram matrix $\Psi_m$ satisfies: $\left\|\Psi_m^{-1}\right\|_{\mathrm{op}} = \mathrm{O}\left(m^k\right)$. In our case, the transition density of the process $X$ is bounded from below by a function that converges exponentially toward zero as $x$ tends to infinity (see Equation~\eqref{eq:approx-gobet}), leading to $\left\|\Psi_m^{-1}\right\|_{\mathrm{op}} = \mathrm{O}\left(\exp(Cm)\right)$, where the constant $C>0$ depends on the diffusion coefficients $b$ and $\sigma$ (see the supplementary material of \cite{ella2024nonparametric}). This situation motivated the construction of a projection estimator of the drift coefficient $b$ on the real line based on the truncation of the dimension proposed in \cite{comte2021drift}. This method is not considered in this paper and may be subject to further investigation.

\section{Conclusion}
\label{sec:conclusion}

In this paper, we have studied minimax convergence rates for minimum contrast estimators of the square $\sigma^2$ of the diffusion coefficient of a time-homogeneous stochastic differential equation. We showed that the projection estimator of $\sigma^2$ on a compact interval and from a single diffusion path reaches a minimax rate of order $n^{-\beta/(2\beta+1)}$. We extended the study to the case where we have $N$ independent observations of the process $X$ and derived a minimax rate of order $(Nn)^{-\beta/(2\beta+1)}$. However, for the case of repeated observations ($N \rightarrow \infty$), we observed that setting a lower bound on the estimation risk of order $(Nn)^{-\beta/(2\beta+1)}$ over the H\"older space requires a strong assumption on the time step $\Delta_n$, that is, $\Delta_n \propto N^{-2}$, whereas in \cite{ella2024nonparametric}, \textit{Corollary 4.2}, the establishment of an upper bound on the estimation risk of order $(Nn)^{-\beta/(2\beta+1)}$ imposes that $\Delta_n \propto N^{-1}$. The assumption $\Delta_n \propto N^{-2}$ for the study of the lower bound of the estimation risk was necessary since the transition density $p_X$ of the process $X$, essential for this study, involves derivatives of the hypotheses whose construction is inspired by a $\beta-$H\"older function provided in \cite{tsybakov2008introduction}, \textit{Chapter 2}. Thus, condition $\Delta_n \propto N^{-1}$ may be suitable for a simpler model such as the following Gaussian model:
\begin{equation*}
    dX_t = b(t)dt + \sigma(t)dW_t, ~~ X_0=x=0,
\end{equation*}
as the transition density of the solution $X$ of the above equation is simpler and does not involve derivatives of the diffusion coefficients.

We faced some challenges when studying the minimax rate for the estimation of $\sigma^2$ on the real line. We showed that the risk of estimation reaches a lower bound of order $(Nn)^{-\beta/(2\beta+1)}$. However, an upper bound of the same order (up to a logarithmic factor) could only be reached on a growing interval of the form $[-A_N, A_N]$, with $A_N \rightarrow \infty$ as $N \rightarrow \infty$. Moreover, setting these two bounds required different assumptions on the time step $\Delta_n$ as in the case of estimation of $\sigma^2$ on a compact interval. A similar result may be obtained on the real line if we consider a minimum contrast estimator of $\sigma^2$ based on the truncation of the dimension of the approximation space, as proposed in \cite{comte2021drift}. 

\section{Proofs}
\label{sec:proofs}

The section is devoted to the proofs of the main results of the paper. We denote the constants by $C>0$, $c>0$ and $C_{\theta}$, where $\theta$ is a constant or a parameter, and these constants can change from one line to another. Finally, for simplicity, we denote the time step $\Delta_n$ by $\Delta$. 

\subsection{Proof of Lemma~\ref{lm:OperatorNorm-Psi}}

\begin{proof}
    As we know from Assumptions~\ref{ass:LipFunctions} and \ref{ass:Ellipticity}, the unique solution $X$ of Equation~\eqref{eq:model} admits a transition density $(s, t, x, y) \mapsto p_X(s, t, x, y)$. From \cite{dacunha1986estimation}, \textit{Lemma 2}, we have
    \begin{multline}\label{eq:TransDens}
        p_X(s, t,x,y) = \dfrac{1}{\sqrt{2\pi(t-s)\sigma^{2}(y)}}\exp\left[-\dfrac{(S(y) - S(x))^2}{2(t-s)} + H(y) - H(x)\right]\\
        \times \mathbf{E}\left[\exp\left((t-s)\int_{0}^{1}{G\left(uS(x) + (1-u)S(y) + \sqrt{t-s}\widetilde{B}_u\right)du}\right)\right],
    \end{multline} 
    where $\widetilde{B} = (\widetilde{B}_t)_{t \in [0,1]}$ is a Brownian bridge, with $\mathbf{E}(\widetilde{B}_t^2) = t(1-t)$, and the functions $S, H$ and $G$ are given by
    \begin{multline}\label{eq:TransitionParam}
        S(x) = \int_{0}^{x}{\dfrac{1}{\sigma(u)}}du, ~~ H(x) = \int_{0}^{S(x)}{\left(\dfrac{b}{\sigma} - \dfrac{\sigma^{\prime}}{2}\right)\circ S^{-1}(u)du}, ~~~ x \in \mathbb{R}, \\
        G = -\dfrac{1}{2}\left[\left(\dfrac{b}{\sigma} - \dfrac{\sigma^{\prime}}{2}\right)^2\circ S^{-1} + \sigma \circ S^{-1} \times \left(\dfrac{b^{\prime} \sigma - b \sigma^{\prime}}{\sigma^{2}} - \dfrac{\sigma^{\prime\prime}}{2}\right) \circ S^{-1}\right],
    \end{multline}
    where $C>0$ is a constant depending on $b$ and $\sigma$. Under Assumption~\ref{ass:Restrict-Model}, the functions $b, ~ b^{\prime}, ~ \sigma, ~ \sigma^{\prime}$ and $\sigma^{\prime\prime}$ are bounded, and for all $x \in \mathbb{R}, ~ G(x) \geq -\|G\|_{\infty}$. We deduce that for all $(s, t, x, y) \in [0,1] \times [0,1] \times \mathbb{R}\times \mathbb{R}$, such that $t > s$, 
    \begin{equation*}
        p_{X}(s, t, x, y) \geq \dfrac{C}{\sqrt{2\pi(t-s)\sigma^{2}(y)}}\exp\left[-\dfrac{(s(y) - s(x))^2}{2(t-s)} + H(y) - H(x)\right].
    \end{equation*}
    Since $X_0 = 0$, we set $x = 0$ at $s=0$ and obtain
    \begin{equation}\label{eq:LowTransDensity}
        p_{X}(0, t, 0, y) \geq \dfrac{C}{\sqrt{2\pi t\sigma^{2}(y)}}\exp\left[-\dfrac{s^2(y)}{2t} + H(y)\right].
    \end{equation}    
    For all $v \in \mathbb{R}^{m}$ such that $v \neq 0$, we have
    \begin{multline*}
        v^{\prime}\mathbf{\Psi}_{m}v = \sum_{\ell = 0}^{m-1}{\sum_{\ell^{\prime} = 0}^{m-1}{\mathbb{E}_X\left[\dfrac{1}{n}\sum_{k=0}^{n-1}{v_{\ell}\phi_{\ell}(X_{k\Delta})v_{\ell^{\prime}}\phi_{\ell^{\prime}}(X_{k\Delta})}\right]}}\\
        = \mathbb{E}_X\left[\dfrac{1}{n}\sum_{k=0}^{n-1}{\left(\sum_{\ell = 0}^{m-1}{v_{\ell}\phi_{\ell}(X_{k\Delta})}\right)^2}\right] = \|h_v\|_{n}^2,
    \end{multline*}
    where $h_v = \sum_{\ell = 0}^{m-1}{v_{\ell}\phi_{\ell}} \in \mathcal{S}_{m, A, B}$. We obtain
    \begin{equation}\label{eq:MinEigen}
        v^{\prime}\mathbf{\Psi}_{m}v = \|h_v\|_{n}^2 \geq \int_{-A_{N}}^{A_{N}}{h_v^2(y)f_{n}(y)dy},
    \end{equation}
    where $f_{n}: y \in \mathbb{R} \mapsto \dfrac{1}{n}\sum_{k=1}^{n-1}{p_{X}(0, k\Delta, 0, y)}$.
    
    Under Assumption~\ref{ass:Restrict-Model}, the functions $b, \sigma, \sigma^{\prime}$ are bounded, then $H$ is bounded on the interval $[-A_N, A_N]$, and there exists a constant $c>0$ such that for all $x \in [-A_N, A_N]$, we have
    \begin{equation}\label{eq:LowH}
        |H(y)| = \left|\int_{0}^{S(y)}{\left(\dfrac{b}{\sigma} - \dfrac{1}{\sigma^{\prime}}\right)\circ S^{-1}(u)du}\right| \leq c \kappa_0 |S(y)| \leq cA_N,
    \end{equation}
    which implies that for all $y \in \mathbb{R}$, $H(y) \geq -cA_N$. We then obtain from Equations~\eqref{eq:LowH} and \eqref{eq:LowTransDensity},
    \begin{equation*}
        p_{X}(0, t, 0, y) \geq \dfrac{C}{\sqrt{t}}\exp\left(-\dfrac{S^{2}(y)}{2t} + H(y)\right) \geq C\exp\left(-\dfrac{S^{2}(y)}{2t} - cA_{N}\right), ~~ (t, y) \in (0,1] \times [-A_{N}, A_{N}].
    \end{equation*}
    Moreover, for all $y \in [-A_N, A_N]$ such that $y \neq 0$, the function $t \mapsto C\exp\left(-\dfrac{S^{2}(y)}{2t} - cA_{N}\right)$ is increasing on the interval $(0,1]$ since for all $y \in [-A_N, A_N]$ such that $y \neq 0$ and for all $t \in (0,1]$,
    \begin{equation*}
        \dfrac{\partial}{\partial t}\left[C\exp\left(-\dfrac{S^{2}(y)}{2t} - cA_{N}\right)\right] = \dfrac{S^2(y)}{2t^2}\exp\left(-\dfrac{S^2(y)}{2t} - cA_N\right) > 0.
    \end{equation*}
    Then, setting $\xi(s) = k\Delta$ for all $s \in [k\Delta, (k+1)\Delta)$, we obtain for all $y \in [-A_N, A_N]$ such that $y \neq 0$,
    \begin{multline*}
        f_{n}(y) = \int_{\Delta}^{(n-1)\Delta}{p_{X}(0, \xi(s), 0, y)ds} \geq C\int_{\Delta}^{(n-1)\Delta}{\exp\left(-\dfrac{S^2(y)}{2\xi(s)} - cA_N\right)ds} \\ 
        \geq C\int_{\Delta}^{(n-1)\Delta}{\exp\left(-\dfrac{S^2(y)}{2(s-\Delta)} - cA_N\right)ds}\\
        = C\int_{0}^{(n-2)\Delta}{\exp\left(-\dfrac{S^2(y)}{2s} - cA_N\right)ds}.
    \end{multline*}
    For $n \geq N$ and $N$ large enough,
    \begin{multline}\label{eq:LowerBoundDensity}
        f_{n}(y) \geq C\int_{1-\log^{-1}(N)}^{1-2^{-1}\log^{-1}(N)}{\exp\left(-\dfrac{S^2(y)}{2s} - cA_N\right)ds} \\
        \geq \dfrac{C}{\log(N)}\exp\left(-\dfrac{S^2(A_{N})}{2(1 - \log^{-1}(N))} - cA_{N}\right).
    \end{multline}
    From Equations~\eqref{eq:MinEigen} and \eqref{eq:LowerBoundDensity}, we obtain
    \begin{equation*}
        v^{\prime}\mathbf{\Psi}_{m}v \geq \dfrac{C}{\log(N)}\exp\left(-\dfrac{1}{2(1 - \log^{-1}(N))}\left(\int_{0}^{A_{N}}{\dfrac{1}{\sigma(u)}du}\right)^2 - cA_{N}\right)\|h_v\|^2.
    \end{equation*}
    We now focus on the lower bounding of 
    $$\|h_v\|^2 = \int_{-A_N}^{A_N}{\left(\sum_{\ell = 0}^{m-1}{v_{\ell}\phi_{\ell}(x)}\right)^2dx}.$$
    The result depends on the basis $\left(\phi_0, \ldots, \phi_{m-1}\right)$ that generated the approximation space $\mathcal{S}_{m, A, B}$. Then, for the spline basis ($[\mathbf{B}]$), we have
    \begin{equation}\label{eq:bound-spline}
        \mathcal{L}(m) := \underset{x \in [-A_N, A_N]}{\sup}{\sum_{\ell = 0}^{m-1}{\phi_{\ell}^2(x)}} = \underset{x \in \mathbb{R}}{\sup}{\sum_{\ell = -M}^{K-1}{B_{\ell}^2(x)}} \leq 1,
    \end{equation}
    and from \cite{denis2020ridge}, \textit{Lemma 2.6}, there exists a constant $C>0$ such that $\|h_v\|^2 \geq \dfrac{CA_{N}}{m}\|v\|_2^2$. Then, there exists a constant $C>0$ such that
    \begin{equation}\label{eq:eig-spline}
        v^{\prime}\mathbf{\Psi}_{m}v \geq \dfrac{CA_{N}}{m\log(N)}\exp\left(-\dfrac{1}{2(1 - \log^{-1}(N))}\left(\int_{0}^{A_{N}}{\dfrac{1}{\sigma(u)}du}\right)^2  - cA_{N}\right)\|v\|_{2}^2.
    \end{equation}
    Then, for all $v \in \mathbb{R}^m$ such that $v \neq 0, ~~ v^{\prime}\Psi_{m}v > 0$, which implies that $\Psi_{m}$ is invertible, being a symmetrical and positive definite matrix, and we obtain from Equations~\eqref{eq:bound-spline} and \eqref{eq:eig-spline} that
    \begin{equation*}
        \mathcal{L}(m)\left\|\mathbf{\Psi}^{-1}_{m}\right\|_{\mathrm{op}} \leq \dfrac{m\log(N)}{CA_{N}}\exp\left(\dfrac{1}{2(1 - \log^{-1}(N))}\left(\int_{0}^{A_{N}}{\dfrac{1}{\sigma(u)}du}\right)^2  + cA_{N}\right).
    \end{equation*}
    For the Fourier basis ($[\mathbf{F}]$), we have the following:
    \begin{equation*}
        \mathcal{L}(m) := \underset{x \in [-A_N, A_N]}{\sup}{\sum_{\ell = 0}^{m-1}{\phi_{\ell}^2(x)}} \leq \left(\underset{\ell = 0, \ldots, m-1}{\max}{\|\phi_{\ell}\|_{\infty}^2}\right)m, ~~ \mathrm{and} ~~ \|h_v\|^2 = \|v\|_2^2.
    \end{equation*}
    Then, $\Psi_m$ is also a symmetrical and positive definite matrix since, for any $v \in \mathbb{R}^m$ such that $v \neq 0$, we have 
    \begin{equation*}
        v^{\prime}\mathbf{\Psi}_{m}v \geq \dfrac{CA_{N}}{\log(N)}\exp\left(-\dfrac{1}{2(1 - \log^{-1}(N))}\left(\int_{0}^{A_{N}}{\dfrac{1}{\sigma(u)}du}\right)^2  - cA_{N}\right)\|v\|_2^2 > 0,
    \end{equation*}
    and there exists a constant $C>0$ such that
    \begin{equation*}
        \mathcal{L}(m)\left\|\mathbf{\Psi}^{-1}_{m}\right\|_{\mathrm{op}} \leq \dfrac{m\log(N)}{CA_{N}}\exp\left(\dfrac{1}{2(1 - \log^{-1}(N))}\left(\int_{0}^{A_{N}}{\dfrac{1}{\sigma(u)}du}\right)^2  + cA_{N}\right).
    \end{equation*}
\end{proof}

\subsection{Proof of Theorem~\ref{thm:upper-bound}}

\begin{proof} 
Let $\Omega_{n,N,K_N}$ be the random event in which the empirical norm $\|.\|_n$ and the pseudo-norm $\|.\|_{n,N}$ are equivalent and given as follows:
\begin{equation*}
    \Omega_{n,N,K_N} := \left\{h \in \mathcal{S}_{K_N} \setminus \{0\}: \left|\dfrac{\|h\|_{n,N}^2}{\|h\|_{n}^2} - 1\right| \leq \dfrac{1}{2}\right\}.
\end{equation*}
On the event $\Omega_{n,N,K_N}$, we have for all $h \in \mathcal{S}_{K_N}$, $\dfrac{1}{2}\|h\|_n^2 \leq \|h\|_{n,N}^2 \leq \dfrac{3}{2}\|h\|_n^2$. Under the assumptions of the theorem, we obtain from \cite{denis2024nonparametric}, \textit{proof of Lemma 7, Equation (A13)} that
\begin{equation}\label{eq:Prob-OmegaComp}
    \mathbb{P}\left(\Omega^{c}_{n,N,K_N}\right) \leq 2K_N\exp\left(-C_1\log^{3/2}(N)\right).
\end{equation}
Recall that for all $j \in \{1,\ldots, N\}$,
$\zeta^{j}_{k\Delta}=\zeta^{j,1}_{k\Delta}+\zeta^{j,2}_{k\Delta}+\zeta^{j,3}_{k\Delta}$ is the error term of the regression model, with:
\begin{equation}
\label{eq:zeta 1}
    \zeta^{j,1}_{k\Delta}=\frac{1}{\Delta}\left[\left(\int_{k\Delta}^{(k+1)\Delta}{\sigma(X^{j}_{s})dW^{j}_{s}}\right)^2-\int_{k\Delta}^{(k+1)\Delta}{\sigma^{2}(X^{j}_{s})ds}\right],
\end{equation}
\begin{equation}
\label{eq:zeta 2}
    \zeta^{j,2}_{k\Delta}=\frac{2}{\Delta}\int_{k\Delta}^{(k+1)\Delta}{((k+1)\Delta-s)\sigma^{\prime}(X^{j}_{s})\sigma^{2}(X^{j}_{s})dW^{j}_{s}},
\end{equation}
\begin{equation}
\label{eq:zeta 3}
\zeta^{j,3}_{k\Delta}=2b(X^{j}_{k\Delta})\int_{k\Delta}^{(k+1)\Delta}{\sigma\left(X^{j}_{s}\right)dW^{j}_{s}}.
\end{equation}
In addition, $R^{j}_{k\Delta}
=R^{j,1}_{k\Delta}+R^{j,2}_{k\Delta} + R^{j,3}_{k\Delta}$ is a negligible residual with:
\begin{equation}
\label{eq: R 1}
    R^{j,1}_{k\Delta}=\frac{1}{\Delta}\left(\int_{k\Delta}^{(k+1)\Delta}{b(X^{j}_{s})ds}\right)^2, ~~ R^{j,2}_{k\Delta} = \frac{1}{\Delta}\int_{k\Delta}^{(k+1)\Delta}{((k+1)\Delta-s)\Phi(X^{j}_{s})ds}
\end{equation}
\begin{equation}
\label{eq: R 2}
R^{j,3}_{k\Delta}=\frac{2}{\Delta}\left(\int_{k\Delta}^{(k+1)\Delta}{\left(b(X^{j}_{s})-b(X^{j}_{k\Delta})\right)ds}\right)\left(\int_{k\Delta}^{(k+1)\Delta}{\sigma(X^{j}_{s})dW^{j}_{s}}\right),
\end{equation}
where 
\begin{equation*}
\Phi:=2b\sigma^{\prime}\sigma+\left[\sigma^{\prime\prime}\sigma+\left(\sigma^{\prime}\right)^2\right]\sigma^{2}.
\end{equation*}
By the definition of the projection estimator $\widehat{\sigma}^{2}_{K_N}$ with $K_N = N^{2/(2\beta+1)}\log^{-5/2}(N)$, for all $h \in \mathcal{S}_{K_N, A_N}$, we have the following:
\begin{equation}
\label{eq:Ineq-Gamma}
    \gamma_{n,N}\left(\widehat{\sigma}^{2}_{K_N}\right)-\gamma_{n,N}(\sigma^{2}_{A_N})\leq\gamma_{n,N}(h)-\gamma_{n,N}(\sigma^{2}_{A_N}).
\end{equation}
Furthermore, for all $h \in \mathcal{S}_{K_N, A_N}$, 
$$\gamma_{n,N}(h)-\gamma_{n,N}(\sigma^{2}_{A_N})=\left\|\sigma^{2}_{A_N}-h\right\|^{2}_{n,N}+2\nu_{1}(\sigma^{2}_{A_N}-h)+2\nu_{2}(\sigma^{2}_{A_N}-h)+2\nu_{3}(\sigma^{2}_{A_N}-h)+2\mu(\sigma^{2}_{A_N}-h),$$
where 
\begin{equation}
    \label{eq: nu and mu}
    \nu_{i}\left(h\right) = \frac{1}{n}\sum_{j=1}^{N}{\sum_{k=0}^{n-1}{h(X^{j}_{k\Delta})\zeta^{j,i}_{k\Delta}}}, \ \ i\in\{1,2,3\}, \ \ \ \ \mu(h)=\frac{1}{n}\sum_{j=1}^{N}{\sum_{k=0}^{n-1}{h(X^{j}_{k\Delta})R^{j}_{k\Delta}}}, 
\end{equation}
and $\zeta^{j,1}_{k\Delta}, \ \zeta^{j,2}_{k\Delta}, \ \zeta^{j,3}_{k\Delta}$ are given by Equations \eqref{eq:zeta 1}, \eqref{eq:zeta 2}, \eqref{eq:zeta 3}, and finally, $R^{j}_{k\Delta} = R^{j,1}_{k\Delta} + R^{j,2}_{k\Delta} + R^{j,3}_{k\Delta}$ given in Equations \eqref{eq: R 1} and \eqref{eq: R 2}. Then, for all $h \in \mathcal{S}_{K_N, A_N}$, we obtain from Equation~\eqref{eq:Ineq-Gamma} that
\begin{equation*}
\left\|\widehat{\sigma}^{2}_{K_N}-\sigma^{2}_{A_N}\right\|^{2}_{n,N}\leq \left\|h-\sigma^{2}_{A_N}\right\|^{2}_{n,N}+2\nu\left(\widehat{\sigma}^{2}_{K_N}-h\right)+2\mu\left(\widehat{\sigma}^{2}_{K_N}-h\right), \ \ \mathrm{with} \ \ \nu=\nu_1+\nu_2+\nu_3.
\end{equation*}
Then
\begin{equation}
    \label{eq:equation2-proof1}
    \mathbb{E}\left[\left\|\widehat{\sigma}^{2}_{K_N}-\sigma^{2}_{A_N}\right\|^{2}_{n,N}\right] \leq \underset{h\in\mathcal{S}_{K_N, A_N}}{\inf}{\left\|h-\sigma^{2}_{A_N}\right\|^{2}_{n}}+2\mathbb{E}\left[\nu\left(\widehat{\sigma}^{2}_{K_N}-h\right)\right]+2\mathbb{E}\left[\mu\left(\widehat{\sigma}^{2}_{K_N}-h\right)\right].
\end{equation}
Besides, for any $a,d>0$, using the inequality $2xy \leq \eta x^2 + y^2/\eta$ with $\eta = a, d$, we obtain from Equation \eqref{eq: nu and mu} and on the event $\Omega_{n,N,K_N}$,
\begin{multline}\label{eq:mj-nu}
    2\nu\left(\widehat{\sigma}^{2}_{K_N}-h\right) = 2\left\|\widehat{\sigma}^{2}_{K_N}-h\right\|_n\nu\left(\dfrac{\widehat{\sigma}^{2}_{K_N}-h}{\left\|\widehat{\sigma}^{2}_{K_N}-h\right\|_n}\right)\\
    \leq \dfrac{1}{a}\left\|\widehat{\sigma}^{2}_{K_N}-h\right\|_n^2 + a\underset{h\in\mathcal{S}_{K_N}, \ \|h\|_{n}=1}{\sup}{\nu^{2}(h)}\\
    \leq \dfrac{2}{a}\left\|\widehat{\sigma}^{2}_{K_N} - \sigma_{A_N}^2\right\|_{n,N}^2 + \dfrac{2}{a}\left\|h - \sigma_{A_N}^2\right\|_{n,N}^2 + a\underset{h\in\mathcal{S}_{K_N}, \ \|h\|_{n}=1}{\sup}{\nu^{2}(h)},
\end{multline}
\begin{multline}\label{eq:mj-mu}
    2\mu\left(\widehat{\sigma}^{2}_{K_N}-h\right) \leq \dfrac{1}{d}\left\|\widehat{\sigma}^{2}_{K_N}-h\right\|_n^2 + \dfrac{d}{Nn}\sum_{j=1}^{N}{\sum_{k=0}^{n-1}{\left(R_{k\Delta}^j\right)^2}}\\
    \leq \dfrac{2}{d}\left\|\widehat{\sigma}^{2}_{K_N}-\sigma^{2}_{A_N}\right\|^{2}_{n,N} + \dfrac{2}{d}\left\|h-\sigma^{2}_{A_N}\right\|^{2}_{n,N} + \dfrac{d}{Nn}\sum_{j=1}^{N}{\sum_{k=0}^{n-1}{(R^{j}_{k\Delta})^2}}.
\end{multline}
From \cite{ella2024nonparametric}, \textit{Proof of Theorem 3.2}, there exists a constant $C>0$ such that
\begin{equation}
\label{eq:UpperBound-TimeStep}
\mathbb{E}\left[\frac{1}{n}\sum_{k=1}^{n}{(R^{j}_{k\Delta})^2}\right]\leq C\Delta^2.
\end{equation}
We set $a = 8$ and $d = 8$, and, using the following relation 
\begin{align*}
    \mathbb{E}\left[\left\|\widehat{\sigma}^{2}_{K_N}-\sigma^{2}_{A_N}\right\|^{2}_{n,N}\right] = &~ \mathbb{E}\left[\left\|\widehat{\sigma}^{2}_{K_N}-\sigma^{2}_{A_N}\right\|^{2}_{n,N}\mathds{1}_{\Omega_{n,N,K_N}}\right] + \mathbb{E}\left[\left\|\widehat{\sigma}^{2}_{K_N}-\sigma^{2}_{A_N}\right\|^{2}_{N,n}\mathds{1}_{\Omega^{c}_{n,N,K_N}}\right] \\
    \leq &~ \mathbb{E}\left[\left\|\widehat{\sigma}^{2}_{K_N}-\sigma^{2}_{A_N}\right\|^{2}_{n,N}\mathds{1}_{\Omega_{n,N,K_N}}\right] + 2K_NA_N^2\log(N)\mathbb{P}\left(\Omega^{c}_{n,N,K_N}\right),
\end{align*}
we deduce from Equations~\eqref{eq:Bias-Term}, \eqref{eq:Prob-OmegaComp}, \eqref{eq:equation2-proof1}, \eqref{eq:mj-nu}, \eqref{eq:mj-mu} and \eqref{eq:UpperBound-TimeStep} and for $N$ large enough that,
\begin{multline}
\label{eq:equation3-proof1}
	\mathbb{E}\left[\left\|\widehat{\sigma}^{2}_{K_N}-\sigma^{2}_{A_N}\right\|^{2}_{n,N}\right] \leq 3\underset{h\in\mathcal{S}_{K_N, A_N}}{\inf}{\left\|h-\sigma^{2}_{A_N}\right\|^{2}_{n}} + C\mathbb{E}\left(\underset{h \in \mathcal{S}_{K_N}, ~ \|h\|_{n} = 1}{\sup}{\nu^{2}(h)}\right) + C\Delta^2\\
    + 4A_N^2K_N\log(N)\exp\left(-C_1\log^{3/2}(N)\right)\\
      \leq C\left[\left(\dfrac{A_N}{K_N}\right)^{2\beta} + \mathbb{E}\left(\underset{h\in\mathcal{S}_{K_N}, ~ \|h\|_{n} = 1}{\sup}{\nu^{2}(h)}\right) + \Delta^2 + \exp\left(-C_2\log^{3/2}(N)\right)\right],
\end{multline}
where $C>0$ is a constant depending on $\kappa_1$, and $C_2 > 0$ is a constant depending on $C_1$.

\subsection*{Upper bound of $\mathbb{E}\left(\underset{h\in\mathcal{S}_{K_N}, \ \|h\|_{n}^2 = 1}{\sup}{\nu^{2}(h)}\right)$}

For all $h = \sum_{\ell = -M}^{K_N-1}{a_{\ell}B_{\ell}} \in \mathcal{S}_{K_N}$, we have $\|h\|^{2}_{n} = {^t}\mathbf{a}\mathbf{\Psi}_{K_N}\mathbf{a} = \left\|\mathbf{\Psi}^{1/2}_{K_N}\mathbf{a}\right\|^{2}_{2,K_N+M}$. Then,
$$\|h\|^{2}_{n} = 1 \iff \exists ~ \mathbf{u} \in \mathbb{R}^{K_N}: ~ \mathbf{a} = \mathbf{\Psi}^{-1/2}_{K_N}\mathbf{u} ~~ \mathrm{and} ~~ \|\mathbf{u}\|_{2,K_N+M} = 1.$$
Thus, for all $h \in \mathcal{S}_{K_N}$ such that $\|h\|^{2}_{n} = 1$, there exists a vector $\mathbf{u} \in \mathbb{R}^{K_N+M}$ such that $\|\mathbf{u}\|_{2,K_N+M} = 1$ and $\mathbf{\Psi}^{-1/2}_{K_N}\mathbf{u}$ is the coordinate vector of $h$, and we have
\begin{align*}
    \nu^{2}(h) = &~ \left(\dfrac{1}{Nn}\sum_{j=1}^{N}{\sum_{k=0}^{n-1}{\sum_{\ell = -M}^{K_N-1}{\left[\mathbf{\Psi}^{-1/2}_{K_N}\mathbf{u}\right]_{\ell}B_{\ell}(X^{j}_{k\Delta})\zeta^{j}_{k\Delta}}}}\right)^2\\ 
    = &~ \dfrac{1}{N^2n^2}\left(\sum_{\ell^{\prime} = -M}^{K_N-1}{u_{\ell^{\prime}}\sum_{j=1}^{N}{\sum_{k=0}^{n-1}{\sum_{\ell = -M}^{K_N-1}{\left[\mathbf{\Psi}^{-1/2}_{K_N}\right]_{\ell,\ell^{\prime}}}B_{\ell}(X^{j}_{k\Delta})\zeta^{j}_{k\Delta}}}}\right)^2,
\end{align*}
and using the Cauchy-Schwarz inequality and the fact that $\|\mathbf{u}\|_{2,K_N+M} = 1$, we obtain
\begin{equation}\label{eq:nu-bound1}
    \nu^{2}(h) \leq \dfrac{1}{N^2n^2}\sum_{\ell^{\prime} = -M}^{K_N-1}{\left(\sum_{j=1}^{N}{\sum_{k=0}^{n-1}{\sum_{\ell = -M}^{K_N-1}{\left[\mathbf{\Psi}^{-1/2}_{K_N}\right]_{\ell,\ell^{\prime}}}B_{\ell}(X^{j}_{k\Delta})\zeta^{j}_{k\Delta}}}\right)^2}.
\end{equation}
Recall that for all $j = 1, \ldots, N$ and for all $k = 0, \ldots, n-1$, we have $\zeta^{j}_{k\Delta} = \zeta^{j,1}_{k\Delta} + \zeta^{j,2}_{k\Delta} + \zeta^{j,3}_{k\Delta}$. 
From Equation~\eqref{eq:zeta 1}, we have $\zeta^{j,1}_{k\Delta}=\frac{1}{\Delta}\left[\left(\int_{k\Delta}^{(k+1)\Delta}{\sigma(X^{j}_{s})}dW_s^j\right)^2-\int_{k\Delta}^{(k+1)\Delta}{\sigma^{2}(X^{j}_{s})ds}\right]$.
     We fix an initial time $s\in[0,1)$ and set $M^{s}_t=\int_{s}^{t}{\sigma(X^{j}_u)dW_u^j}, \ \forall t\geq s$. $(M^{s}_t)_{t\geq s}$ is a martingale with respect to the natural filtration $\left(\mathcal{F}_{t \in [0,1]}\right)$ of the process $X$, and for all $t\in[s,1]$, we have:  
    \begin{align*}
      \left<M^{s},M^{s}\right>_t=\int_{s}^{t}{\sigma^{2}\left(X^{j}_{u}\right)du}.
    \end{align*}
     Then, $\zeta^{j,1}_{k\Delta}=\frac{1}{\Delta}\left(M^{k\Delta}_{(k+1)\Delta}\right)^2-\left<M^{k\Delta},M^{k\Delta}\right>_{(k+1)\Delta}$ is also a $\mathcal{F}_{k\Delta}$-martingale, and, using the Burkholder-Davis-Gundy inequality, we obtain for all $k\in[\![0,n-1]\!]$,
   \begin{equation}
     \label{eq:martingale-burkholder-davis-gundy}
      \mathbb{E}\left[\zeta^{j,1}_{k\Delta}|\mathcal{F}_{k\Delta}\right]=0, \ \ \ \mathbb{E}\left[\left(\zeta^{j,1}_{k\Delta}\right)^2|\mathcal{F}_{k\Delta}\right]\leq \frac{C}{\Delta^2}\mathbb{E}\left[\left(\int_{k\Delta}^{(k+1)\Delta}{\sigma^{2}(X^{j}_u)du}\right)^2\right]\leq C\kappa^{4}_{1}.
   \end{equation}
We also have from Equations~\eqref{eq:zeta 2} and \eqref{eq:zeta 3} that
\begin{equation}\label{eq:Exp-zero}
    \forall ~ k \in \{0, \ldots, n-1\}, ~~ \mathbb{E}\left[\zeta^{j,i}_{k\Delta} | \mathcal{F}_{k\Delta}\right] = 0, ~~ i = 2, 3. 
\end{equation}
Then, we obtain for all $k \in \{0,\ldots,n-1\}$, $\mathbb{E}\left[\zeta^{j}_{k\Delta} | \mathcal{F}_{k\Delta}\right] = 0$, and using the independence between the diffusion paths $X^j, ~~ j = 1, \ldots, N$, we have:
\begin{multline}\label{eq:nu-bound2}
    \mathbb{E}\left[\left(\sum_{j=1}^{N}{\sum_{k=0}^{n-1}{\sum_{\ell = -M}^{K_N-1}{\left[\mathbf{\Psi}^{-1/2}_{K_N}\right]_{\ell,\ell^{\prime}}}B_{\ell}(X^{j}_{k\Delta})\zeta^{j}_{k\Delta}}}\right)^2\right] = \mathrm{Var}\left(\sum_{j=1}^{N}{\sum_{k=0}^{n-1}{\sum_{\ell = -M}^{K_N-1}{\left[\mathbf{\Psi}^{-1/2}_{K_N}\right]_{\ell,\ell^{\prime}}}B_{\ell}(X^{j}_{k\Delta})\zeta^{j}_{k\Delta}}}\right)\\
    = \sum_{j=1}^{N}{\mathrm{Var}\left(\sum_{k=0}^{n-1}{\sum_{\ell = -M}^{K_N-1}{\left[\mathbf{\Psi}^{-1/2}_{K_N}\right]_{\ell,\ell^{\prime}}}B_{\ell}(X^{j}_{k\Delta})\zeta^{j}_{k\Delta}}\right)}\\
    = \sum_{j=1}^{N}{\mathbb{E}\left[\left(\sum_{k=0}^{n-1}{\sum_{\ell = -M}^{K_N-1}{\left[\mathbf{\Psi}^{-1/2}_{K_N}\right]_{\ell,\ell^{\prime}}}B_{\ell}(X^{j}_{k\Delta})\zeta^{j}_{k\Delta}}\right)^2\right]}.
\end{multline}
We deduce from Equations \eqref{eq:nu-bound1} and \eqref{eq:nu-bound2},
\begin{multline*}
    \mathbb{E}\left(\underset{h\in\mathcal{S}_{K_N}, \ \|h\|_{n}=1}{\sup}{\nu^{2}(h)}\right) \leq \dfrac{1}{N^2n^2}\sum_{\ell^{\prime} = -M}^{K_N-1}{\mathbb{E}\left[\left(\sum_{j=1}^{N}{\sum_{k=0}^{n-1}{\sum_{\ell = -M}^{K_N-1}{\left[\mathbf{\Psi}^{-1/2}_{K_N}\right]_{\ell,\ell^{\prime}}}B_{\ell}(X^{j}_{k\Delta})\zeta^{j}_{k\Delta}}}\right)^2\right]}\\
    \leq \dfrac{1}{N^2n^2}\sum_{\ell^{\prime} = -M}^{K_N-1}{\sum_{j=1}^{N}{\mathbb{E}\left[\left(\sum_{k=0}^{n-1}{\sum_{\ell = -M}^{K_N-1}{\left[\mathbf{\Psi}^{-1/2}_{K_N}\right]_{\ell,\ell^{\prime}}}B_{\ell}(X^{j}_{k\Delta})\zeta^{j}_{k\Delta}}\right)^2\right]}}\\
    \leq \dfrac{1}{N^2n^2}\sum_{j=1}^{N}{\sum_{\ell^{\prime} = -M}^{K_N-1}{\sum_{\ell = -M}^{K_N-1}{\sum_{\ell^{\prime\prime} = -M}^{K_N-1}{\left[\mathbf{\Psi}^{-1/2}_{K_N}\right]_{\ell,\ell^{\prime}}\left[\mathbf{\Psi}^{-1/2}_{K_N}\right]_{\ell^{\prime\prime},\ell^{\prime}}}\mathbb{E}\left[\sum_{k,k^{\prime} = 0}^{n-1}{B_{\ell}(X^{j}_{k\Delta})B_{\ell^{\prime\prime}}(X^{j}_{k^{\prime}\Delta})\zeta^{j}_{k\Delta}\zeta^{j}_{k^{\prime}\Delta}}\right]}}}.
\end{multline*}
Then
\begin{multline}\label{eq:mj-nu3}
   \mathbb{E}\left(\underset{h\in\mathcal{S}_{K_N}, \ \|h\|_{n}=1}{\sup}{\nu^{2}(h)}\right) \leq \dfrac{1}{N^2n^2}\sum_{j=1}^{N}{\sum_{\ell^{\prime} = -M}^{K_N-1}{\sum_{\ell = -M}^{K_N-1}{\sum_{\ell^{\prime\prime} = -M}^{K_N-1}{\left[\mathbf{\Psi}^{-1/2}_{K_N}\right]_{\ell,\ell^{\prime}}\left[\mathbf{\Psi}^{-1/2}_{K_N}\right]_{\ell^{\prime\prime},\ell^{\prime}} T_1^j(\ell, \ell^{\prime\prime})}}}}\\
     + \dfrac{1}{N^2n^2}\sum_{j=1}^{N}{\sum_{\ell^{\prime} = -M}^{K_N-1}{\sum_{\ell = -M}^{K_N-1}{\sum_{\ell^{\prime\prime} = -M}^{K_N-1}{\left[\mathbf{\Psi}^{-1/2}_{K_N}\right]_{\ell,\ell^{\prime}}\left[\mathbf{\Psi}^{-1/2}_{K_N}\right]_{\ell^{\prime\prime},\ell^{\prime}} T_2^j(\ell, \ell^{\prime\prime})}}}}, 
\end{multline}
where, for all $j \in \{0,\ldots, N\}$ and for all $\ell, \ell^{\prime\prime} \in \{0, m-1\}$, we set 
$$T_1^j(\ell, \ell^{\prime\prime}) = \mathbb{E}\left[\sum_{k = 0}^{n-1}{B_{\ell}(X^{j}_{k\Delta})B_{\ell^{\prime\prime}}(X^{j}_{k\Delta})\left(\zeta^{j}_{k\Delta}\right)^2}\right], ~~~ T_2^j(\ell, \ell^{\prime\prime}) = \mathbb{E}\left[\sum_{k \neq k^{\prime}}{B_{\ell}(X^{j}_{k\Delta})B_{\ell^{\prime\prime}}(X^{j}_{k^{\prime}\Delta})\zeta^{j}_{k\Delta}\zeta^{j}_{k^{\prime}\Delta}}\right].$$

\subsubsection*{Upper-bound of $T_1^j(\ell, \ell^{\prime\prime})$}

We have the following:
\begin{equation}\label{eq:T1}
    T_1^j(\ell, \ell^{\prime\prime}) \leq T_{1,1}^j(\ell, \ell^{\prime\prime}) + T_{1,2}^j(\ell, \ell^{\prime\prime}) + T_{1,3}^j(\ell, \ell^{\prime\prime}),
\end{equation}
where 
\begin{align*}
    T_{1,1}^j(\ell, \ell^{\prime\prime}) = &~ 3\mathbb{E}\left[\sum_{k = 0}^{n-1}{B_{\ell}(X^{j}_{k\Delta})B_{\ell^{\prime\prime}}(X^{j}_{k\Delta})\left(\zeta^{j,1}_{k\Delta}\right)^2}\right], ~~~ T_{1,2}^j(\ell, \ell^{\prime\prime}) = 3\mathbb{E}\left[\sum_{k = 0}^{n-1}{B_{\ell}(X^{j}_{k\Delta})B_{\ell^{\prime\prime}}(X^{j}_{k\Delta})\left(\zeta^{j,2}_{k\Delta}\right)^2}\right]\\
    T_{1,3}^j(\ell, \ell^{\prime\prime}) = &~ 3\mathbb{E}\left[\sum_{k = 0}^{n-1}{B_{\ell}(X^{j}_{k\Delta})B_{\ell^{\prime\prime}}(X^{j}_{k\Delta})\left(\zeta^{j,3}_{k\Delta}\right)^2}\right],
\end{align*}
and $\zeta^{j,1}_{k\Delta}, ~ \zeta^{j,2}_{k\Delta}$ and $\zeta^{j,3}_{k\Delta}$ are given in Equations~\eqref{eq:zeta 1}, \eqref{eq:zeta 2} and \eqref{eq:zeta 3}.

\begin{enumerate}
    \item Upper-bound of $T_{1,1}^j(\ell, \ell^{\prime\prime})$
    
Using Equation~\eqref{eq:martingale-burkholder-davis-gundy} we have:
\begin{multline*}
     T_{1,1}^j(\ell, \ell^{\prime\prime}) = 3\mathbb{E}\left[\sum_{k = 0}^{n-1}{B_{\ell}(X^{j}_{k\Delta})B_{\ell^{\prime\prime}}(X^{j}_{k\Delta})\left(\zeta^{j,1}_{k\Delta}\right)^2}\right]\\
     = 3\mathbb{E}\left[\sum_{k = 0}^{n-1}{B_{\ell}(X^{j}_{k\Delta})B_{\ell^{\prime\prime}}(X^{j}_{k\Delta})\mathbb{E}\left[\left(\zeta^{j,1}_{k\Delta}\right)^2|\mathcal{F}_{k\Delta}\right]}\right]\\
     \leq Cn\left[\mathbf{\Psi}_{K_N}\right]_{\ell, \ell^{\prime\prime}},
\end{multline*}
where the constant $C>0$ depends on $\sigma$.

\item Upper-bound of $T_{1,2}^j(\ell, \ell^{\prime\prime})$

From Equation~\eqref{eq:zeta 2}, we have $\zeta^{j,2}_{k\Delta}=\frac{2}{\Delta}\int_{k\Delta}^{(k+1)\Delta}{\left((k+1)\Delta-s\right)\sigma^{\prime}(X^{j}_{s})\sigma^{2}(X^{j}_{s})dW_s}$ and,
\begin{multline*}
    T_{1,2}^j(\ell, \ell^{\prime\prime}) =  3\mathbb{E}\left[\sum_{k = 0}^{n-1}{B_{\ell}(X^{j}_{k\Delta})B_{\ell^{\prime\prime}}(X^{j}_{k\Delta})\left(\zeta^{j,2}_{k\Delta}\right)^2}\right]\\
    = \dfrac{12}{\Delta^2}\mathbb{E}\left[\sum_{k=0}^{n-1}{B_{\ell}(X^{j}_{k\Delta})B_{\ell^{\prime\prime}}(X^{j}_{k\Delta})\mathbb{E}\left[\left(\int_{k\Delta}^{(k+1)\Delta}{(k+1)\Delta-s)}\sigma^{\prime}\left(X^{j}_{s}\right)\sigma^{2}\left(X^{j}_{s}\right)dW_s^j\right)^2|\mathcal{F}_{k\Delta}\right]}\right]\\
    = \dfrac{12}{\Delta^2}\mathbb{E}\left[\sum_{k=0}^{n-1}{B_{\ell}(X^{j}_{k\Delta})B_{\ell^{\prime\prime}}(X^{j}_{k\Delta})\mathbb{E}\left[\int_{k\Delta}^{(k+1)\Delta}{(k+1)\Delta-s)^2}\sigma^{\prime 2}\left(X^{j}_{s}\right)\sigma^{4}\left(X^{j}_{s}\right)ds|\mathcal{F}_{k\Delta}\right]}\right]\\
    \leq C\left[\mathbf{\Psi}_{K_N}\right]_{\ell, \ell^{\prime\prime}},
\end{multline*}
where the constant $C>0$ depends on $\sigma$.

\item Upper-bound of $T_{1,3}^j(\ell, \ell^{\prime\prime})$

From Equation~\eqref{eq:zeta 3} and under Assumption~\ref{ass:Restrict-Model}, we have $\zeta^{j,3}_{k\Delta}=2b(X^{j}_{k\Delta})\int_{k\Delta}^{(k+1)\Delta}{\sigma(X^{j}_{s})dW_s^j}$ and,
\begin{multline*}
    T_{1,3}^j(\ell, \ell^{\prime\prime}) = 3\mathbb{E}\left[\sum_{k = 0}^{n-1}{B_{\ell}(X^{j}_{k\Delta})B_{\ell^{\prime\prime}}(X^{j}_{k\Delta})\left(\zeta^{j,3}_{k\Delta}\right)^2}\right]\\
    = 12\mathbb{E}\left[\sum_{k = 0}^{n-1}{B_{\ell}(X^{j}_{k\Delta})B_{\ell^{\prime\prime}}(X^{j}_{k\Delta})b^{2}(X^{j}_{k\Delta})\mathbb{E}\left[\left(\int_{k\Delta}^{(k+1)\Delta}{\sigma(X^{j}_{s})dW_s^j}\right)^2| \mathcal{F}_{k\Delta}\right]}\right]\\
    \leq C\left[\mathbf{\Psi}_{K_N}\right]_{\ell, \ell^{\prime\prime}},
\end{multline*}
where the constant $C>0$ depends on $b$ and $\sigma$.
\end{enumerate}

We deduce from Equation~\eqref{eq:T1} and the respective upper bounds of $T_{1,1}^j(\ell, \ell^{\prime\prime}), ~ T_{1,2}^j(\ell, \ell^{\prime\prime})$ and $T_{1,3}^j(\ell, \ell^{\prime\prime})$ that there exists a constant $C>0$ such that
\begin{equation}\label{eq:mj-nu4}
    T_1^j(\ell, \ell^{\prime\prime}) \leq Cn\left[\mathbf{\Psi}_{K_N}\right]_{\ell,\ell^{\prime\prime}}, ~~ j = 1, \ldots, N, ~~ \ell, \ell^{\prime\prime} = -M, \ldots, K_N - 1.
\end{equation}

\subsubsection*{Upper-bound of $T_2^j(\ell, \ell^{\prime\prime})$}

For all $j \in \{1, \ldots, N\}$ and for all $\ell, \ell^{\prime\prime} \in \{-M, \ldots, K_N-1\}$, we have 
\begin{align*}
    T_2^j(\ell, \ell^{\prime\prime}) = \mathbb{E}\left[\sum_{k \neq k^{\prime}}{B_{\ell}(X^{j}_{k\Delta})B_{\ell^{\prime\prime}}(X^{j}_{k^{\prime}\Delta})\zeta^{j}_{k\Delta}\zeta^{j}_{k^{\prime}\Delta}}\right] = 2\mathbb{E}\left[\sum_{0 \leq k < k^{\prime} \leq n-1}{B_{\ell}(X^{j}_{k\Delta})B_{\ell^{\prime\prime}}(X^{j}_{k^{\prime}\Delta})\zeta^{j}_{k\Delta}\zeta^{j}_{k^{\prime}\Delta}}\right].
\end{align*}
For all $k,k^{\prime} \in \{0, \ldots, n-1\}$ such that $k < k^{\prime}$, the random variables $B_{\ell}(X^{j}_{k\Delta}), B_{\ell^{\prime\prime}}(X^{j}_{k^{\prime}\Delta})$ and $\zeta^{j}_{k\Delta}$ are measurable with respect to $\mathcal{F}_{k^{\prime}\Delta}$. Then we obtain the following:
\begin{equation}\label{eq:mj-nu5}
    T_2^j(\ell, \ell^{\prime\prime}) = 2\mathbb{E}\left[\sum_{0 \leq k < k^{\prime} \leq n-1}{B_{\ell}(X^{j}_{k\Delta})B_{\ell^{\prime\prime}}(X^{j}_{k^{\prime}\Delta})\zeta^{j}_{k\Delta}\mathbb{E}\left(\zeta^{j}_{k^{\prime}\Delta}|\mathcal{F}_{k^{\prime}\Delta}\right)}\right] = 0, ~~ j = 1, \ldots, N
\end{equation}
since for all $k^{\prime} \in \{1, \ldots, n-1\}, ~ j \in \{1, \ldots, N\}$ and from Equations~\eqref{eq:martingale-burkholder-davis-gundy} and \eqref{eq:Exp-zero},
$$\mathbb{E}\left(\zeta^{j}_{k^{\prime}\Delta}|\mathcal{F}_{k^{\prime}\Delta}\right) =  \mathbb{E}\left(\zeta^{j,1}_{k^{\prime}\Delta}|\mathcal{F}_{k^{\prime}\Delta}\right) + \mathbb{E}\left(\zeta^{j,2}_{k^{\prime}\Delta}|\mathcal{F}_{k^{\prime}\Delta}\right) + \mathbb{E}\left(\zeta^{j,3}_{k^{\prime}\Delta}|\mathcal{F}_{k^{\prime}\Delta}\right) = 0.$$
We finally deduce from Equations~\eqref{eq:mj-nu3}, \eqref{eq:mj-nu4} and \eqref{eq:mj-nu5} that 
\begin{multline}\label{eq:mj-nu6}
    \mathbb{E}\left(\underset{h\in\mathcal{S}_{K_N}, \ \|h\|_{n}=1}{\sup}{\nu^{2}(h)}\right) \leq \dfrac{C}{Nn}\sum_{\ell^{\prime} = -M}^{K_N-1}{\sum_{\ell = -M}^{K_N-1}{\sum_{\ell^{\prime\prime} = -M}^{K_N-1}{\left[\mathbf{\Psi}^{-1/2}_{K_N}\right]_{\ell,\ell^{\prime}}\left[\mathbf{\Psi}^{-1/2}_{K_N}\right]_{\ell^{\prime\prime},\ell^{\prime}}\left[\mathbf{\Psi}_{K_N}\right]_{\ell,\ell^{\prime\prime}}}}}\\
    = \dfrac{C}{Nn}\mathrm{Tr}\left(\mathbf{\Psi}^{-1}_{K_N}\mathbf{\Psi}_{K_N}\right)\\
    = C\dfrac{K_N}{Nn}.
\end{multline}
We then obtain from Equations~\eqref{eq:mj-nu6} and \eqref{eq:equation3-proof1},
\begin{equation*}
    \mathbb{E}\left[\left\|\widehat{\sigma}^{2}_{K_N}-\sigma^{2}_{A_N}\right\|^{2}_{n,N}\right] \leq C\left[\dfrac{A_N^{2\beta}}{K_N^{2\beta}} + \dfrac{K_N}{Nn} + \Delta^2 + \exp\left(-C_2\log^{3/2}(N)\right)\right],
\end{equation*}
where $N,n \rightarrow \infty$. Using the equivalence relation between the empirical norm $\|.\|_n$ and the pseudo-norm $\|.\|_{n,N}$, we obtain for $N$ large enough,
\begin{multline*}
    \mathbb{E}\left[\left\|\widehat{\sigma}^{2}_{K_N}-\sigma^{2}_{A_N}\right\|^{2}_{n}\right] = \mathbb{E}\left[\left\|\widehat{\sigma}^{2}_{K_N}-\sigma^{2}_{A_N}\right\|^{2}_{n}\mathds{1}_{\Omega_{n,N,K_N}}\right] + \mathbb{E}\left[\left\|\widehat{\sigma}^{2}_{K_N}-\sigma^{2}_{A_N}\right\|^{2}_{n}\mathds{1}_{\Omega^{c}_{n,N,K_N}}\right] \\
    \leq 2\mathbb{E}\left[\left\|\widehat{\sigma}^{2}_{K_N} - h\right\|^{2}_{n}\mathds{1}_{\Omega_{n,N,K_N}}\right] + 2\left\|h - \sigma^{2}_{A_N}\right\|_n^2 +  4mA_N^2\log(N)\mathbb{P}\left(\Omega^{c}_{n,N,K_N}\right)\\
    \leq 8\mathbb{E}\left[\left\|\widehat{\sigma}^{2}_{K_N} - \sigma_{A_N}^2\right\|^{2}_{n,N}\right] + 10\left\|h - \sigma^{2}_{A_N}\right\|_n^2 +  4K_NA_N^2\log(N)\mathbb{P}\left(\Omega^{c}_{n,N,K_N}\right)\\
    \leq C\left[\dfrac{A_N^{2\beta}}{K_N^{2\beta}} + \dfrac{K_N}{Nn} + \Delta^2 + \exp\left(C_2\log^{3/2}(N)\right)\right].
\end{multline*}
From Equation~\eqref{eq:Prob-OmegaComp}, and for $A_N \propto \sqrt{\log(N)}, ~ K_N = N^{2/(2\beta+1)\log^{-5/2}(N)}$ with $n \propto N$, there exists a constant $C > 0$ depending on $\beta \geq 1$ such that
\begin{equation*}
    \mathbb{E}\left[\left\|\widehat{\sigma}^{2}_{K_N}-\sigma^{2}_{A_N}\right\|_{n}^2\right] \leq C\log^{\beta}(N)(Nn)^{-2\beta/(2\beta+1)}.
\end{equation*}
\end{proof}

\subsection{Proof of Proposition~\ref{prop:Rate-realline}}

\begin{proof}
    From Equation~\eqref{eq:RelateTrunc}, we have the following:
    \begin{equation*}
       \mathbb{E}\left[\left\|\widetilde{\sigma}^{2}_{K_N}-\sigma^{2}\right\|_{n}\right] \leq \mathbb{E}\left[\left\|\widehat{\sigma}^{2}_{K_N}-\sigma^{2}\right\|_{n}\right] \leq \mathbb{E}\left[\left\|\widehat{\sigma}^{2}_{K_N}-\sigma^{2}_{A_N}\right\|_{n}\right] + 2\log(N)\underset{t \in [0,1]}{\mathrm{sup}}{\mathbb{P}(|X_t| > A_N)}.
    \end{equation*}
    From Theorem~\ref{thm:upper-bound} and the assumptions therein, there exists a constant $C > 0$ such that
    \begin{equation*}
        \mathbb{E}\left[\left\|\widehat{\sigma}^{2}_{K_N}-\sigma^{2}_{A_N}\right\|_{n}\right] \leq \left(\mathbb{E}\left[\left\|\widehat{\sigma}^{2}_{K_N}-\sigma^{2}_{A_N}\right\|_{n}^2\right]\right)^{1/2} \leq C\log^{\beta/2}(N)(Nn)^{-\beta/(2\beta+1)}.
    \end{equation*}
    We deduce that there exists a constant $C>0$ such that
    \begin{equation}\label{eq:Rate}
        \mathbb{E}\left[\left\|\widetilde{\sigma}^{2}_{K_N}-\sigma^{2}\right\|_{n}\right] \leq C\left[\log^{\beta/2}(N)(Nn)^{-\beta/(2\beta+1)} + \log(N)\underset{t \in [0,1]}{\mathrm{sup}}{\mathbb{P}(|X_t| > A_N)}\right]. 
    \end{equation}
    Under Assumptions~\ref{ass:LipFunctions}, \ref{ass:Ellipticity} and \ref{ass:Restrict-Model} and from Equations~\eqref{eq:TransDens} and \eqref{eq:TransitionParam}, with $S$ an increasing function that satisfies $S(0) = 0$, for all $t \in (0,1]$, we have
    \begin{multline}\label{eq:Exit1}
        \mathbb{P}(|X_t| > A_N) = \int_{|x| > A_N}{p_X(0, t, 0, x)dx}\\
        \leq \dfrac{C}{\sqrt{t}}\int_{A_N}^{+\infty}{\exp\left(-\dfrac{S^2(x)}{2t} + cS(x)\right)} + \dfrac{C}{\sqrt{t}}\int_{A_N}^{+\infty}{\exp\left(-\dfrac{S^2(-x)}{2t} + cS(-x)\right)}\\
        \leq \dfrac{C}{\sqrt{t}}\int_{A_N}^{+\infty}{\exp\left(-\dfrac{S^2(x)}{2t} + cS(x)\right)} + \dfrac{C}{\sqrt{t}}\int_{A_N}^{+\infty}{\exp\left(-\dfrac{S^2(-x)}{2t}\right)},
    \end{multline}
    where $C,c>0$ are constants. Focusing on the two last terms of equation~\eqref{eq:Exit1}, for $N$ large enough, we have, on one side,
    \begin{multline}\label{eq:Exit2}
        \dfrac{C}{\sqrt{t}}\int_{A_{N}}^{+\infty}{\exp\left[-\dfrac{S^2(x)}{2t} + cS(x)\right]dx}\\ 
        \leq C\kappa_1^{2}\int_{A_{N}}^{+\infty}{\left(\dfrac{1}{t\sigma(x)}\left(\int_{0}^{x}{\dfrac{1}{\sigma(u)}du}\right) - \dfrac{c}{\sigma(x)}\right)\exp\left[-\dfrac{S^2(x)}{2t} + cS(x)\right]dx}\\
        = C\kappa_1^{2}\exp\left[-\dfrac{1}{2}\left(\int_{0}^{A_{N}}{\dfrac{1}{\sigma(u)}du}\right)^2 + cS(A_N)\right]\\
        \leq C\exp\left(cA_N\right)\exp\left[-\dfrac{1}{2}\left(\int_{0}^{A_{N}}{\dfrac{1}{\sigma(u)}du}\right)^2\right],
    \end{multline}
    where $C>0$ and $c>0$ are new constants. On the other hand, for $N$ large enough, we have the following
    \begin{multline}\label{eq:Exit3}
        \dfrac{C}{\sqrt{t}}\int_{A_{N}}^{+\infty}{\exp\left(-\dfrac{S^2(-x)}{2t}\right)dx}\\ 
        \leq \kappa_1^{2}\int_{A_{N}}^{+\infty}{\dfrac{-1}{t\sigma(-x)}\left(\int_{0}^{-x}{\dfrac{1}{\sigma(u)}du}\right)\exp\left[-\dfrac{1}{2t}\left(\int_{0}^{-x}{\dfrac{1}{\sigma(u)}du}\right)^2\right]dx}\\
        = \kappa_1^{2}\exp\left[-\dfrac{1}{2}\left(\int_{0}^{-A_{N}}{\dfrac{1}{\sigma(u)}du}\right)^2\right].
    \end{multline}
    We finally obtain from Equations~\eqref{eq:Exit3}, \eqref{eq:Exit2} and \eqref{eq:Exit1} that there exists a constant $C>0$ such that for all $t \in [0,1]$,
    \begin{equation*}
        \mathbb{P}(|X_t| > A_N) \leq C\exp\left(cA_N\right)\left(\exp\left[-\dfrac{1}{2}\left(\int_{0}^{-A_{N}}{\dfrac{1}{\sigma(u)}du}\right)^2\right] + \exp\left[-\dfrac{1}{2}\left(\int_{0}^{-A_{N}}{\dfrac{1}{\sigma(u)}du}\right)^2\right]\right).
    \end{equation*}
    From Assumption~\ref{ass:Restrict-Model} and for $A_N = \sqrt{\frac{3\beta}{2\beta+1}\log(N)}$, we obtain the following:
    \begin{equation}\label{eq:Exit4}
        \underset{t \in [0,1]}{\sup}{~\mathbb{P}(|X_t| > A_N)} \leq 2C\exp\left(cA_N\right)\exp\left(-\dfrac{A_N^2}{2}\right) \leq C\exp\left(c\sqrt{\log(N)}\right)N^{-3\beta/2(2\beta+1)}.
    \end{equation}
    The final result is derived from Equations~\eqref{eq:Rate} and \eqref{eq:Exit4}.
\end{proof}

\subsection{Proof of Theorem~\ref{thm:lower-bound-compact}}

\begin{proof}
For the sake of simplicity and without loss of generality, the drift function is set to $b = 0$. We denote by $\widehat{\sigma}^2$ any estimator of $\sigma_{|I}^2$ built from the sample $\mathfrak{S}_{n,N}$ given by Equation~\eqref{eq:sample}, where $I = [A, B]$ is a compact interval, and we consider the estimation of $\sigma^2$ from a single diffusion path ($N=1$ and $n \rightarrow \infty$) and the estimation of $\sigma^2$ from repeated observations ($N,n \rightarrow \infty$). We want to prove that there exists a constant $c>0$ such that
    \begin{equation}\label{eq:L-bound1}
        \forall~ \w{\sigma}^{2}:~ \underset{\sigma^2\in\Sigma}{\sup}{~\E\left[\left\|\w{\sigma}^{2} - \sigma^{2}_{|I}\right\|^2\right]} \geq c\psi^{2}_{n,N},
    \end{equation}
    where $\Sigma = \Sigma_I(\beta, R)$ is the H\"older space and $\psi_{n,N} = (Nn)^{-\beta/(2\beta+1)}$, which is equivalent to
    \begin{equation*}
        \forall~ \w{\sigma}^{2},~~ \underset{\sigma^2\in\Sigma}{\sup}{~\E\left[\left\|\psi^{-1}_{n,N}(\w{\sigma}^{2} - \sigma^{2}_{|I})\right\|^2\right]} \geq c.
    \end{equation*}  
    Let $\Lambda_0 > 0$ be a constant to be determined later. We have 
    \begin{align*}
        \forall~ \w{\sigma}^{2},~~ \E\left[\left\|\psi^{-1}_{n,N}(\w{\sigma}^{2} - \sigma^{2}_{|I})\right\|^2\right] & \geq \E\left[\left\|\psi^{-1}_{n,N}(\w{\sigma}^{2} - \sigma^{2}_{|I})\right\|^2\one_{\left\|\psi^{-1}_{n,N}(\w{\sigma}^{2} - \sigma^{2}_{|I})\right\| \geq \Lambda_0}\right] \\
        & \geq \Lambda^{2}_0 \P\left(\left\|\w{\sigma}^{2} - \sigma^{2}_{|I}\right\| \geq s_{n,N}\right),
    \end{align*}
    where $s_{n,N} = \Lambda_0 \psi_{n,N}$. Then, to obtain Equation~\eqref{eq:L-bound1}, it suffices to prove the following result:
    \begin{equation*}
       \exists ~ p>0, ~~ \forall ~ \widehat{\sigma}^2, ~~ \underset{\sigma^2 \in \Sigma}{\sup}{\mathbb{P}\left(\left\|\widehat{\sigma}^2 - \sigma_{|I}^2\right\| \geq s_{n,N}\right)} > p.
    \end{equation*}
    Moreover, let $M \geq 2$ be an integer, and $\Sigma^{M} = \left\{\sigma_0^2,\cdots,\sigma_M^2\right\} \subset \Sigma$ a well-chosen finite set of $M+1$ hypotheses. Then, it suffices to prove that there exists a numerical constant $p_M>0$ depending on $M$ such that
    \begin{equation}
        \label{eq:minimax-probability}
        \forall~ \w{\sigma}^{2}:~~ \underset{\sigma^2 \in \Sigma^{M}}{\sup}{\P\left(\left\|\w{\sigma}^{2} - \sigma^{2}_{|I}\right\| \geq s_{n,N}\right)} \geq p_M.
    \end{equation}
    For each $j \in \{0,\cdots,M\}$, let $P_j = P_{\sigma_j^2} = \P_{\sigma_j^2 | \mathcal{F}_X}$ be the distribution of the diffusion process $X$ of diffusion coefficient $\sigma_j^2$ on the $\sigma-$algebra $\mathcal{F}_{X}$. The determination of the minimum probability $p_M$ follows the scheme presented in \cite{tsybakov2008introduction} {\it Theorem 2.5} . In fact, under the following three conditions:
    \begin{enumerate}
        \item $\Sigma^{M} \subset \Sigma_I(\beta,R), ~~ \beta \geq 1$, ~~ {\it i.e}, ~~ $\sigma_j^2 \in \Sigma, \ j=0,\cdots, M$,
        \item $\left\|\sigma^{2}_j - \sigma^{2}_k\right\| \geq 2s_{n,N} > 0, ~~ 0 \leq j < k \leq M$,
        \item $\frac{1}{M}\sum_{j=1}^{M}{K_{\mathrm{div}}(P^{\otimes N}_j,P^{\otimes N}_0)} \leq \alpha \log M$,
    \end{enumerate}
    with $P_j^{\otimes N} = \mathcal{L}\left(\bar{X}^1, \ldots, \bar{X}^N\right)$, the joint distribution of the sample $\mathfrak{S}_{n,N}$, Theorem 2.5 implies that
    \begin{equation}
    \label{eq:LowerBound-MinimaxProba}
        p_M \geq \frac{\sqrt{M}}{1+\sqrt{M}}\left(1 - 2\alpha - \sqrt{\frac{2\alpha}{\log M}}\right) > 0,
    \end{equation}
    where the real number $\alpha \in (0,1/8)$ and $K_{\mathrm{div}}$ is the Kullback divergence given for all probability measures $P$ and $Q$ by
    \begin{equation*} 
        K_{\mathrm{div}}(P,Q) = \begin{cases} \int{\log \frac{dP}{dQ} dP} = \mathrm{E}^{P}\left[\log \frac{dP}{dQ}\right], \ \ \mathrm{if} \ \ P \ll Q, \\\\ +\infty, \hspace*{1.6cm} \mathrm{otherwise},  \end{cases}
    \end{equation*}
    with $\mathrm{E}^{P}$ the expectation linked to $P$. Denote by $\nu$ the Lebesgue measure and suppose that $P \ll \nu$ and $Q \ll \nu$ and set 
    $$ p = \frac{dP}{d\nu}, ~~ q = \frac{dQ}{d\nu}. $$
    If $p > 0 ~ a.e.$ and $q > 0 ~ a.e.$, then the Kullback divergence is given by
    $$ K_{\mathrm{div}}(P,Q) = \int{\log\frac{p}{q} pd\nu} = \mathrm{E}^{P}\left[\log \frac{p}{q}\right].$$
    
    \subsection*{Construction of the set of hypotheses $\Sigma^{M}$}

    To satisfy the three conditions for the fulfillment of Theorem 2.5 in \cite{tsybakov2008introduction}, we propose a set of hypotheses $\Sigma^M = \left\{\sigma_j^2, ~~ j = 0, \ldots, M\right\}$ such that $\Sigma^M \subset \Sigma_I(\beta, R)$. For this purpose, let $c_0 > 0$ be a real number and define \begin{align*}
        m = &~ \left\lceil c_0 (Nn)^{\frac{1}{2\beta+1}} \right\rceil, ~~ h = \dfrac{1}{m}, ~~ x_k = \dfrac{k - 1/2}{m},
    \end{align*}
    \begin{equation}\label{eq:phi-k}
        \phi_k(x) = Rh^{\beta}K\left(\dfrac{x - x_k}{h}\right), ~~ k = 1, \ldots, m, ~~ x \in (0,1),
    \end{equation}
    where the function $K: \mathbb{R} \longrightarrow [0,\infty)$ is given by
    \begin{equation*}
        K(x) = aK_0(2x), ~ \mathrm{where} ~ K_0(x) = \exp\left(-\dfrac{1}{1-x^2}\right)\mathds{1}_{(-1,1)}(x),
    \end{equation*}
    with $a>0$ a sufficiently small real number. The function $K$ satisfies the following:
    \begin{equation}\label{eq:Fun-K}
        K \in \Sigma(\beta,1/2) \cap \mathcal{C}^{\infty}(\mathbb{R}), ~~ \mathrm{and} ~~ K(u) > 0 \iff u \in (-1/2, 1/2),
    \end{equation}
    and the functions $\phi_k$ belong to the H\"older class $\Sigma(\beta, R)$ on the interval $[0,1]$, and satisfy 
    \begin{equation}\label{eq:CompactSupp-Phi-k}
        \forall ~ k = 1, \ldots, m, ~~ \phi_k(x) > 0 \iff x \in \left(\dfrac{k-1}{m}, \dfrac{k}{m}\right),
    \end{equation}
    (see \cite{tsybakov2008introduction}, page 92). We extend the compact interval [0,1] considered in \cite{tsybakov2008introduction} to the compact interval $[A,B]$ by defining the following functions:
    \begin{equation}\label{eq:FunEta}
        \eta_k(x) = \phi_k\left(\dfrac{x-A}{B-A}\right), ~~ k = 1, \ldots, m, ~~ x \in [A,B].
    \end{equation}
    The functions $\eta_k$ are compactly supported with
    \begin{equation}\label{eq:support-eta}
       \mathrm{Supp}(\eta_k) = \mathcal{I}_{k,A,B} = \left(A + \dfrac{k-1}{m}(B-A), A + \dfrac{k}{m}(B-A)\right), ~~ k = 1, \ldots, m.    
    \end{equation}
   We then deduce that for all $k \in \left\{1, \ldots, m\right\}$, the function $\eta_k$ belongs to the H\"older class $\Sigma_I(\beta, R)$ on the interval $I = [A,B]$. We choose the finite set $\Sigma^{M}$ as follows:
   \begin{equation}\label{eq:set-hypothesis}
        \Sigma^{M} := \left\{\sigma_j^2 = 1 + \Gamma \sum_{k = 1}^{m}{w^{j}_k\eta_k}, ~ w^{j} = (w^{j}_{1},\cdots,w^{j}_{m-1})\in \{0,1\}^{m}, ~ \Gamma > 0, ~ j = 0, \cdots, M\right\},
    \end{equation}
    where $\sigma_0^2 \equiv 1$ is related to $w^{0} = (0,\ldots,0)$, and $\Gamma > 0$ is a numerical constant. Assumptions~\ref{ass:LipFunctions} and \ref{ass:Ellipticity} are satisfied with  
    $$\kappa_0^2 = 1 < \kappa_1^2 ~~ \mathrm{and} ~~ 0 < \Gamma \leq \dfrac{\kappa_1^2 - 1}{Rh^{\beta}\|K\|_{\infty}},$$
    and the hypotheses $\sigma_j^2$ belong to the H\"older class $\Sigma=\Sigma_I(\beta,R)$ by construction, which satisfies the first condition of {\it Theorem 2.5} in \cite{tsybakov2008introduction}. 

    \subsection*{Conclusion}

    Once the set of hypotheses $\Sigma^M$ is built such that $\Sigma^M \subset \Sigma$, the other two conditions required to fulfill Theorem 2.5 in \cite{tsybakov2008introduction} are satisfied through the following lemmas.

    \begin{lemma}\label{lm:condition2}
        Set
        \begin{equation*}
             \Lambda_0 = \frac{R\Gamma\|K\|}{2^{\beta+1}c_0^{\beta}}\sqrt{B-A} ~~ \mathrm{and} ~~ \Gamma = \dfrac{\kappa_1^2 - 1}{R\|K\|_{\infty}}.
       \end{equation*} 
       Then, for all $\sigma_j^2 , \sigma_{\ell}^2 \in \Sigma^M$, with $j, \ell \in \{0, \ldots, M\}$ such that $j \neq \ell$, we have
       \begin{equation*}
           \left\|\sigma_j^2 - \sigma_{\ell}^2\right\| \geq 2s_{n,N},
       \end{equation*}
       where $s_{n,N} = \Lambda_0\psi_{n,N} = \Lambda_0(Nn)^{-\beta/(2\beta+1)}$ with $\beta \geq 1$.
    \end{lemma}

    \begin{lemma}\label{lm:condition3}
        Suppose that $n \rightarrow \infty$ and $M \geq 2^{m/8}$. The constant $c_0>0$ is chosen so that 
        \begin{align*}
            &~ \dfrac{2C\Gamma^2R^2\|K\|_{\infty}^2}{c_0^{2\beta+1}} = \dfrac{1}{16} \times \dfrac{\log(2)}{8} ~~ \mathrm{for} ~~ N = 1 ~~ \mathrm{and} ~~ \beta > 2,\\\\
            &~ C\left(\dfrac{\Gamma^{2}R^{2}\|K\|^{2}_{\infty}}{c_0^{2\beta+1}} + \dfrac{2^{\beta+2}}{c_0^{\beta-1}}\right) = \dfrac{1}{16} \times \dfrac{\log(2)}{8} ~~ \mathrm{for} ~~ N \rightarrow \infty, ~~ n \propto N^2 ~~ \mathrm{and} ~~ \beta \geq 4.
        \end{align*}
        There following holds:
        \begin{equation*}
            \frac{1}{M}\sum_{j=1}^{M}{K_{\mathrm{div}}(P^{\otimes N}_{j},P^{\otimes N}_{0})} \leq
            \dfrac{1}{16}\log(M),
        \end{equation*}
        where, for all $j \in \{0, \ldots, M\}$, $P_j^{\otimes N} = \mathcal{L}\left(\bar{X}^1, \ldots, \bar{X}^N\right)$ is the joint distribution of the $N$ independent copies $\left\{\bar{X}^1, \ldots, \bar{X}^N\right\}$ of the discrete-time observation $\bar{X} = (X_{k\Delta})_{1\leq k\leq n}$ of the diffusion process $X$, solution of the stochastic differential equation $dX_t = \sigma_j(X_t)dW_t$.
    \end{lemma}

    The result of Lemma~\ref{lm:condition2} implies that the hypotheses $\sigma_j^2$ satisfy the second condition. Moreover, from the result of Lemma~\ref{lm:condition3},
    \begin{align*}
        \frac{1}{M}\sum_{j=1}^{M}{K_{\mathrm{div}}(P_{j}, P_{0})} \leq &~ \dfrac{1}{16}\log(M) ~~ \mathrm{for} ~~ N=1,\\\\
        \frac{1}{M}\sum_{j=1}^{M}{K_{\mathrm{div}}(P^{\otimes N}_{j},P^{\otimes N}_{0})} \leq &~ \dfrac{1}{16}\log(M) ~~ \mathrm{for} ~~ N \rightarrow \infty,
    \end{align*}
    which satisfies the third condition. Finally, applying Theorem 2.5 in \cite{tsybakov2008introduction}, we deduce from Equations~\eqref{eq:minimax-probability} and \eqref{eq:LowerBound-MinimaxProba} that
    \begin{align*}
       \underset{\w{\sigma}^{2}}{\inf}{\underset{\sigma^2 \in \Sigma^{M}}{\sup}{~\P\left(\|\w{\sigma}^{2}-\sigma^{2}_{|I}\| \geq s_{n,N}\right)}} \geq p_M \geq \frac{1}{2}\left(\dfrac{7}{8}-\dfrac{1}{2\sqrt{2\log(2)}}\right) > 0.
   \end{align*}
   Finally, there exists a constant $c>0$ depending on $\|K\|, \beta, c_0, \Gamma$ and $R$ such that
   \begin{align*}
       \underset{\w{\sigma}^{2}}{\inf}{\underset{\sigma^2 \in \Sigma^{M}}{\sup}{~\E\left[\|\w{\sigma}^{2}-\sigma^{2}_{|I}\|^{2}\right]}} \geq &~ cn^{-2\beta/(2\beta+1)}, ~~ \mathrm{for} ~~ N=1, ~~ n\rightarrow \infty ~~ \mathrm{and} ~~ \beta>2,\\
       \underset{\w{\sigma}^{2}}{\inf}{\underset{\sigma^2 \in \Sigma^{M}}{\sup}{~\E\left[\|\w{\sigma}^{2}-\sigma^{2}_{|I}\|^{2}\right]}} \geq &~ c(Nn)^{-2\beta/(2\beta+1)} ~~ \mathrm{for} ~~ N,n \rightarrow \infty ~~ \mathrm{and} ~~ \beta \geq 4.
   \end{align*}
\end{proof}

\subsubsection*{Proof of Lemma~\ref{lm:condition2}}

    \begin{proof}
    For all $\sigma_j^2, \sigma_{\ell}^2 \in \Sigma^M$, with $j, \ell \in \{0, \ldots, M\}$  such that $j \neq \ell$,
    \begin{align*}
        \left\|\sigma_j^2 - \sigma_{\ell}^2\right\|^2 = &~ \int_{A}^{B}{\left(\sigma_j^2(x) - \sigma_{\ell}^2(x)\right)^2dx}\\
        = &~ \Gamma^2\int_{A}^{B}{\left(\sum_{k=1}^{m}{(w_k^j - w_k^{\ell})\phi_k\left(\dfrac{x-A}{B-A}\right)}\right)^2dx}\\
        = &~ (B-A)\Gamma^2\int_{0}^{1}{\left(\sum_{k=1}^{m}{(w_k^j - w_k^{\ell})\phi_k(x)}\right)^2dx}.
    \end{align*}
    From Equation~\eqref{eq:CompactSupp-Phi-k}, the functions $\phi_k$ are compactly supported and these compact intervals are disjoint two-by-two. Then we obtain the following:
    \begin{align*}
        \left\|\sigma_{j}^2 - \sigma_{\ell}^2\right\|^2 = &~ (B-A)\Gamma^2\int_{0}^{1}{\left(\sum_{k=1}^{m}{(w_k^j - w_k^{\ell})\phi_k(x)}\right)^2dx}\\
        = &~ (B-A)\Gamma^2\sum_{k=1}^{m}{(w_k^j - w_k^{\ell})^2\int_{(k-1)/m}^{k/m}{\phi_k^2(x)dx}}\\
        = &~ (B-A)\Gamma^2\sum_{k=1}^{m}{\mathds{1}_{w_k^j \neq w_k^{\ell}}\int_{(k-1)/m}^{k/m}{\phi_k^2(x)dx}}.
    \end{align*}
   From the definition of the functions $\phi_k$ (Equation~\eqref{eq:phi-k}), we have for each $k = 1, \ldots, m$,
   \begin{align*}
       \int_{(k-1)/m}^{k/m}{\phi_k^2(x)dx} =&~ R^2h^{2\beta}\int_{(k-1)/m}^{k/m}{K^2\left(\dfrac{x-x_k}{h}\right)dx} = R^2h^{2\beta+1}\|K\|^2.
   \end{align*}
   We finally deduce from the precedent results that
    \begin{equation}
    \label{eq:DistanceHypotheses}
        \left\|\sigma_j^2 - \sigma_{\ell}^2\right\| \geq R\Gamma\|K\|\sqrt{B-A} h^{\beta + 1/2}\sqrt{\rho(w^j,w^{\ell})},
    \end{equation}
    where $\rho\left(w^{j},w^{\ell}\right) = \sum_{k=1}^{m}{\one_{w^{j}_{k} \neq w^{\ell}_{k}}}$ is the Hamming distance between the binary sequences $w^{j}$ and $w^{\ell}$. Moreover, for $m \geq 8$, we suppose that all binary sequences $w\in\Omega$ associated with the hypotheses of the finite set $\Sigma^{M}$ satisfy Lemma 2.9 in \cite{tsybakov2008introduction}, that is,
    \begin{equation}
    \label{eq:LowerBound-HammingDistance}
       M \geq 2^{m/8} ~~ \mathrm{and} ~~ \rho(w^{j}, w^{\ell}) \geq \frac{m}{8}, ~~ \forall~ 0<j<\ell<M.
    \end{equation}
    The condition $m \geq 8$ is satisfied for all $N \geq 1$ and $n \rightarrow \infty$ such that $Nn \geq (7/c_0)^{2\beta+1}$. 
    From Equations~\eqref{eq:DistanceHypotheses} and~\eqref{eq:LowerBound-HammingDistance}, we deduce that
    \begin{equation*}
        \left\|\sigma^{2}_j - \sigma^{2}_{\ell}\right\| \geq \frac{R\Gamma\|K\|}{2^{\beta}c_0^{\beta}}\sqrt{B-A}\psi_{n,N} = 2\Lambda_0\psi_{n,N} = 2s_{n,N},
    \end{equation*}    
    where
   \begin{equation*}
       \Lambda_0 = \frac{R\Gamma\|K\|}{2^{\beta+1}c_0^{\beta}}\sqrt{B-A}, ~~ \Gamma = \dfrac{\kappa_1^2 - 1}{R\|K\|_{\infty}}.
   \end{equation*} 
    \end{proof}

\subsection{Proof of Lemma~\ref{lm:condition3}}

\begin{proof}
    Under Assumptions~\ref{ass:LipFunctions}, \ref{ass:Ellipticity} and \ref{ass:RegularityBis}, for all $j\in \{0, \ldots, M\}$, we obtain from \cite{dacunha1986estimation}, \textit{Lemma 2} that the diffusion process $X$, solution of the stochastic differential equation $dX_t = \sigma_j(X_t)dW_t$, admits a transition density $p_j$ which is strictly positive on $\R$ and given, for all $(s, t,x,y) \in [0,1] \times [0,1] \times \mathbb{R}\times\mathbb{R}$ such that $t > s$, by
    \begin{multline}\label{eq:Transition-Density}
        p_j(s, t, x, y) = \dfrac{1}{\sqrt{2\pi(t-s)}}\dfrac{1}{\sigma_j(y)}\exp\left(H_j(y) - H_j(x) - \dfrac{(S_j(y) - S_j(x))^2}{2(t-s)}\right)\\
        \times \mathbf{E}\left[\exp\left((t-s)\int_{0}^{1}{G_{j,h}\left(z_u(S_j(x), S_j(y))+\sqrt{t-s}\widetilde{W}_u\right)du}\right)\right],
    \end{multline}    
    where
    \begin{equation}\label{eq:Functions-Sj-Hj-Zu}
        \forall x,y \in \mathbb{R}, ~~ S_j(x) =\int_{0}^{x}{\dfrac{du}{\sigma_j(u)}}, ~~ H_j(y) - H_j(x) = \log\left(\sqrt{\dfrac{\sigma_j(x)}{\sigma_j(y)}}\right), ~~ z_u(x,y) = (1-u)x + uy, 
    \end{equation}
    \begin{equation}\label{eq:Fonction-gjm}
        G_{j,h} = \dfrac{1}{4}\left[\sigma_j\sigma_j^{\prime\prime} - \dfrac{1}{2}(\sigma_j^{\prime})^2\right]\circ S_j^{-1},
    \end{equation}    
     and $\widetilde{W}$ is a Brownian bridge such that $\mathbf{E}(\widetilde{W}_t^2) = t(1-t), ~~ t \in [0,1]$. For the case $j=0$, we have $\sigma_0^2 = 1$ and the corresponding diffusion process coincides with the standard Brownian motion whose transition density $p_0$ is given, for all $(s,t,x,y) \in [0,1] \times [0,1]\times\mathbb{R}\times\mathbb{R}$ such that $t>s$, by
    $$p_0(s,t,x,y) = \dfrac{1}{\sqrt{2\pi(t-s)}}\exp\left(-\dfrac{(y-x)^2}{2(t-s)}\right).$$
    For all $j\in \{1, \ldots, M\}$ and $\ell \in \{1, \ldots, N\}$, we have
    \begin{equation}\label{eq:KullBack}
        \dfrac{1}{M}\sum_{j=1}^{M}{K_{\mathrm{div}}(P^{\otimes N}_j, P^{\otimes N}_0)} = \dfrac{1}{M}\sum_{j=1}^{M}{\sum_{\ell = 1}^{N}{\mathbb{E}^{P_j}\left[\log \frac{dP_j}{dP_0}(\bar{X}^{\ell})\right]}},
    \end{equation}   
    where, for each $j \in \{1, \ldots, M\}$, $\bar{X}^{\ell} = (X^{\ell}_{k\Delta})_{0 \leq k \leq n}$ is a discrete observation of the process $X$ solution of the stochastic differential equation $dX_t = \sigma_j(X_t)dW_t$. Then for all $j \in \{1, \ldots, M\}$, for all $\ell \in \{1, \ldots, N\}$, we obtain the following likelihood ratio:
    \begin{multline*}
        \dfrac{dP_j}{dP_0}(\bar{X}^{\ell}) = \prod_{k=1}^{n-1}{\dfrac{p_j(k\Delta, (k+1)\Delta, X^{\ell}_{k\Delta}, X^{\ell}_{(k+1)\Delta})}{p_0(k\Delta, (k+1)\Delta, X^{\ell}_{k\Delta}, X^{\ell}_{(k+1)\Delta})}}\\
        = \prod_{k=1}^{n-1}{\dfrac{1}{\sigma_j(X_{k\Delta}^{\ell})}\exp\left(\log\sqrt{\dfrac{\sigma_j(X^{\ell}_{k\Delta})}{\sigma_j(X^{\ell}_{(k+1)\Delta})}}-\frac{\left((S_j(X^{\ell}_{(k+1)\Delta}) - S_j(X^{\ell}_{k\Delta})\right)^{2}}{2 \Delta} + \dfrac{\left(X_{(k+1)\Delta}^{\ell} - X_{k\Delta}^{\ell}\right)^2}{2\Delta}\right)}\\
        \times \prod_{k=1}^{n-1}{\mathbf{E}\left[\exp\left(\Delta\int_{0}^{1}{G_{j,h}\left(z_u(S_j(X_{k\Delta}^{\ell}), S_j(X_{(k+1)\Delta}^{\ell})) + \Delta^{1/2}\widetilde{W}_u\right)du}\right)\right]}.
    \end{multline*}   
    We deduce that
    \begin{multline}\label{eq:Log-Kull1}
        \log\dfrac{dP_j}{dP_0}(\bar{X}^{\ell}) = \sum_{k=1}^{n-1}{\log\left(\mathbf{E}\left[\exp\left(\Delta\int_{0}^{1}{G_{j,h}\left(z_u(S_j(X_{k\Delta}^{\ell}), S_j(X_{(k+1)\Delta}^{\ell})) + \Delta^{1/2}\widetilde{W}_u\right)du}\right)\right]\right)}\\
        + \sum_{k=1}^{n-1}{\log\dfrac{1}{\sigma_j(X_{k\Delta}^{\ell})}} + \dfrac{1}{2}\sum_{k=1}^{n-1}{\log\dfrac{\sigma_j(X_{k\Delta}^{\ell})}{\sigma_j(X_{(k+1)\Delta}^{\ell})}} \\
        + \dfrac{1}{2}\sum_{k=1}^{n-1}{\left[\dfrac{\left(X_{(k+1)\Delta}^{\ell} - X_{k\Delta}^{\ell}\right)^2}{\Delta} - \dfrac{\left(S_j(X^{\ell}_{(k+1)\Delta}) - S_j(X^{\ell}_{k\Delta})\right)^{2}}{\Delta}\right]}.
    \end{multline}    
    We focus on the term $\left(S_j(X^{\ell}_{(k+1)\Delta}) - S_j(X^{\ell}_{k\Delta})\right)^{2}$ for $k \in [\![1,n-1]\!]$, with $S_j$ given by Equation~\eqref{eq:Functions-Sj-Hj-Zu}. For each $k \in [\![1,n-1]\!]$ and for all $\ell \in \{1, \ldots, N\}$, using the It\^o formula from the stochastic differential equation $dX_t = \sigma_j(X_t^{\ell})dW_t^{\ell}$, we have
    \begin{equation*}
        S_j(X_{k\Delta}^{\ell}) = W_{k\Delta}^{\ell} - \dfrac{1}{2}\int_{0}^{k\Delta}{\sigma_j^{\prime}(X_u^{\ell})du},
    \end{equation*}
    and
    \begin{multline*}
        \left(S(X_{(k+1)\Delta}^{\ell}) - S(X_{k\Delta}^{\ell})\right)^2 = \left(W_{(k+1)\Delta}^{\ell} - W_{k\Delta}^{\ell}\right)^2 - \dfrac{1}{2}\left(W_{(k+1)\Delta}^{\ell} - W_{k\Delta}^{\ell}\right)\int_{k\Delta}^{(k+1)\Delta}{\sigma^{\prime}(X_u^{\ell})du}\\
        + \dfrac{1}{4}\left(\int_{k\Delta}^{(k+1)\Delta}{\sigma^{\prime}(X_u^{\ell})du}\right)^2.
    \end{multline*}
    Then we obtain
    \begin{multline}\label{eq:kull1}
        \dfrac{\left(X_{(k+1)\Delta}^{\ell} - X_{k\Delta}^{\ell}\right)^2}{\Delta} - \dfrac{\left(S(X_{(k+1)\Delta}^{\ell}) - S(X_{k\Delta}^{\ell})\right)^2}{\Delta} = \dfrac{\left(X_{(k+1)\Delta}^{\ell} - X_{k\Delta}^{\ell}\right)^2}{\Delta} - \dfrac{\left(W_{(k+1)\Delta}^{\ell} - W_{k\Delta}^{\ell}\right)^2}{\Delta}\\
        + \dfrac{1}{2\Delta}\left(W_{(k+1)\Delta}^{\ell} - W_{k\Delta}^{\ell}\right)\int_{k\Delta}^{(k+1)\Delta}{\sigma^{\prime}(X_u^{\ell})du} - \dfrac{1}{4\Delta}\left(\int_{k\Delta}^{(k+1)\Delta}{\sigma^{\prime}(X_u^{\ell})du}\right)^2.
    \end{multline}
    From Equations~\ref{eq:kull1} and \eqref{eq:Log-Kull1}, we have for all $j \in \{1, \ldots, M\}$ and $\ell \in \{1, \ldots, N\}$,
    \begin{equation}\label{eq:Log-Kull3bis}
        \mathbb{E}^{P_j}\left[\log\dfrac{dP_j}{dP_0}(\bar{X}^{\ell})\right] = \mathbb{E}^{P_j}\left[T_n^{(1)}(\ell, j)\right] + \mathbb{E}^{P_j}\left[T_n^{(2)}(\ell, j)\right] + \dfrac{1}{2}\mathbb{E}^{P_j}\left[\sum_{k=1}^{n-1}{\log\dfrac{\sigma_j(X_{k\Delta}^{\ell})}{\sigma_j(X_{(k+1)\Delta}^{\ell})}}\right],
    \end{equation}
     where
     \begin{equation}\label{eq:bound-T1}
         T_n^{(1)}(\ell, j) = \sum_{k=1}^{n-1}{\log\left(\mathbf{E}\left[\exp\left(\Delta\int_{0}^{1}{G_{j,h}\left(z_u(S_j(X_{k\Delta}^{\ell}), S_j(X_{(k+1)\Delta}^{\ell})) + \Delta^{1/2}\widetilde{W}_u\right)du}\right)\right]\right)},
     \end{equation}
     and
    \begin{multline}\label{eq:bound-T2}
        T_n^{(2)}(\ell, j) = \sum_{k=1}^{n-1}{\log\dfrac{1}{\sigma_j(X_{k\Delta}^{\ell})}} + \dfrac{1}{2}\sum_{k=1}^{n-1}{\left[\dfrac{\left(X_{(k+1)\Delta}^{\ell} - X_{k\Delta}^{\ell}\right)^2}{\Delta} - \dfrac{\left(W_{(k+1)\Delta}^{\ell} - W_{k\Delta}^{\ell}\right)^2}{\Delta}\right]}\\
        + \dfrac{1}{2\Delta}\left(W_{(k+1)\Delta}^{\ell} - W_{k\Delta}^{\ell}\right)\int_{k\Delta}^{(k+1)\Delta}{\sigma_j^{\prime}(X_u^{\ell})du} - \dfrac{1}{4\Delta}\left(\int_{k\Delta}^{(k+1)\Delta}{\sigma_j^{\prime}(X_u^{\ell})du}\right)^2.
    \end{multline}
    The following lemmas give upper bounds for the first and second terms on the right-hand side of Equation~\eqref{eq:Log-Kull3bis}.
    \begin{lemma}\label{lm:bound-T1}
        For $\beta > 2, ~ \ell \in \{1, \ldots, N\}$ and $j \in \{1, \ldots, M\}$, there exists a constant $C>0$ such that
        \begin{equation*}
            \mathbb{E}^{P_j}\left[T_n^{(1)}(\ell, j)\right] \leq C h^{\beta-2}.
        \end{equation*}
    \end{lemma}

    \begin{lemma}\label{lm:bound-T2}
        For $\beta > 2, ~ \ell \in \{1, \ldots, N\}$ and $j \in \{1, \ldots, M\}$, the following holds:
        \begin{equation*}
           \mathbb{E}^{P_j}\left[T^{(2)}_{n}(\ell, j)\right] \leq C\left(\Gamma^{2}R^{2}\|K\|^{2}_{\infty}h^{2\beta}n + h^{\beta-2}\right),
       \end{equation*}
       where the constant $C>0$ depends on $R, ~ \Gamma, ~ \|K\|_{\infty}, ~ \|K^{\prime}\|_{\infty}$ and $\|K^{\prime\prime}\|_{\infty}$. 
    \end{lemma}

   We deduce from the respective results of Lemma~\ref{lm:bound-T1} and Lemma~\ref{lm:bound-T2}, and from Equation~\eqref{eq:Log-Kull3bis} that there exists a constant $C>0$ such that
   \begin{equation}\label{eq:Log-Kull5}
       \mathbb{E}^{P_j}\left[\log\dfrac{dP_j}{dP_0}\left(\bar{X}^{\ell}\right)\right] \leq C\left(\Gamma^{2}R^{2}\|K\|^{2}_{\infty}h^{2\beta}n + h^{\beta-2}\right) + \dfrac{1}{2}\mathbb{E}^{P_j}\left[\sum_{k=1}^{n-1}{\log\dfrac{\sigma_j(X_{k\Delta}^{\ell})}{\sigma_j(X_{(k+1)\Delta}^{\ell})}}\right].
   \end{equation}
   Focusing on the last term on the right-hand side of the above equation, we have the following:
   \begin{multline*}
       \dfrac{1}{2}\mathbb{E}^{P_j}\left[\sum_{k=1}^{n-1}{\log\dfrac{\sigma_j(X_{k\Delta}^{\ell})}{\sigma_j(X_{(k+1)\Delta}^{\ell})}}\right] = \dfrac{1}{2}\mathbb{E}^{P_j}\left[\sum_{k=1}^{n-1}{\log\left(\sigma_j(X_{k\Delta}^{\ell})\right) - \log\left(\sigma_j(X_{(k+1)\Delta}^{\ell})\right)}\right]\\
       = \dfrac{1}{2}\mathbb{E}^{P_j}\left[\log\left(\sigma_j(X_{\Delta}^{\ell})\right) - \log\left(\sigma_j(X_{n\Delta}^{\ell})\right)\right].
   \end{multline*}
   We apply the It\^o formula on $\log(\sigma_j^2(X_t^{\ell}))$ and obtain the following result.
   \begin{multline*}
       d\log(\sigma_j^2(X_t^{\ell})) = \dfrac{(\sigma_j^2)^{\prime}(X_t^{\ell})}{\sigma_j^2(X_t^{\ell})}dX_t^{\ell} + \dfrac{(\sigma_j^2)^{\prime\prime}(X_t^{\ell})\sigma_j^2(X_t^{\ell}) - (\sigma_j^2)^{\prime 2}(X_t^{\ell})}{2(\sigma_j^2(X_t^{\ell}))^2}\sigma_j^2(X_t^{\ell})dt\\
       = \dfrac{(\sigma_j^2)^{\prime}(X_t^{\ell})}{\sigma_j^2(X_t^{\ell})}\sigma_j(X_t^{\ell})dW_t^{\ell} + \dfrac{(\sigma_j^2)^{\prime\prime}(X_t^{\ell})\sigma_j^2(X_t^{\ell}) - (\sigma_j^2)^{\prime 2}(X_t^{\ell})}{2\sigma_j^2(X_t^{\ell})}dt. 
   \end{multline*}
   Since for all $j \in \{1, \ldots, M\}, ~ \sigma_j^2 \geq 1$, we deduce that
   \begin{multline}\label{eq:third-term-log}
       \dfrac{1}{2}\left|\mathbb{E}^{P_j}\left[\sum_{k=1}^{n-1}{\log\dfrac{\sigma_j(X_{k\Delta}^{\ell})}{\sigma_j(X_{(k+1)\Delta}^{\ell})}}\right]\right| = \dfrac{1}{2}\left|\mathbb{E}^{P_j}\left[\int_{\Delta}^{1}{\dfrac{(\sigma_j^2)^{\prime\prime}(X_t^{\ell})\sigma_j^2(X_t^{\ell}) - (\sigma_j^2)^{\prime 2}(X_t^{\ell})}{2\sigma_j^2(X_t^{\ell})}dt}\right]\right|\\
       \leq \dfrac{1}{4}\mathbb{E}^{P_j}\left[\int_{\Delta}^{1}{\left|(\sigma_j^2)^{\prime\prime}(X_t^{\ell})\sigma_j^2(X_t^{\ell}) - (\sigma_j^2)^{\prime 2}(X_t^{\ell})\right|dt}\right].
   \end{multline}
   For all $j \in \{1, \ldots, M\}$ and $\ell \in \{1, \ldots, N\}$, we have
   \begin{equation*}
       (\sigma_j^2)^{\prime\prime}(X_t^{\ell})\sigma_j^2(X_t^{\ell}) - (\sigma_j^2)^{\prime 2}(X_t^{\ell}) = 2\sigma_j^2(X_t^{\ell})\left[\sigma_j^{\prime\prime}(X_t^{\ell})\sigma_j(X_t^{\ell}) - 2\left(\sigma_j^{\prime}\right)^2(X_t^{\ell})\right],
   \end{equation*}
   and from Equations~\eqref{eq:sigma-prime}, \eqref{eq:Fun-fj}, \eqref{eq:sigma-primeprime} and \eqref{eq:Fun-gj} (proof of Lemma~\ref{lm:bound-T1}), we obtain the following:
   \begin{multline*}
       \left|(\sigma_j^2)^{\prime\prime}(X_t^{\ell})\sigma_j^2(X_t^{\ell}) - (\sigma_j^2)^{\prime 2}(X_t^{\ell})\right| \leq 2\left\|\sigma_j^2\right\|_{\infty}^2\left(\left\|\sigma_j^{\prime\prime}\right\|_{\infty}\left\|\sigma_j\right\|_{\infty} + 2\left\|\sigma_j^{\prime}\right\|_{\infty}^2\right)\\
       \leq 2\kappa_1^2h^{\beta-2}\left(\left\|\sigma_j\right\|_{\infty}\left\|g_j\right\|_{\infty} + 2\left\|f_j\right\|_{\infty}\right),
   \end{multline*}
   and from Equation~\eqref{eq:sup-Gjh} (proof of Lemma~\ref{lm:bound-T1}), we obtain for $h$ sufficiently small,
   \begin{equation*}
       \left|(\sigma_j^2)^{\prime\prime}(X_t^{\ell})\sigma_j^2(X_t^{\ell}) - (\sigma_j^2)^{\prime 2}(X_t^{\ell})\right| \leq \dfrac{6\kappa_1^3\Gamma R\left\|K^{\prime\prime}\right\|_{\infty}}{2(B-A)^2}h^{\beta-2}.
   \end{equation*}
   Then, back to Equation~\eqref{eq:third-term-log}, we obtain:
   \begin{equation}\label{eq:Kullback-term1}
       \dfrac{1}{2}\left|\mathbb{E}^{P_j}\left[\sum_{k=1}^{n-1}{\log\dfrac{\sigma_j(X_{k\Delta}^{\ell})}{\sigma_j(X_{(k+1)\Delta}^{\ell})}}\right]\right| \leq Ch^{\beta-2},
   \end{equation}
   where $C = \dfrac{3\kappa_1^3\Gamma R\left\|K^{\prime\prime}\right\|_{\infty}}{8(B-A)^2}$. We deduce from Equations~\eqref{eq:Log-Kull5} and \eqref{eq:Kullback-term1} that for $\beta > 2$ and $n \rightarrow \infty$, there exists a constant $C>0$ such that
   \begin{equation}\label{eq:Kullback-Inequality1}
       \mathbb{E}^{P_j}\left[\log\dfrac{dP_j}{dP_0}\left(\bar{X}^{\ell}\right)\right] \leq C\left(\Gamma^{2}R^{2}\|K\|^{2}_{\infty}h^{2\beta}n + h^{\beta-2}\right).
   \end{equation}
   From Equations~\eqref{eq:Kullback-Inequality1} and \eqref{eq:KullBack}, we have
   \begin{multline}\label{eq:Kullback-Inequality2}
       \frac{1}{M}\sum_{j=1}^{M}{K_{\mathrm{div}}(P^{\otimes N}_{j},P^{\otimes N}_{0})} = \frac{1}{M}\sum_{j=1}^{M}{\sum_{\ell=1}^{N}{\mathrm{E}^{P_j}\left[\log \frac{dP_j}{dP_0}(\bar{X}^{\ell})\right]}}\\
       \leq C\left(\Gamma^{2}R^{2}\|K\|^{2}_{\infty}h^{2\beta}Nn + Nh^{\beta-2}\right).
   \end{multline}
   To conclude the proof, we distinguish the two cases $N = 1$ and $N \rightarrow \infty$. Then, for the case $N = 1$ and $n \rightarrow \infty$, since $m = \lceil c_0n^{1/(2\beta+1)} \rceil \geq c_0n^{1/(2\beta+1)}, ~ h = m^{-1}$, $\beta > 2$ and $n \rightarrow \infty$, we have the following:
   \begin{multline*}
       \dfrac{1}{M}\sum_{j=1}^{M}{K_{\mathrm{div}}(P_j, P_0)} \leq C\left(\Gamma^{2}R^{2}\|K\|^{2}_{\infty}h^{2\beta}n + h^{\beta-2}\right)\\
       \leq C\Gamma^{2}R^{2}\|K\|^{2}_{\infty}\dfrac{m^{2\beta+1}}{c_0^{2\beta+1}} \times m^{-2\beta} + Cm^{2-\beta} \leq \dfrac{2C\Gamma^2R^2\|K\|_{\infty}^2}{c_0^{2\beta+1}}m.
   \end{multline*}
   Since $c_0>0$ is chosen so that
   \begin{equation*}
       \dfrac{2C\Gamma^2R^2\|K\|_{\infty}^2}{c_0^{2\beta+1}} = \dfrac{1}{16} \times \dfrac{\log(2)}{8},
   \end{equation*}
   and for $M \geq 2^{m/8}$, we obtain that:
   \begin{equation*}
       \dfrac{1}{M}\sum_{j=1}^{M}{K_{\mathrm{div}}(P_j, P_0)} \leq \dfrac{1}{16}\times \dfrac{\log(2)}{8}m \leq \dfrac{1}{16}\log(M). 
   \end{equation*}
   For the case $N,n \rightarrow \infty$, we have $m = \lceil c_0(Nn)^{1/(2\beta+1)} \rceil \geq c_0(Nn)^{1/(2\beta+1)}$, and
   \begin{multline*}
       \dfrac{1}{M}\sum_{j=1}^{M}{K_{\mathrm{div}}(P^{\otimes N}_{j},P^{\otimes N}_{0})} \leq C\left(\Gamma^{2}R^{2}\|K\|^{2}_{\infty}h^{2\beta}Nn + Nh^{\beta-2}\right)\\
       \leq C\Gamma^{2}R^{2}\|K\|^{2}_{\infty}\dfrac{m^{2\beta+1}}{c_0^{2\beta+1}}m^{-2\beta} + C\dfrac{m^{2\beta+1}}{nc_0^{2\beta+1}}m^{2-\beta}\\
       = C\Gamma^{2}R^{2}\|K\|^{2}_{\infty}\dfrac{m}{c_0^{2\beta+1}} + C\dfrac{m^{\beta+2}}{c_0^{2\beta+1}n}m.
   \end{multline*}
   Since $n \propto N^2$ and $\beta \geq 4$, there exists a constant $C>0$ such that
   \begin{multline*}
       \dfrac{1}{M}\sum_{j=1}^{M}{K_{\mathrm{div}}(P^{\otimes N}_{j},P^{\otimes N}_{0})} \leq C\left(\dfrac{\Gamma^{2}R^{2}\|K\|^{2}_{\infty}}{c_0^{2\beta+1}} + \dfrac{2^{\beta+2}N^{(3\beta+6)/(2\beta+1)}}{c_0^{\beta-1}N^2}\right)m\\
       = C\left(\dfrac{\Gamma^{2}R^{2}\|K\|^{2}_{\infty}}{c_0^{2\beta+1}} + \dfrac{2^{\beta+2}N^{(4-\beta)/(2\beta+1)}}{c_0^{\beta-1}}\right)m\\
       \leq C\left(\dfrac{\Gamma^{2}R^{2}\|K\|^{2}_{\infty}}{c_0^{2\beta+1}} + \dfrac{2^{\beta+2}}{c_0^{\beta-1}}\right)m.
   \end{multline*}
   Since the constant $c_0>0$ is chosen so that
   \begin{equation*}
       C\left(\dfrac{\Gamma^{2}R^{2}\|K\|^{2}_{\infty}}{c_0^{2\beta+1}} + \dfrac{2^{\beta+2}}{c_0^{\beta-1}}\right) = \dfrac{1}{16} \times \dfrac{\log(2)}{8},
   \end{equation*}
   and for $M \geq 2^{m/8}$, we finally obtain the following: 
   \begin{align*}
       \dfrac{1}{M}\sum_{j=1}^{M}{K_{\mathrm{div}}(P^{\otimes N}_{j},P^{\otimes N}_{0})} \leq \dfrac{1}{16}\times \dfrac{\log(2)}{8}m \leq \dfrac{1}{16}\log(M).
   \end{align*}
\end{proof}

\subsection{Proof of Lemma~\ref{lm:bound-T1}}

\begin{proof}
     From Equation~\eqref{eq:bound-T1}, we have for all $j \in \{1, \ldots, M\}$ and $\ell \in \{1, \ldots, N\}$,
     \begin{equation*}
         T_n^{(1)}(\ell, j) = \sum_{k=1}^{n-1}{\log\left(\mathbf{E}\left[\exp\left(\Delta\int_{0}^{1}{G_{j,h}\left(z_u(S_j(X_{k\Delta}^{\ell}), S_j(X_{(k+1)\Delta}^{\ell})) + \Delta^{1/2}\widetilde{W}_u\right)du}\right)\right]\right)},
     \end{equation*}
     where
     \begin{equation*}
         G_{j,h} = \dfrac{1}{4}\left[\sigma_j\sigma^{\prime\prime} - \dfrac{1}{2}\left(\sigma_j^{\prime}\right)^{2}\right] \circ S_j^{-1}, ~~ S_j: x \in \mathbb{R} \mapsto \int_{0}^{x}{\dfrac{1}{\sigma_j(u)}du} ~~ \mathrm{and} ~~ z_u: (x,y) \mapsto ux+(1-u)y.
     \end{equation*}
      Since for any $j \in \{1, \ldots, M\}$, the function $\sigma_j$ is given by $\sigma_j = (\sigma_j^2)^{1/2}$, we obtain from Equations~\eqref{eq:set-hypothesis}, \eqref{eq:FunEta}, \eqref{eq:CompactSupp-Phi-k} and \eqref{eq:phi-k} that for all $x \in [A,B]$,
    \begin{multline*}
        \sigma_j^{\prime}(x) = \dfrac{1}{2}(\sigma_j^2)^{\prime}(x)(\sigma_j^2(x))^{-1/2} = \dfrac{\Gamma}{2}(\sigma_j^2(x))^{-1/2} \sum_{k=1}^{m}{w_k^j \eta_k^{\prime}(x)}\\
        = \dfrac{\Gamma R h^{\beta-1}}{2(B-A)}(\sigma_j^2(x))^{-1/2}\sum_{k=1}^{m}{w_k^jK^{\prime}\left(\dfrac{x-A-(B-A)x_k}{(B-A)h}\right)\mathds{1}_{\mathcal{I}_{k,A,B}}(x)},
    \end{multline*}
    where for $k \in [\![1, m]\!]$, $\mathcal{I}_{k,A,B}$ is given by Equation~\eqref{eq:support-eta}. Since the functions $\eta_k$ are compactly supported, and their respective supports $\mathcal{I}_{k,A,B}$ are two by two joint intervals (see Equations~\eqref{eq:CompactSupp-Phi-k} and \eqref{eq:FunEta}), we obtain the following:
    \begin{equation}\label{eq:sigma-prime}
        \left[\sigma_j^{\prime}(x)\right]^2 = \dfrac{1}{4}\left[(\sigma_j^2)^{\prime}(x)\right]^2(\sigma_j^2(x))^{-1} = \dfrac{\Gamma^2}{4}(\sigma_j^2(x))^{-1} \sum_{k=1}^{m}{w_k^j \left[\eta_k^{\prime}(x)\right]^2} = h^{\beta-2}f_j(x),
    \end{equation}
    where
    \begin{equation}\label{eq:Fun-fj}
        f_j(x) = \dfrac{\Gamma^2 R^2 h^{\beta}}{4(B-A)^2}(\sigma_j^2(x))^{-1}\sum_{k=1}^{m}{w_k^jK^{\prime 2}\left(\dfrac{x-A-(B-A)x_k}{(B-A)h}\right)\mathds{1}_{\mathcal{I}_{k,A,B}}(x)}.
    \end{equation}
    We also have
    \begin{equation}\label{eq:sigma-primeprime}
        \sigma_j^{\prime\prime}(x) = \dfrac{1}{2}(\sigma_j^2)^{\prime\prime}(x)(\sigma_j^2(x))^{-1/2} - \dfrac{1}{4}\left[(\sigma_j^2)^{\prime}(x)\right]^2(\sigma_j^2(x))^{-3/2} = h^{\beta-2}g_j(x),
    \end{equation}
    where
    \begin{multline}\label{eq:Fun-gj}
        g_j(x) = \dfrac{\Gamma R}{2(B-A)^2}(\sigma_j^2(x))^{-1/2}\sum_{k=1}^{m}{w_k^jK^{\prime\prime}\left(\dfrac{x-A-(B-A)x_k}{(B-A)h}\right)\mathds{1}_{\mathcal{I}_{k,A,B}}(x)}\\
        - \dfrac{\Gamma^2R^2h^{\beta}}{16(B-A)^2}(\sigma_j^2(x))^{-5/2}\sum_{k=1}^{m}{w_k^j K^{\prime 2}\left(\dfrac{x-A-(B-A)x_k}{(B-A)h}\right)\mathds{1}_{\mathcal{I}_{k,A,B}}(x)}.
    \end{multline}
    We deduce from Equations~\eqref{eq:Fonction-gjm}, \eqref{eq:sigma-prime} and \eqref{eq:sigma-primeprime} that for all $x \in \mathbb{R}$,
    \begin{equation*}
        G_{j,h}(x) = \dfrac{1}{4}\left[\sigma_j\sigma_j^{\prime\prime} - \dfrac{1}{2}(\sigma_j^{\prime})^2\right] \circ S_j^{-1}(x) = h^{\beta-2}T_{j,h} \circ S_j^{-1}(x),
    \end{equation*}
    where
    \begin{equation}\label{eq:Tjh}
        T_{j,h} = \dfrac{1}{4}\left[\sigma_jg_j - \dfrac{1}{2}f_j\right]. 
    \end{equation}
    Since $\beta > 2$, from the definition of the function $K$ and from Equations~\eqref{eq:Fun-gj} and \eqref{eq:Fun-fj}, we have for $h = 1/m$ small enough,
    \begin{multline}\label{eq:sup-Gjh}
        \left\|G_{j,h}\right\|_{\infty} \leq h^{\beta-2}\left\|T_{j,h}\right\|_{\infty} \leq h^{\beta-2}\left(\left\|\sigma_j\right\|_{\infty}\left\|g_j\right\|_{\infty} + \left\|f_j\right\|_{\infty}\right)\\
        \leq h^{\beta-2}\left[\kappa_1\left(\dfrac{\Gamma R\left\|K^{\prime\prime}\right\|_{\infty}}{2(B-A)^2} + \dfrac{\Gamma^2R^2h^{\beta}\left\|K^{\prime}\right\|_{\infty}^2}{16(B-A)^2}\right) + \dfrac{\Gamma^2R^2h^{\beta}\left\|K^{\prime}\right\|_{\infty}^2}{4(B-A)^2}\right]\\
        \leq \dfrac{3\kappa_1\Gamma R\left\|K^{\prime\prime}\right\|_{\infty}}{2(B-A)^2}h^{\beta-2},
    \end{multline}
    which implies that $\|G_{j,h}\|_{\infty} \rightarrow 0$ as $h \rightarrow 0$. Moreover, $\Delta \rightarrow 0$ as $n \rightarrow \infty$. We then deduce that for all $x,y \in \mathbb{R}$,
    \begin{multline}\label{eq:Factor-Density}
        \mathbf{E}\left[\exp\left(\Delta\int_{0}^{1}{G_{j,h}\left(z_u(S_j(x), S_j(y)) + \Delta^{1/2}\widetilde{W}_u\right)}\right)\right] = 1 \\
        + \mathbf{E}\left[\Delta\int_{0}^{1}{G_{j,h}\left(z_u(S_j(x), S_j(y)) + \Delta^{1/2}\widetilde{W}_u\right)}\right]\\
        + o\left(\mathbf{E}\left[\Delta\int_{0}^{1}{G_{j,h}\left(z_u(S_j(x), S_j(y)) + \Delta^{1/2}\widetilde{W}_u\right)}\right]\right).
    \end{multline}
   Moreover, from Equation~\eqref{eq:sup-Gjh}, we obtain
   \begin{equation*}
       \forall ~ x,y \in \mathbb{R}, ~~ \mathbf{E}\left[\Delta\int_{0}^{1}{G_{j,h}\left(z_u(S_j(x), S_j(y)) + \Delta^{1/2}\widetilde{W}_u\right)}\right] \leq \dfrac{3\kappa_1\Gamma R\left\|K^{\prime\prime}\right\|_{\infty}}{2(B-A)^2}\Delta h^{\beta-2}.
   \end{equation*}
   Then, from Equation~\eqref{eq:Factor-Density}, there exists a constant $C>0$ depending on $\kappa_1, ~ \Gamma, ~ R$ and $\|K^{\prime\prime}\|_{\infty}$ such that
   \begin{equation*}
       \forall ~ x,y \in \mathbb{R}, ~~ \mathbf{E}\left[\exp\left(\Delta\int_{0}^{1}{G_{j,h}\left(z_u(S_j(x), S_j(y)) + \Delta^{1/2}\widetilde{W}_u\right)}\right)\right] \leq 1 + C\Delta h^{\beta-2}.
   \end{equation*}
   We finally deduce that for all $j \in \{1, \ldots, M\}$ and $\ell \in \{1, \ldots, N\}$,
   \begin{multline*}
       \sum_{k=1}^{n-1}{\log\left(\mathbf{E}\left[\exp\left(\Delta\int_{0}^{1}{G_{j,h}\left(z_u(S_j(X_{k\Delta}^{\ell}), S_j(X_{(k+1)\Delta}^{\ell})) + \Delta^{1/2}\widetilde{W}_u\right)du}\right)\right]\right)}\\
       \leq \sum_{k=1}^{n-1}{\log(1+C\Delta h^{\beta-2})} \leq Ch^{\beta-2},
   \end{multline*}
   which concludes the proof.
\end{proof}

\subsection{Proof of Lemma~\ref{lm:bound-T2}}

\begin{proof}
    Recall that for all $j \in \{1, \ldots, M\}$ and $\ell \in \{1, \ldots, N\}$,
    \begin{multline*}
        T_n^{(2)}(\ell, j) = \sum_{k=1}^{n-1}{\log\dfrac{1}{\sigma_j(X_{k\Delta}^{\ell})}} + \dfrac{1}{2}\sum_{k=1}^{n-1}{\left[\dfrac{\left(X_{(k+1)\Delta}^{\ell} - X_{k\Delta}^{\ell}\right)^2}{\Delta} - \dfrac{\left(W_{(k+1)\Delta}^{\ell} - W_{k\Delta}^{\ell}\right)^2}{\Delta}\right]}\\
        + \dfrac{1}{2\Delta}\left(W_{(k+1)\Delta}^{\ell} - W_{k\Delta}^{\ell}\right)\int_{k\Delta}^{(k+1)\Delta}{\sigma_j^{\prime}(X_u^{\ell})du} - \dfrac{1}{4\Delta}\left(\int_{k\Delta}^{(k+1)\Delta}{\sigma_j^{\prime}(X_u^{\ell})du}\right)^2.
    \end{multline*}
    Since for $k \in \{1, \ldots n-1\}$, $\left(X_{(k+1)\Delta}^{\ell} - X_{k\Delta}^{\ell}\right)^2 = \int_{k\Delta}^{(k+1)\Delta}{\sigma_j(X_s^{\ell})dW_s}$, we have
    \begin{multline}\label{eq:bound-T20}
        \mathbb{E}^{P_j}\left[T_n^{(2)}(\ell, j)\right] \leq \sum_{k=1}^{n-1}{\mathbb{E}^{P_j}\left[\log\dfrac{1}{\sigma_j(X_{k\Delta}^{\ell})}\right]} + \dfrac{1}{2}\sum_{k=1}^{n-1}{\left[\dfrac{1}{\Delta}\mathbb{E}^{P_j}\left(\int_{k\Delta}^{(k+1)\Delta}{\sigma_j^2(X_s^{\ell})ds}\right) - 1\right]}\\
        + \dfrac{1}{2\Delta}\sum_{k=1}^{n-1}{\mathbb{E}^{P_j}\left[\sigma_j^{\prime}(X_{k\Delta}^{\ell})\left(W_{(k+1)\Delta}^{\ell} - W_{k\Delta}^{\ell}\right)\right]}\\
        + \dfrac{1}{2\Delta}\sum_{k=1}^{n-1}{\mathbb{E}^{P_j}\left[\left(W_{(k+1)\Delta}^{\ell} - W_{k\Delta}^{\ell}\right)\int_{k\Delta}^{(k+1)\Delta}{\left(\sigma_j^{\prime}(X_u^{\ell}) - \sigma_j^{\prime}(X_{k\Delta}^{\ell})\right)du}\right]}. 
    \end{multline}
    On the one hand, considering the natural filtration $\left(\mathcal{F}_t\right)_{t \in [0,1]}$ of the diffusion process $X$, and since for all $k \in \{1, \ldots, n-1\}$, $W_{(k+1)\Delta}^{\ell} - W_{k\Delta}^{\ell}$ is independent of $\mathcal{F}_{k\Delta}$, we have the following: 
    \begin{multline}\label{eq:bound-T21}
        \sum_{k=1}^{n-1}{\mathbb{E}^{P_j}\left[\sigma_j^{\prime}(X_{k\Delta}^{\ell})\left(W_{(k+1)\Delta}^{\ell} - W_{k\Delta}^{\ell}\right)\right]} = \sum_{k=1}^{n-1}{\mathbb{E}^{P_j}\left[\sigma_j^{\prime}(X_{k\Delta}^{\ell})\mathbb{E}^{P_j}\left(W_{(k+1)\Delta}^{\ell} - W_{k\Delta}^{\ell} \bigm\vert \mathcal{F}_{k\Delta}\right)\right]}\\
        = \sum_{k=1}^{n-1}{\mathbb{E}^{P_j}\left[\sigma_j^{\prime}(X_{k\Delta}^{\ell})\mathbb{E}^{P_j}\left(W_{(k+1)\Delta}^{\ell} - W_{k\Delta}^{\ell}\right)\right]} = 0.
    \end{multline}
    On the other hand, using the Cauchy-Schwarz inequality and Equations~\eqref{eq:sigma-primeprime} and \eqref{eq:Cons-Ass2.1}, we have the following.
    \begin{multline}\label{eq:bound-T22}
        \dfrac{1}{2\Delta}\sum_{k=1}^{n-1}{\mathbb{E}^{P_j}\left[\left(W_{(k+1)\Delta}^{\ell} - W_{k\Delta}^{\ell}\right)\int_{k\Delta}^{(k+1)\Delta}{\left(\sigma_j^{\prime}(X_u^{\ell}) - \sigma_j^{\prime}(X_{k\Delta}^{\ell})\right)du}\right]}\\
        \leq \dfrac{1}{2\Delta}\sum_{k=1}^{n-1}{\left(\mathbb{E}^{P_j}\left[\left(W_{(k+1)\Delta}^{\ell} - W_{k\Delta}^{\ell}\right)^2\right]\right)^{1/2}\left(\mathbb{E}^{P_j}\left[\Delta\int_{k\Delta}^{(k+1)\Delta}{\left(\sigma_j^{\prime}(X_u^{\ell}) - \sigma_j^{\prime}(X_{k\Delta}^{\ell})\right)^2du}\right]\right)^{1/2}}\\
        \leq \dfrac{1}{2}\sum_{k=1}^{n-1}{\left\|\sigma_j^{\prime\prime}\right\|_{\infty}\left(\int_{k\Delta}^{(k+1)\Delta}{\mathbb{E}^{P_j}\left[\left|X_u^{\ell} - X_{k\Delta}^{\ell}\right|^2\right]}\right)^{1/2}} \leq Ch^{\beta-2},
    \end{multline}
    where the constant $C>0$ depends on $R, ~\Gamma, ~\|K^{\prime}\|_{\infty}$ and $\|K^{\prime\prime}\|_{\infty}$. We deduce from Equations~\eqref{eq:bound-T20}, \eqref{eq:bound-T21} and \eqref{eq:bound-T22} that
    \begin{multline}\label{eq:bound-T23}
        \mathbb{E}^{P_j}\left[T_n^{(2)}(\ell, j)\right] \leq \sum_{k=1}^{n-1}{\mathbb{E}^{P_j}\left[\log\dfrac{1}{\sigma_j(X_{k\Delta}^{\ell})}\right]} + \dfrac{1}{2}\sum_{k=1}^{n-1}{\left[\dfrac{1}{\Delta}\mathbb{E}^{P_j}\left(\int_{k\Delta}^{(k+1)\Delta}{\sigma_j^2(X_s^{\ell})ds}\right) - 1\right]} + Ch^{\beta-2}\\
        \leq \sum_{k=1}^{n-1}{\mathbb{E}^{P_j}\left[\log\dfrac{1}{\sigma_j(X_{k\Delta}^{\ell})}\right]} + \dfrac{1}{2}\sum_{k=1}^{n-1}{\left(\mathbb{E}^{P_j}\left[\sigma_j^2(X_{k\Delta}^{\ell})\right] - 1\right)}\\
        + \dfrac{1}{2\Delta}\sum_{k=1}^{n-1}{\int_{k\Delta}^{(k+1)\Delta}{\mathbb{E}^{P_j}\left[\sigma_j^2(X_s^{\ell}) - \sigma_j^2(X_{k\Delta}^{\ell})\right]}} + Ch^{\beta-2}.
    \end{multline}
    For all $\ell \in [\![1,N]\!]$ and for all $k \in [\![1,n-1]\!]$, since $\beta > 2$ and $h = 1/m$, we set   $$\mathcal{G}_m(X^{\ell}_{k\Delta}) = \Gamma\sum_{r=1}^{m}{w^{j}_{r}\eta_{r}(X^{\ell}_{k\Delta})} \longrightarrow 0 ~~ \mathrm{as} ~~ m \rightarrow \infty, ~~ \mathrm{and} ~~ 0 \leq \mathcal{G}_m(X^{\ell}_{k\Delta}) \leq 2 R\Gamma\|K\|_{\infty}h^{\beta}.$$
   So, for all $j = 1, \ldots, M, ~~ \sigma_j^2(X^{\ell}_{k\Delta}) = 1 + \mathcal{G}_m(X^{\ell}_{k\Delta})$. Then, using the second-order Taylor-Young expansion, 
   \begin{align*}
       \sigma_j(X_{k\Delta}^{\ell}) = \sqrt{1+\mathcal{G}_m(X^{\ell}_{k\Delta})} = &~ 1 + \dfrac{1}{2}\mathcal{G}_m(X^{\ell}_{k\Delta}) - \dfrac{1}{8}\mathcal{G}^{2}_{m}(X^{\ell}_{k\Delta}) + \mathrm{O}\left(\mathcal{G}^{2}_{m}(X^{\ell}_{k\Delta})\right)\\
       \dfrac{1}{\sigma_j(X^{\ell}_{k\Delta})} = &~ 1 - \dfrac{1}{2}\mathcal{G}_m(X^{\ell}_{k\Delta}) + \dfrac{3}{8}\mathcal{G}_m^2(X^{\ell}_{k\Delta}) + \mathrm{O}\left(\mathcal{G}_m^2(X_{k\Delta}^{\ell})\right)\\
       \log\dfrac{1}{\sigma_j(X^{\ell}_{k\Delta})} = &~ - \dfrac{1}{2}\mathcal{G}_m(X^{\ell}_{k\Delta}) + \dfrac{1}{4}\mathcal{G}_m^2(X^{\ell}_{k\Delta}) + \mathrm{O}\left(\mathcal{G}_m^2(X_{k\Delta}^{\ell})\right).
   \end{align*} 
   We deduce from Equation~\eqref{eq:bound-T23}, for each $\ell \in [\![1,N]\!]$ and for $j \in [\![1, M]\!]$, that:
   \begin{multline}\label{eq:bound-T24}
       \mathbb{E}^{P_j}\left[T_n^{(2)}(\ell, j)\right] \leq -\dfrac{1}{2}\mathbb{E}^{P_j}\left[\sum_{k=1}^{n-1}{\mathcal{G}_m(X_{k\Delta}^{\ell})}\right] + \dfrac{1}{4}\sum_{k=1}^{n-1}{\mathbb{E}^{P_j}\left[\mathcal{G}_m^2(X_{k\Delta}^{\ell})\right]} + \mathrm{O}\left(\sum_{k=1}^{n-1}{\mathbb{E}^{P_j}\left[\mathcal{G}_m^2(X_{k\Delta}^{\ell})\right]}\right)\\
       + \dfrac{1}{2}\mathbb{E}^{P_j}\left[\sum_{k=1}^{n-1}{\mathcal{G}_m(X_{k\Delta}^{\ell})}\right] + \dfrac{1}{2\Delta}\sum_{k=1}^{n-1}{\int_{k\Delta}^{(k+1)\Delta}{\mathbb{E}^{P_j}\left[\sigma_j^2(X_s^{\ell}) - \sigma_j^2(X_{k\Delta}^{\ell})\right]}} + Ch^{\beta-2}\\
       = \dfrac{1}{4}\sum_{k=1}^{n-1}{\mathbb{E}^{P_j}\left[\mathcal{G}_m^2(X_{k\Delta}^{\ell})\right]} + \mathrm{O}\left(\sum_{k=1}^{n-1}{\mathbb{E}^{P_j}\left[\mathcal{G}_m^2(X_{k\Delta}^{\ell})\right]}\right)\\
       + \dfrac{1}{2\Delta}\sum_{k=1}^{n-1}{\int_{k\Delta}^{(k+1)\Delta}{\mathbb{E}^{P_j}\left[\sigma_j^2(X_s^{\ell}) - \sigma_j^2(X_{k\Delta}^{\ell})\right]}} + Ch^{\beta-2}. 
   \end{multline}
   For all $\ell \in \{1, \ldots, N\}, ~ j \in \{1, \ldots, M\}$ and $k \in \{1, \ldots, n\}$, we have 
   \begin{align*}
    \mathcal{G}_m^2(X^{\ell}_{k\Delta}) = &~ \Gamma^2\left(\sum_{r=1}^{m}{ w_r^j\eta_r(X_{k\Delta}^{\ell})}\right)^2 = \Gamma^2\sum_{r=1}^{m}{ w_r^j\eta_r^2(X_{k\Delta}^{\ell})}\\
    = &~ \Gamma^2R^2h^{2\beta}\sum_{r=1}^{m}{w_r^jK^2\left(\dfrac{X_{k\Delta}^{\ell} - A - (B-A)x_r}{(B-A)h}\right)\mathds{1}_{\left(\frac{r-1}{m}(B-A)+A, \frac{r}{m}(B-A)+A\right)}(X_{k\Delta}^{\ell})}\\
    \leq &~ \Gamma^2R^2h^{2\beta}\|K\|_{\infty}^2.
   \end{align*}
   We deduce that $\mathbb{E}^{P_j}\left[\mathcal{G}_m^2(X^{\ell}_{k\Delta})\right] \leq \Gamma^2R^2h^{2\beta}\|K\|_{\infty}^2$ and we obtain from Equation~\eqref{eq:bound-T24} that for $n$ large enough, there exists a constant $C>0$ such that
   \begin{equation}\label{eq:bound-T25}
       \mathbb{E}^{P_j}\left[T^{(2)}_{n}(\ell, j)\right] \leq C\Gamma^{2}R^{2}\|K\|^{2}_{\infty}h^{2\beta}n + Ch^{\beta-2}
       + \dfrac{1}{2\Delta}\sum_{k=1}^{n-1}{\int_{k\Delta}^{(k+1)\Delta}{\mathbb{E}^{P_j}\left[\sigma_j^2(X_s^{\ell}) - \sigma_j^2(X_{k\Delta}^{\ell})\right]}}.
   \end{equation}
   Finally, using the It\^o formula, we have the following:
   \begin{equation*}
       d\sigma_j^2(X_t^{\ell}) = (\sigma_j^2)^{\prime}(X_t^{\ell})\sigma_j(X_t^{\ell})dW_t^{\ell} + \dfrac{1}{2}(\sigma_j^2)^{\prime\prime}(X_t^{\ell})\sigma_j^2(X_t^{\ell})dt,
   \end{equation*}
   and for all $s \in [k\Delta, (k+1)\Delta]$, there exists a constant $C>0$ such that
   \begin{multline*}
       \mathbb{E}^{P_j}\left[\sigma_j^2(X_s^{\ell}) - \sigma_j^2(X_{k\Delta}^{\ell})\right] = \dfrac{1}{2}\int_{k\Delta}^{s}{(\sigma_j^2)^{\prime\prime}(X_t^{\ell})\sigma_j^2(X_t^{\ell})dt} \leq \left\|(\sigma_j^2)^{\prime\prime}\right\|_{\infty}\left\|\sigma_j^2\right\|_{\infty}(s-k\Delta) \leq C\Delta h^{\beta-2},
   \end{multline*}
   and from Equation~\eqref{eq:bound-T25} and for all $\ell \in \{1, \ldots, N\}$ and $j \in \{1, \ldots, M\}$, there exists a constant $C>0$ such that
   \begin{equation*}
       \mathbb{E}^{P_j}\left[T^{(2)}_{n}(\ell, j)\right] \leq C\left(\Gamma^{2}R^{2}\|K\|^{2}_{\infty}h^{2\beta}n + h^{\beta-2}\right).
   \end{equation*}
\end{proof}

\subsection{Proof of Theorem~\ref{thm:lower-bound-R-Nn}}

\begin{proof}
    For the sake of simplicity, the drift function is set to $b=0$. The goal is to prove the existence of a constant $c>0$ such that
    \begin{equation}\label{eq:lowerbound1}
        \forall ~ \widehat{\sigma}^2: \underset{\sigma^2 \in \Sigma}{\sup}{~\mathbb{E}\left[\left\|\widehat{\sigma}^2 - \sigma^2\right\|_n^2\right] \geq c\psi^{2}_{n,N}},
    \end{equation}
    where $\Sigma = \Sigma(\beta, R)$ is the H\"older space, and $\psi_{n,N} = (Nn)^{-\beta/(2\beta+1)}$. Equation~\eqref{eq:lowerbound1} is equivalent to 
    \begin{equation*}
        \forall \widehat{\sigma}^2, ~ \underset{\sigma^2 \in \Sigma}{\sup}{\mathbb{E}\left[\left\|\psi^{-1}_{n,N}(\widehat{\sigma}^2 - \sigma^2)\right\|_n^2\right]} \geq c.
    \end{equation*}
    Using a line of reasoning similar to the one used in the proof of Theorem~\ref{thm:lower-bound-compact}, we obtain the following result.
    \begin{align*}
        \forall ~ \widehat{\sigma}^2, ~ \underset{\sigma^2 \in \Sigma}{\sup}{\mathbb{E}\left[\left\|\psi^{-1}_{n,N}(\widehat{\sigma}^2 - \sigma^2)\right\|_n^2\right]} \geq \Lambda_0^2\underset{\sigma^2 \Sigma^{M}}{\sup}{~\mathbb{P}\left(\left\|\widehat{\sigma}^2 - \sigma^2\right\|_n \geq s_{n,N}\right)},
    \end{align*}
    where $\Lambda_0 > 0$ is a constant to be determined later, $s_{n,N} = \Lambda_0\psi_{n,N}$, and $\Sigma^M = \left\{\sigma_0^2, \ldots, \sigma_M^2\right\} \in \Sigma$ a finite set of $M+1$ hypotheses with $M \geq 2$. It suffices to prove the existence of a constant $p_M>0$ such that
    \begin{equation}\label{eq:lowerboundé-pM}
        \forall \sigma^2 \in \Sigma^{M}, ~ \forall \widehat{\sigma}^2: ~ \mathbb{P}\left(\left\|\widehat{\sigma}^2 - \sigma^2\right\|_n\right) \geq p_M.
    \end{equation}  
    For this purpose, we use the set of hypotheses $\Sigma^M \subset \Sigma$ defined in the proof of Theorem~\ref{thm:lower-bound-compact} with $B = -A = 1$, $\Gamma = \dfrac{\kappa_1^2 - 1}{R\|K\|_{\infty}}$, and $m = \left\lceil c_0(Nn)^{1/(2\beta+1)}\right\rceil$ where the constant $c_0>0$ should be chosen appropriately. Note that any function $\sigma_j^2 \in \Sigma^M$ takes values in the interval $[\kappa_0^2, \kappa_1^2]$ and satisfies
    \begin{align*}
        \sigma_j^2(x) = \begin{cases} 1 + \Gamma\sum_{k=1}^{m}{w_k^j\phi_k(x)} ~~ \mathrm{if} ~~ x \in (-1,1) \\\\ 1 ~~ \mathrm{if} ~~ x \notin (-1,1).
        \end{cases}
    \end{align*}
    Consequently, for all $j,\ell = \{0, \ldots, M\}$, we have the following:
    \begin{align*}
        \left\|\sigma_j^2 - \sigma_{\ell}^2\right\|_n^2 \geq \int_{-\infty}^{+\infty}{\left(\sigma_j^2 - \sigma_{\ell}^2\right)^2(x)f_n(x)dx} = \int_{-1}^{1}{\left(\sigma_j^2 - \sigma_{\ell}^2\right)^2(x)f_n(x)dx},
    \end{align*} 
    where $f_n: x \mapsto \dfrac{1}{n}\sum_{k=1}^{n-1}{p_X(0, k\Delta, 0,  x)}$ is a function that is lower bounded from below in the compact interval $[-1,1]$. Then, there exists a constant $c>0$ such that
    \begin{multline*}
        \left\|\sigma_j^2 - \sigma_{\ell}^2\right\|_n^2 \geq c\int_{-\infty}^{+\infty}{\left(\sigma_j^2 - \sigma_{\ell}^2\right)^2(x)dx} \\
        = c\Gamma^2\sum_{k=1}^{m}{(w_k^j-w_k^{\ell})\int_{-1}^{1}{\phi_k^2\left(\dfrac{x+1}{2}\right)dx}}\\
        \geq 2c\Gamma^2R^2h^{2\beta}\sum_{k=1}^{m}{\mathds{1}_{w_k^j \neq w_k^{\ell}}\int_{(k-1)/m}^{k/m}{K^2\left(\dfrac{x-x_k}{h}\right)}dx},
    \end{multline*}  
    where $h = 1/m$, for all $x\in (0,1)$ and for all $k \in \{1, \ldots, m\}$, $\phi_k(x) = 2Rh^{\beta}K\left(\dfrac{x-x_k}{h}\right)$, and $x_k = (k-1/2)/m$. We finally obtain:
    \begin{align*}
        \left\|\sigma_j^2 - \sigma_{\ell}^2\right\|_n^2 \geq 2c\Gamma^2R^2\|K\|^2h^{2\beta+1}\rho(w^j,w^{\ell}),
    \end{align*} 
    where $\rho(w^j,w^{\ell}) = \sum_{k=1}^{m}{\mathds{1}_{w_k^j \neq w_k^{\ell}}}$ is the Hammer distance. We consider sequences $w \in \{0,1\}^m$ associated with the set $\Sigma^M$ that satisfy Equation~\eqref{eq:LowerBound-HammingDistance}. We deduce that $M \geq 2^{m/8}$ and
    \begin{equation}\label{eq:cond1}
        \left\|\sigma_j^2 - \sigma_{\ell}^2\right\|_n \geq \dfrac{1}{2}c^{1/2}R\|K\|h^{\beta} \geq 2s_{n,N},
    \end{equation}
    where $s_{n,N} = \Lambda_0\psi_{n,N}$ and 
    \begin{equation*}
        \Lambda_0 = \dfrac{c^{1/2}\Gamma R\|K\|}{2^{\beta+1}c_0^{\beta}}.
    \end{equation*}  
    From the proof of Theorem~\ref{thm:lower-bound-compact} with $B = -A = 1$ and the constant $c_0 > 0$ chosen so that 
    $$C\dfrac{\Gamma^2R^2\|K\|_{\infty}^2}{c_0^{2\beta+1}\log(2)}= \dfrac{1}{16} \times \dfrac{\log(2)}{8},$$ 
    we have for $n \propto N^2$ and $\beta \geq 4$,
    \begin{equation}\label{eq:cond2}
         \dfrac{1}{M}\sum_{j=1}^{M}{K_{\mathrm{div}}(P_j,P_0)} = \frac{1}{M}\sum_{j=1}^{M}{\sum_{\ell=1}^{N}{\mathrm{E}^{P_j}\left[\log \frac{dP_j}{dP_0}(\bar{X}^{\ell})\right]}} \leq \dfrac{1}{16}\log(M), ~~ \dfrac{1}{16} \in (0,1/8),
    \end{equation} 
    where $K_{\mathrm{div}}$ is the Kullback divergence between $P_j$ and $P_0$, the distributions of the diffusion processes whose diffusion coefficients are, respectively, $\sigma_j^2$ and $\sigma_0^2$. Then, from Equations~\eqref{eq:cond1} and \eqref{eq:cond2} and since $\Sigma^M \subset \Sigma$, we obtain from Theorem 2.5 in \cite{tsybakov2008introduction} that
    \begin{align*}
      & \underset{\w{\sigma}^{2}}{\inf}{\underset{\sigma^2 \in \Sigma^{M}}{\sup}{~\P\left(\|\w{\sigma}^{2}-\sigma^{2}\|_n \geq s_{n,N}\right)}} \geq p_M \geq \frac{1}{2}\left(\dfrac{7}{8}-\dfrac{1}{2\sqrt{2\log(2)}}\right) > 0, \\
      & \underset{\w{\sigma}^{2}}{\inf}{\underset{\sigma^2 \in \Sigma^{M}}{\sup}{~\E\left[\|\w{\sigma}^{2}-\sigma^{2}\|_n^{2}\right]}} \geq c(Nn)^{-2\beta/(2\beta+1)}.
   \end{align*}
   Moreover, for $N$ large enough such that $[-1, 1] \subset [-A_N, A_N]$, we also obtain
   \begin{equation*}
       \underset{\w{\sigma}^{2}}{\inf}{\underset{\sigma^2 \in \Sigma^{M}}{\sup}{~\E\left[\|\w{\sigma}^{2}-\sigma^{2}_{A_N}\|_n^{2}\right]}} \geq c(Nn)^{-2\beta/(2\beta+1)},
   \end{equation*}
   where $\sigma_{A_N}^2 = \sigma^{2}\mathds{1}_{[-A_N, A_N]}$.
\end{proof}


\bibliographystyle{ScandJStat}
\bibliography{mabiblio.bib}

\begin{thebibliography}{33}
\providecommand{\natexlab}[1]{#1}
\providecommand{\url}[1]{\texttt{#1}}
\providecommand{\urlprefix}{URL }
\expandafter\ifx\csname urlstyle\endcsname\relax
  \providecommand{\doi}[1]{doi:\discretionary{}{}{}#1}\else
  \providecommand{\doi}{doi:\discretionary{}{}{}\begingroup
  \urlstyle{rm}\Url}\fi

\bibitem[{Barron \emph{et~al.}(1999)Barron, Lucien \& Pascal}]{barron1999risk}
Barron, A., Lucien, B. \& Pascal, M. (1999).
\newblock Risk bounds for model selection via penalization.
\newblock \emph{Probability theory and related fields} \textbf{113}, 301--413.

\bibitem[{Bressloff(2024)}]{bressloff2024cellular}
Bressloff, P.~C. (2024).
\newblock Cellular diffusion processes in singularly perturbed domains.
\newblock \emph{Journal of Mathematical Biology} \textbf{89}, 58.

\bibitem[{Cairns(2004)}]{cairns2004interest}
Cairns, A.~J. (2004).
\newblock \emph{Interest rate models: an introduction}, vol.~10.
\newblock Princeton University Press.

\bibitem[{Clement(1997)}]{clement1997}
Clement, E. (1997).
\newblock Estimation of diffusion processes by simulated moment methods.
\newblock \emph{Scandinavian Journal of Statistics} \textbf{24}, 353--369.

\bibitem[{Comte(2020)}]{comte2019regression}
Comte, F. (2020).
\newblock From regression function to diffusion drift estimation in
  nonparametric setting.
\newblock \emph{ESAIM: Proceedings and Surveys} \textbf{68}, 20--34.

\bibitem[{Comte \& Genon-Catalot(2020)}]{comte2020nonparametric}
Comte, F. \& Genon-Catalot, V. (2020).
\newblock Nonparametric drift estimation for iid paths of stochastic
  differential equations.
\newblock \emph{The Annals of Statistics} \textbf{48}, 3336--3365.

\bibitem[{Comte \& Genon-Catalot(2021)}]{comte2021drift}
Comte, F. \& Genon-Catalot, V. (2021).
\newblock Drift estimation on non compact support for diffusion models.
\newblock \emph{Stochastic Processes and their Applications} \textbf{134},
  174--207.

\bibitem[{Comte \emph{et~al.}(2007)Comte, Genon-Catalot, Rozenholc
  \emph{et~al.}}]{comte2007penalized}
Comte, F., Genon-Catalot, V., Rozenholc, Y. \emph{et~al.} (2007).
\newblock Penalized nonparametric mean square estimation of the coefficients of
  diffusion processes.
\newblock \emph{Bernoulli} .

\bibitem[{Crow(2017)}]{crow2017introduction}
Crow, J.~F. (2017).
\newblock \emph{An introduction to population genetics theory}.
\newblock Scientific Publishers.

\bibitem[{Dacunha-Castelle \& Florens-Zmirou(1986)}]{dacunha1986estimation}
Dacunha-Castelle, D. \& Florens-Zmirou, D. (1986).
\newblock Estimation of the coefficients of a diffusion from discrete
  observations.
\newblock \emph{Stochastics: An International Journal of Probability and
  Stochastic Processes} \textbf{19}, 263--284.

\bibitem[{De~Boor(1978)}]{deboor78}
De~Boor, C. (1978).
\newblock \emph{A practical guide to splines}, vol.~27.
\newblock springer-verlag New York.

\bibitem[{Denis \emph{et~al.}(2024)Denis, Dion-Blanc, Ella-Mintsa \&
  Tran}]{denis2024nonparametric}
Denis, C., Dion-Blanc, C., Ella-Mintsa, E. \& Tran, V.~C. (2024).
\newblock Nonparametric plug-in classifier for multiclass classification of sde
  paths.
\newblock \emph{Scandinavian Journal of Statistics} .

\bibitem[{Denis \emph{et~al.}(2020)Denis, Dion-Blanc \&
  Martinez}]{denis2020consistent}
Denis, C., Dion-Blanc, C. \& Martinez, M. (2020).
\newblock Consistent procedures for multiclass classification of discrete
  diffusion paths.
\newblock \emph{Scandinavian Journal of Statistics} \textbf{47}, 516--554.

\bibitem[{Denis \emph{et~al.}(2021)Denis, Dion-Blanc \&
  Martinez}]{denis2020ridge}
Denis, C., Dion-Blanc, C. \& Martinez, M. (2021).
\newblock A ridge estimator of the drift from discrete repeated observations of
  the solutions of a stochastic differential equation.
\newblock \emph{Bernoulli} .

\bibitem[{El~Karoui \emph{et~al.}(1997)El~Karoui, Peng \&
  Quenez}]{el1997backward}
El~Karoui, N., Peng, S. \& Quenez, M.~C. (1997).
\newblock Backward stochastic differential equations in finance.
\newblock \emph{Mathematical finance} \textbf{7}, 1--71.

\bibitem[{Ella-Mintsa(2024)}]{ella2024nonparametric}
Ella-Mintsa, E. (2024).
\newblock Nonparametric estimation of the diffusion coefficient from iid sde
  paths.
\newblock \emph{Statistical Inference for Stochastic Processes} pp. 1--56.

\bibitem[{Florens-Zmirou(1993)}]{florens1993estimating}
Florens-Zmirou, D. (1993).
\newblock On estimating the diffusion coefficient from discrete observations.
\newblock \emph{Journal of applied probability} \textbf{30}, 790--804.

\bibitem[{Gadat \emph{et~al.}(2020)Gadat, Gerchinovitz \&
  Marteau}]{gadat2020optimal}
Gadat, S., Gerchinovitz, S. \& Marteau, C. (2020).
\newblock Optimal functional supervised classification with separation
  condition.
\newblock \emph{Bernoulli} \textbf{26}, 1797--1831.

\bibitem[{Genon-Catalot \& Jacod(1993)}]{genon1993estimation}
Genon-Catalot, V. \& Jacod, J. (1993).
\newblock On the estimation of the diffusion coefficient for multi-dimensional
  diffusion processes.
\newblock In \emph{Annales de l'IHP Probabilit{\'e}s et statistiques}, vol.~29,
  pp. 119--151.

\bibitem[{Gloter(2000)}]{gloter2000discrete}
Gloter, A. (2000).
\newblock Discrete sampling of an integrated diffusion process and parameter
  estimation of the diffusion coefficient.
\newblock \emph{ESAIM: Probability and Statistics} \textbf{4}, 205--227.

\bibitem[{Gobet(2002)}]{gobet2002lan}
Gobet, E. (2002).
\newblock Lan property for ergodic diffusions with discrete observations.
\newblock \emph{Annales de l'Institut Henri Poincare (B) Probability and
  Statistics} \textbf{38}, 711--737.

\bibitem[{Gy{\"o}rfi \emph{et~al.}(2006)Gy{\"o}rfi, Kohler, Krzyzak \&
  Walk}]{gyorfi2006distribution}
Gy{\"o}rfi, L., Kohler, M., Krzyzak, A. \& Walk, H. (2006).
\newblock \emph{A distribution-free theory of nonparametric regression}.
\newblock Springer Science \& Business Media.

\bibitem[{Hoffmann(1997)}]{hoffmann1997minimax}
Hoffmann, M. (1997).
\newblock Minimax estimation of the diffusion coefficient through irregular
  samplings.
\newblock \emph{Statistics \& probability letters} \textbf{32}, 11--24.

\bibitem[{Hoffmann(1999)}]{hoffmann1999lp}
Hoffmann, M. (1999).
\newblock Lp estimation of the diffusion coefficient.
\newblock \emph{Bernoulli} pp. 447--481.

\bibitem[{Jacod(1993)}]{jacod1993random}
Jacod, J. (1993).
\newblock Random sampling in estimation problems for continuous gaussian
  processes with independent increments.
\newblock \emph{Stochastic processes and their applications} \textbf{44},
  181--204.

\bibitem[{Lamberton \& Lapeyre(2011)}]{lamberton2011introduction}
Lamberton, D. \& Lapeyre, B. (2011).
\newblock \emph{Introduction to stochastic calculus applied to finance}.
\newblock Chapman and Hall/CRC.

\bibitem[{Massart \emph{et~al.}(2006)Massart, N{\'e}d{\'e}lec
  \emph{et~al.}}]{massart2006risk}
Massart, P., N{\'e}d{\'e}lec, E. \emph{et~al.} (2006).
\newblock Risk bounds for statistical learning.
\newblock \emph{The Annals of Statistics} \textbf{34}, 2326--2366.

\bibitem[{Nagai \& Mimura(1983)}]{nagai1983asymptotic}
Nagai, T. \& Mimura, M. (1983).
\newblock Asymptotic behavior for a nonlinear degenerate diffusion equation in
  population dynamics.
\newblock \emph{SIAM Journal on Applied Mathematics} \textbf{43}, 449--464.

\bibitem[{Romanczuk \emph{et~al.}(2012)Romanczuk, B{\"a}r, Ebeling, Lindner \&
  Schimansky-Geier}]{romanczuk2012active}
Romanczuk, P., B{\"a}r, M., Ebeling, W., Lindner, B. \& Schimansky-Geier, L.
  (2012).
\newblock Active brownian particles: From individual to collective stochastic
  dynamics.
\newblock \emph{The European Physical Journal Special Topics} \textbf{202},
  1--162.

\bibitem[{Sbalzarini(2006)}]{sbalzarini2006analysis}
Sbalzarini, I.~F. (2006).
\newblock Analysis, modeling, and simulation of diffusion processes in cell
  biology.
\newblock Ph.D. thesis, ETH Zurich.

\bibitem[{S{\o}rensen(2002)}]{sorensen2002estimation}
S{\o}rensen, H. (2002).
\newblock Estimation of diffusion parameters for discretely observed diffusion
  processes.
\newblock \emph{Bernoulli} pp. 491--508.

\bibitem[{Tsybakov(2008)}]{tsybakov2008introduction}
Tsybakov, A.-B. (2008).
\newblock \emph{Introduction to nonparametric estimation}.
\newblock Springer Science \& Business Media.

\bibitem[{Wu(2019)}]{wu2019interest}
Wu, L. (2019).
\newblock \emph{Interest rate modeling: Theory and practice}.
\newblock CRC Press.

\end{thebibliography}

\end{document}